\documentclass[12pt,fleqn]{article}

\usepackage{amsmath}
\usepackage{amsthm}
\usepackage{amsfonts,dsfont}
\usepackage{amssymb}
\usepackage[titletoc]{appendix}
\usepackage{color}
\usepackage{comment}
\usepackage{enumerate}
\usepackage{fancyhdr} 
\usepackage{geometry}
\usepackage[pdfauthor={Alexander Bihlo, Francis Valiquette},%
pdftitle={Symmetry-preserving numerical schemes,pdftex}]{hyperref}%
\usepackage{stmaryrd}
\usepackage{tikz}
\usetikzlibrary{shapes.misc}
\tikzset{cross/.style={cross out, draw=black, minimum size=2*(#1-\pgflinewidth), inner sep=0pt, outer sep=0pt},
cross/.default={5pt}}
\usetikzlibrary{shapes,arrows}
\usepackage{verbatim}
\usepackage[all,cmtip]{xy}
\usepackage{subcaption}
\usepackage{caption}
\hypersetup{colorlinks, linkcolor=black, citecolor=black, urlcolor=blue}

\theoremstyle{definition}\newtheorem{definition}{Definition}[section]
\theoremstyle{definition}\newtheorem{proposition}[definition]{Proposition}
\theoremstyle{definition}\newtheorem{remark}[definition]{Remark}
\theoremstyle{definition}
\theoremstyle{definition}
\theoremstyle{definition}
\theoremstyle{definition}\newtheorem{example}[definition]{Example}
\theoremstyle{definition}\newtheorem{exercise}[definition]{Exercise}
\theoremstyle{definition}

\newcommand{\pp}[2]{\frac{\partial #1}{\partial #2}}
\newcommand{\vv}{\mathbf{v}}

\newcommand{\g}{\mathfrak{g}}

\newcommand{\J}{{\rm J}}

\newcommand{\sn}{^{[\ell]}}
\newcommand{\n}{^{(\ell)}}
\newcommand{\on}{\ell}

\newcommand{\mJ}{\mathcal{J}}
\newcommand{\mK}{\mathcal{K}}

\newcommand{\mM}{\mathcal{M}}
\newcommand{\mO}{\mathcal{O}}

\newcommand{\mS}{\mathcal{S}}

\def\comp{\raise 1pt \hbox{$\,\scriptstyle\circ\,$}}

\numberwithin{equation}{section}

\begin{document}


\thispagestyle{fancy}
\fancyhead{}
\fancyfoot{}
\renewcommand{\headrulewidth}{0pt}
\cfoot{\thepage}
\rfoot{\today}

\vskip 1cm
\begin{center}
{\Large Symmetry-Preserving Numerical Schemes}
\vskip 1cm

\begin{tabular*}{1.0\textwidth}{@{\extracolsep{\fill}} ll}
Alexander Bihlo
& Francis Valiquette\\
Department of Mathematics and Statistics & Department of Mathematics\\
Memorial University & SUNY New Paltz \\
St. John's, NL, Canada\quad A1C 5S7& New Paltz, NY, USA \quad 12561\\
{\tt abihlo@mun.ca} & {\tt valiquef@newpaltz.edu} \\
{\tt http://www.math.mun.ca/$\sim$abihlo} &{\tt http://www2.newpaltz.edu/$\sim$valiquef}
\end{tabular*}
\end{center}


\vskip 1cm

\abstract{}In these lectures we review two procedures for constructing finite difference numerical schemes that preserve symmetries of differential equations.  The first approach is based on Lie's infinitesimal symmetry generators, while the second method uses the novel theory of equivariant moving frames.  The advantages of both techniques are discussed and illustrated with the Schwarzian differential equation, the Korteweg--de Vries equation and Burgers' equation.  Numerical simulations are presented and innovative techniques for obtaining better invariant numerical schemes are introduced.  New research directions and open problems are indicated at the end of these notes.

\tableofcontents

\section{Introduction}

The aim of geometric numerical integration is to develop numerical integrators that preserve geometric properties of the system of differential equations under investigation.  Classical examples include symplectic integrators, \cite{HWL-2006,LR-2004}, energy preserving methods, \cite{QM-2008}, and schemes that preserve a Lie--Poisson structure, \cite{ZM-1988}.  The motivation behind geometric numerical integration is that, as a rule of thumb, such integrators will typically give better global or long term numerical results than standard methods since they incorporate qualitative properties of the system under consideration.  

In mathematical physics, most fundamental differential equations are invariant under a certain collection of symmetry transformations.  These symmetries can be point transformations, contact transformations, or generalized transformations, \cite{O-1993}.  In all cases, the symmetries of a differential equation encapsulate important properties of the equation and its solutions.  Furthermore, Lie group techniques are amongst the most effective methods for obtaining explicit solutions and conservation laws of nonlinear differential equations, \cite{BA-2002,O-1993,O-1982}. 

When discretizing differential equations invariant under a certain symmetry group, there are different incentives for preserving the symmetries of these equations.  From a physical standpoint, discrete spacetime models should preserve the symmetries of their continuous counterparts.  Mathematically, Lie group techniques could then be used to find explicit solutions and compute conservation laws of the discrete models.  From a more practical point of view, symmetry-preserving discretizations should share some exact solutions with the original differential equations, or at least provide better approximations than non-invariant numerical schemes.  

In the last 30 years, the application of Lie group techniques to finite difference equations has become a very active field of research.  To the best of our knowledge, Yanenko and Shokin were the first to use group theoretical methods to study finite difference schemes by introducing first differential approximations of difference equations, \cite{S-1983, YS-1973}.  The application of Lie group methods to finite difference equations, as we know it today, was first introduced by Dorodnitsyn in 1989, \cite{D-1989}.  Early on, one of the main focuses in the field was to construct Lie point symmetry-preserving finite difference approximations of differential equations. Beside Dorodnitsyn, early contributors include Bakirova, Kozlov, Levi, and Winternitz who constructed symmetry-preserving schemes for heat transfer equations, \cite{BD-1994,BDK-1997,DK-2003}, variable coefficient Korteweg--de Vries equations, \cite{DW-2000}, Burgers' equation, \cite{HLW-2000}, the nonlinear Schr\"odinger equation, \cite{BD-2001}, and second-order ordinary differential equations, \cite{DKW-2000}.  Symmetry-preserving approximation of Euler--Lagrange equations and their corresponding Lagrangian have also been considered in \cite{DKW-2003,DKW-2004}, and the application of Noether's theorem to compute  conservation laws has been extensively studied in the discrete setting, \cite{D-2001,H-2014,HM-2011}.  The applications of Lie point symmetries to finite difference equations have also been extended to generalized symmetries, \cite{LY-2001, LY-2009}, $\lambda$-symmetries, \cite{LNR-2012,LR-2010}, and contact transformations, \cite{LSTW-2012}.   

In recent years, more systematic efforts have been directed towards investigating the numerical performance of symmetry-preserving schemes.  For ordinary differential equations, symmetry-preserving schemes have proven to be very promising.  For solutions exhibiting sharp variations or singularities, symmetry-preserving schemes systematically appear to outperform standard numerical schemes, \cite{BCW-2006, BRW-2008, CRW-2016,KO-2004}.  For partial differential equations, the improvement of symmetry-preserving schemes versus traditional integrators is not as clear, \cite{BCW-2015,K-2008,LMW-2015,RV-2013}. On one hand, it was shown in \cite{WKO-2007} that symmetry-preserving schemes do much better in tracking sharp interfaces in Hamilton--Jacobi equations.  On the other hand, invariant numerical schemes for evolution equations generally require the use of time-evolving meshes which can lead to mesh tangling and thereby severely limit the use of symmetry-preserving schemes.  In this case, special techniques have to be developed to avoid mesh singularities.  For example, new ideas relying on r-adaptivity have been implemented to improve the performance of invariant integrators, \cite{BP-2012}. Also, in \cite{BN-2013,BN-2014} an invariant evolution--projection strategy was introduced and invariant meshless discretization schemes were considered in \cite{B-2013}.

The preceding references only provide a short bibliographical overview of the field. Many papers had to be omitted.  More references on the subject can be found in the survey papers \cite{LW-2006,W-2004}, and the books \cite{D-2010,H-2014}.

Given a differential equation with symmetry group $G$, the first step in constructing a symmetry-preserving numerical scheme is to compute difference invariants of the product action of $G$ on a chosen stencil.  There are mainly two approaches for constructing those invariants.  Most of the references cited above use the infinitesimal symmetry generators of the group action and Lie's infinitesimal invariance criterion to construct difference invariants.  Alternatively, the difference invariants can be constructed using the novel method of equivariant moving frames mainly developed by Olver, which was done in \cite{B-2013,CH-2010,KO-2004,O-2001,RV-2013,RV-2015}.  Given sufficiently many difference invariants, an invariant numerical scheme is, in general, obtained by finding a suitable combination of these invariants that converges to the original differential equation in the continuous limit.  When using Lie's infinitesimal generator approach, a suitable combination is found by taking the Taylor expansion of the difference invariants and combining them in such a way to obtained the desired invariant scheme.  With the method of moving frames, a suitable combination is found more systematically by invariantizing a  non-invariant numerical scheme.  Since the symmetry group of a differential equation will, in general, act on both the independent and dependent variables, a symmetry-preserving numerical scheme will usually not be defined on a uniform orthogonal mesh.

The application of Lie groups to finite difference equations is a vast and very dynamic field of study.  While preparing these lecture notes we had to omit many interesting applications and important results.  The focus of these notes will be on the theoretical construction of invariant numerical schemes and their numerical implementation.  At the heart of all our considerations are differential equations, finite difference equations, symmetry groups, and invariants.  
These familiar notions are all reviewed in Sections \ref{preliminary section}, \ref{symmetry section}, and \ref{invariants section}.  As outlined above, there are two different approaches for computing invariants of a Lie group action.   The infinitesimal approach based on Lie's symmetry generators is introduced in Section \ref{invariants - symmetry generators section}, while the equivariant moving frame approach is explained in Section \ref{invariants - moving frame section}. Section~\ref{weakly invariant equation section} is devoted to weakly invariant equations, which can play an important role in the construction of symmetry-preserving schemes. The construction of symmetry-preserving numerical schemes is carefully explained in Section \ref{symmetry-preserving schemes section}.  To illustrate the implementation, we consider the Schwarzian differential equation and the Korteweg--de Vries (KdV) equation. In Section \ref{numerical simulations section} we carry out numerical simulations for the Schwarzian equation, the KdV equation and Burgers' equation.  For partial differential equations, the invariance of a numerical scheme does not, in general, guarantee better numerical results. We will show that symmetry-preserving schemes can lead to mesh tangling, which limit their practical scope.  To circumvent this shortcoming, we discuss new invariant numerical strategies.  For the Korteweg--de Vries equation, we introduce \emph{invariant evolution--projection schemes} and \emph{invariant adaptive numerical schemes}.  Unlike the KdV equation, solutions to Burgers' equation can exhibit shocks.  For these shock solutions we propose a new \emph{invariant finite volume} type scheme.  Finally, in Section \ref{conclusion section} we identify some open problems and challenges in the field of symmetry-preserving numerical integrators.

\section{Differential and difference equations}\label{preliminary section}

In this section we review the definitions of differential equations and finite difference equations.  We take this opportunity to introduce some of the notation used throughout these notes.

\subsection{Differential equations}

Let $M$ be an $m$-dimensional manifold.  For $0 \leq \ell \leq \infty$, let $\J\n = \J\n(M,p)$ denote the \emph{extended $\on^\text{th}$ order jet space} of $1 \leq p < m$ dimensional submanifolds $S \subset M$ defined as the space of equivalence classes of submanifolds under the equivalence relation of $\on^\text{th}$ order contact at a point, \cite{O-1993}.  Local coordinates on $\J\n$ are given by the $\on$-jet 
\begin{equation}\label{standard jet coordinates}
(x,u\n),
\end{equation}
where $x=(x^1,\ldots,x^p)$ correspond to the independent variables and $u\n$ denotes the derivatives
\[
u^\alpha_{x^J} = \frac{\partial^k u^\alpha}{(\partial x^1)^{j_1} \cdots (\partial x^p)^{j_p}}\quad \text{with}\quad 1\leq \alpha \leq q = m-p\quad\text{and}\quad 0 \leq k = j_1+\cdots+j_p \leq \on.
\]
In the above notation, $J=(j_1,\ldots,j_p)$ is an ordered $p$-tuple of non-negative integers, with entries $j_i \geq 0$ indicating the number of derivatives taken in the variable $x^i$.  The order of the multi-index, denoted by $\# J = k$, indicates how many derivatives are being taken.  

\begin{example}\label{graph jet example}
In the case where $p=2$ and $q = 1$, we have two independent variables $(x^1,x^2) = (t,x)$ and one dependent variable $u^1=u$.  Then, the second order jet space is parametrized by
\begin{equation*}
(t,x,u,u_t,u_x,u_{tt},u_{tx},u_{xx}).
\end{equation*}
\end{example}

\begin{definition}
A \emph{differential equation} of order $n$ is the zero locus of a differential map $\Delta\colon \J\n \to \mathbb{R}$. That is,
\begin{equation}\label{differential equation}
\Delta(x,u\n) = 0.
\end{equation}
\end{definition}

For later use, we introduce two regularity requirements on differential equations.

\begin{definition}
A differential equation $\Delta(x,u\n)=0$ is said to be \emph{regular} if the rank of its differential 
\begin{equation*}
d\Delta = \sum_{i=1}^p \pp{\Delta}{x^i} dx^i + \sum_{J} \sum_{\alpha=1}^q \pp{\Delta}{u^\alpha_{x^J}} du^\alpha_{x^J}
\end{equation*}
is constant on the domain of definition of $\Delta\colon \J\n \to \mathbb{R}$.
\end{definition}

\begin{example}
Any evolutionary partial differential equation
\[
\Delta(t,x,u\n) = u_t - f(t,x,u,u_x,u_{xx},\ldots,u_{x^\on}) = 0
\]
is regular since the rank of $d\Delta = du_t - df$ is one.
\end{example}

\begin{definition}\label{regular definition differential equations}
A differential equation $\Delta(x,u\n)=0$ is \emph{locally solvable} at a point $(x_0,u_0\n)$ if there exists a smooth solution $u=f(x)$, defined in the neighborhood of $x_0$, such that $u_0\n=f\n(x_0)$.  A differential equation which is both regular and locally solvable is said to be \emph{fully regular}.
\end{definition}

The above description assumes that a submanifold $S \subset M$ is locally represented as the graph of a function $S = \{(x,f(x))\}$.  Alternatively, when no distinction between independent and dependent variables is made, a submanifold $S \subset M$ is locally parameterized by $p$ variables $s=(s^1,\ldots,s^p)\in \mathbb{R}^p$ such that
\[
(x(s),u(s))\in S.
\]
In numerical analysis, the independent variables $s=(s^1,\ldots,s^p)$ are called \emph{computational variables}, \cite{HR-2011}.  These are the variables that are discretized when finite difference equations are introduced in Section \ref{section:difference equations}.  We let $\mathcal{J}\n$ denote the $\ell^\text{th}$ order jet space of submanifolds $S\subset M$ parametrized by computational variables.  Local coordinates on $\mathcal{J}\n$ are given by
\begin{equation}\label{computational jet coordinates}
(s,x\n,u\n)=(\, \ldots\, s^i\, \ldots\, x^i_{s^A}\, \ldots\, u^\alpha_{s^A}\, \ldots\,),
\end{equation}
with $1\leq i \leq p$, $1\leq \alpha \leq q$, and $0 \leq \# A \leq \on$.  

\begin{remark}
We hope that the jet notations $(x,u\n)$ and $(s,x\n,u\n)$ will not confuse the reader.  The independent variable, that is $x$ and $s$, respectively, indicates with respect to which variables the dependent variables $u$ (and $x$ in the second case) are differentiated in $u\n$.
\end{remark}

\begin{example}\label{computation variables jet example}
In the case where $p=2$ and $m=3$, let $(s^1,s^2) = (\tau, s)$ denote the two computational variables and let $(t,x,u)$ be a local parametrization of $M$.  Then the second order jet space $\mJ^{(2)}$ is parametrized by
\[
(\tau,s,t,x,u,t_\tau,t_s,x_\tau,x_s,u_\tau,u_s,t_{\tau\tau},t_{\tau s}, t_{ss},x_{\tau\tau},x_{\tau s}, x_{ss},u_{\tau\tau},u_{\tau s}, u_{ss}).
\]
\end{example}

The transition between the jet coordinates \eqref{standard jet coordinates} and \eqref{computational jet coordinates} is given by the chain rule.  Provided
\begin{equation}\label{change of variables condition}
\text{det}\bigg( \pp{x^j}{s^i} \bigg) \neq 0,\qquad \text{where}\qquad 1\leq i,j\leq p,
\end{equation}
the implicit total derivative operators
\[
D_{x^i} = \sum_{j=1}^p\> W^j_i\, D_{s^j},\qquad \big(W^j_i\big)=\bigg( \pp{x^j}{s^i} \bigg)^{-1},\qquad i=1,\ldots,p,
\]
are well-defined, and successive application of those operators to the dependent variables $u^\alpha$ will give the coordinate expressions for the $x$-derivatives of $u$ in terms of the $s$-derivatives of $x$ and $u$:
\begin{equation}\label{uxJ}
u^\alpha_{x^J} = (D_{x^1})^{j_1}\cdots (D_{x^p})^{j_p} u^\alpha = \bigg(\sum_{l=1}^p\> W^l_1 \, D_{s^l}\bigg)^{j_1} \cdots \bigg(\sum_{l=1}^p\> W^l_p\, D_{s^l}\bigg)^{j_p} u^\alpha.
\end{equation}
We note that the non-degeneracy constraint \eqref{change of variables condition} implies that the change of variables $x=x(s)$ is invertible.

\begin{example}\label{change of variable example}
Combining Examples \ref{graph jet example} and \ref{computation variables jet example}, assume that $x=x(\tau,s)$ and $t=t(\tau,s)$ are functions of the computational variables $(\tau,s)$.  Provided
\begin{equation}\label{non-degeneracy condition}
\det \begin{bmatrix} t_\tau & t_s \\ x_\tau & x_s \end{bmatrix}  = t_\tau x_s - t_s x_\tau \neq 0,
\end{equation}
the implicit derivative operators 
\begin{equation}\label{D_x,D_t}
D_x = \frac{t_\tau\, D_s - t_{s}\, D_\tau}{x_s t_\tau - x_\tau t_s},\qquad
D_t = \frac{x_s\, D_\tau - x_\tau\, D_s}{x_s t_\tau - x_\tau t_s}.
\end{equation}
are well-defined.  It follows that
\begin{equation}\label{ux-ut}
u_x = \frac{t_\tau\, u_s - t_{s}\, u_\tau}{x_s t_\tau - x_\tau t_s},\qquad
u_t = \frac{x_s\, u_\tau - x_\tau\, u_s}{x_s t_\tau - x_\tau t_s}.
\end{equation}
Relations for the higher order derivatives are obtained by applying \eqref{D_x,D_t} to \eqref{ux-ut}.
\end{example}

Given a differential equation \eqref{differential equation}, the chain rule \eqref{uxJ} can be used to re-express \eqref{differential equation} in terms of $x^i=x^i(s)$, $u^\alpha=u^\alpha(s)$ and their computational derivatives $x^i_{s^A}$, $u^\alpha_{s^A}$:

\begin{subequations}\label{extended system 1}
\begin{equation}\label{computational differential equation}
\overline{\Delta}(s,x\n,u\n)=\Delta(x,u\n)=0.
\end{equation}
Recall that $(s,x\n,u\n) = (s,\, \ldots\, x^i_{s^A}\, \ldots\, u^\alpha_{s^A}\, \ldots\,)\in \mJ\n$ for $\overline{\Delta}=0$ in \eqref{computational differential equation} while $(x,u\n) = (x,\, \ldots\, u^\alpha_{x^J}\, \ldots\,)\in \J\n$ in $\Delta=0$.  Equation \eqref{computational differential equation} can be supplemented by \emph{companion equations}, \cite{M-2010}, 
\begin{equation}\label{companion equations}
\widetilde{\Delta}(s,x\n,u\n)=0.
\end{equation}
\end{subequations}
The latter are introduced to impose restrictions on the change of variables $x=x(s)$.   The system of differential equations \eqref{extended system 1} is called an \emph{extended system} of the differential equation \eqref{differential equation}.  For the extended system of differential equations \eqref{extended system 1} to share the same solution space as the original equation \eqref{differential equation}, the companion equations \eqref{companion equations} cannot introduce differential constraints on the derivatives $u^\alpha_{s^A}$.

Definition \ref{regular definition differential equations} is readily adapted to the computational variable framework.

\begin{definition}\label{regular definition differential equations 2}
A differential equation $\overline{\Delta}(s,x\n,u\n)=0$ (or system of differential equations) is \emph{regular} if the rank of its differential 
\[
d\overline{\Delta} = \sum_{i=1}^p \pp{\overline{\Delta}}{s^i} ds^i + \sum_J \sum_{i=1}^p \pp{\overline{\Delta}}{x^i_{s^J}} dx^i_{s^J} + \sum_J \sum_{\alpha=1}^q \pp{\overline{\Delta}}{u^\alpha_{s^J}} du^\alpha_{s^J}
\]
is constant on the domain of definition.  The equation (or system of equations) is \emph{locally solvable} at a point $(s_0,x_0\n,u_0\n)$ if there exists a smooth solution $x=f(s)$, $u=g(s)$, defined in the neighborhood of $s_0$, such that $x_0\n=f\n(s_0)$ and $u_0\n=g\n(s_0)$.  The differential equation (or system of differential equations) is said to be \emph{fully regular} if it is both regular and locally solvable.
\end{definition}

\begin{example}\label{KdV in computational variables}
As one of our main examples in these notes, we consider the Korteweg--de Vries (KdV) equation
\begin{equation}\label{KdV}
u_t + u u_x + u_{xxx} = 0.
\end{equation}
We introduce the computational variables $(\tau,t)$ so that $x=x(\tau,s)$, $t=t(\tau,s)$.  Then the implicit total derivative operators are given by \eqref{D_x,D_t}.  Before proceeding any further, we assume that
\begin{equation}\label{KdV t companion equations}
t_s = 0,\qquad t_{\tau\tau}=0.
\end{equation}
In other words,
\begin{equation}\label{t solution}
t = k \tau + t^0,
\end{equation}
where $k\neq 0$ and $t^0$ are constants.  The reasons for imposing the constraints \eqref{KdV t companion equations} are explained in Example \ref{KdV symmetry example}.  The operators of implicit differentiation \eqref{D_x,D_t} then simplify to
\[
D_x = \frac{1}{x_s}D_s,\qquad
D_t = \frac{1}{t_\tau}\bigg(D_\tau - \frac{x_\tau}{x_s} D_s\bigg).
\]
Therefore,
\[
u_x = \frac{u_s}{x_s},\qquad u_{xx} = \frac{1}{x_s}\bigg(\frac{u_s}{x_s}\bigg)_s,\qquad
u_{xxx} = \frac{1}{x_s}\bigg(\frac{1}{x_s}\bigg(\frac{u_s}{x_s} \bigg)_s \bigg)_s,\qquad
u_t = \frac{u_\tau}{t_\tau} - \frac{x_\tau}{t_\tau}\cdot \frac{u_s}{x_s}
\]
and the KdV equation \eqref{KdV} becomes
\begin{equation}\label{KdV - general computational variables}
\frac{u_\tau}{t_\tau} + \bigg(u- \frac{x_\tau}{t_\tau}\bigg) \frac{u_s}{x_s} + \frac{1}{x_s}\bigg(\frac{1}{x_s}\bigg(\frac{u_s}{x_s} \bigg)_s \bigg)_s = 0,
\end{equation}
together with the companion equations \eqref{KdV t companion equations}.  The differential equation \eqref{KdV - general computational variables} is reminiscent of the equation one obtains when writing the KdV equation in \emph{Lagrangian form}, \cite{BCW-2015}.  In the classical Lagrangian framework, 
the differential constraint
\begin{equation}\label{KdV x companion equation}
\frac{x_\tau}{t_\tau} = u,
\end{equation}
is also imposed.  The KdV equation then reduces to
\begin{equation}\label{KdV Lagrangian coordinates}
\overline{\Delta} = \frac{u_\tau}{t_\tau} + \frac{1}{x_s}\bigg(\frac{1}{x_s}\bigg(\frac{u_s}{x_s} \bigg)_s \bigg)_s = 0,
\end{equation}
together with the companion equations \eqref{KdV t companion equations}, \eqref{KdV x companion equation}.  In particular, when $k=1$ in \eqref{t solution}, we obtain the system of differential equations
\[
u_\tau + \frac{1}{x_s}\bigg(\frac{1}{x_s}\bigg(\frac{u_s}{x_s} \bigg)_s \bigg)_s = 0,\qquad x_\tau = u.
\]
\end{example}

\subsection{Finite difference equations}\label{section:difference equations}

We now move on to the discrete setting, which is the main focus of these lecture notes.  In the previous section, we introduced two different jets spaces, namely $\J\n$ and $\mJ\n$.  The motivation for introducing computational variables and the corresponding jet space $\mJ\n$ stems from the fact that the discrete framework is more closely related to $\mJ\n$ than $\J\n$.  

Let $N =  (n^1,\ldots,n^p) \in \mathbb{Z}^p$ denote an integer-valued multi-index.  Thinking of the multi-index $N$ as sampling the computational variables $s=(s^1,\ldots,s^p) \in \mathbb{R}^p$ at integer values, the discrete notation $(x_N,u_N)$ should be understood as sampling the submanifold $S = \{(x(s),u(s))\} \subset M$ at the integer-valued points $s=N \in \mathbb{Z}^p \subset \mathbb{R}^p$.  In other words $(x_N,u_N) = (x(N),u(N))$.
To approximate the $\on$-jet $(s,x\n,u\n) \in \mathcal{J}\n$ at $s=N$, we consider a finite collection of points 
\begin{equation}\label{discrete jet}
(N,x_N\sn,u_N\sn)=(N,\ldots, x_{N+K},\ldots,u_{N+K},\ldots), 
\end{equation}
where $K \in \mathbb{Z}^p$.  We require that the point $(x_N,u_N)$ is always included and that 
\[
x_{N+K_1} \neq x_{N+K_2}\qquad \text{whenever}\qquad K_1 \neq K_2,
\]
so that no two discrete independent variables are the same.  We refer to \eqref{discrete jet} as the \emph{$\on^\text{th}$ order discrete jet at $N$}.  In numerical analysis, a point in \eqref{discrete jet} is also called a  \emph{stencil}.  For theoretical purposes, one can assume that the multi-index $K \in (\mathbb{Z}^{\geq 0})^p$ only takes non-negative values and that $0 \leq \# K=k^1+\cdots+k^p \leq \on$.  The latter provides the minimal number of points required to approximate the $\on$-jet $(x,u\n)$ (or $(s,x\n,u\n)$) by first order forward differences.  In applications, especially when constructing numerical schemes,  it is generally preferable to consider points centered around $(x_N,u_N)$ and to include more than the minimum number of points in $(N,x_N\sn,u_N\sn)$ required to approximate $(x,u\n)$ for better numerical accuracy and stability.   From now on, we will assume that a certain stencil \eqref{discrete jet} has been chosen.  We denote by
\[
\mJ\sn = \bigcup_{N \in \mathbb{Z}^p} (N,x\sn_N,u_N\sn)
\]
the union over all the stencils and call $\mJ\sn$ the \emph{$\on^\text{th}$ order discrete jet space} as $\mJ\sn$ provides an approximation of $\mJ\n$.  Since the jet coordinates of $\J\n$ can be expressed in terms of the jet coordinates of $\mJ\n$ using \eqref{uxJ}, it follows that the points in $\mJ\sn$ can be used to approximate $\J\n$.

\begin{example}\label{derivative discretization example}
Consider the case where $p=2$ and the dimension of the manifold $M$ is $\text{dim}\,M = m = 3$.  Let $(t,x,u)$ be local coordinates on $M$.   In the continuous case, see Example \ref{computation variables jet example}, we introduced the computational variables $(\tau, s)$.  In the discrete case, let $N = (n,i) \in \mathbb{Z}^2$, which can be thought of as evaluating the computational variables $(\tau,s)$ at integer values.  To make the multi-index notation more compact, we let
\[
(t_N,x_N,u_N) = (t^n_i,x_i^n,u^n_i)\qquad N=(n,i)\in \mathbb{Z}^2.
\]
Working with forward differences, the simplest first order discrete jet is parametrized by
\[
(t_N^{[1]},x_N^{[1]},u_N^{[1]}) = (t^n_i,x_i^n,u^n_i,t^{n+1}_i,x_i^{n+1},u^{n+1}_i,t^n_{i+1},x_{i+1}^n,u^n_{i+1}).
\]
First order approximations of the first order derivatives $(t_\tau,x_\tau,u_\tau)$ and $(t_s,x_s,u_s)$ on a grid with unit spacing are then given by
\begin{equation}\label{first order computation variable derivatives}
\begin{aligned}
(t_\tau,x_\tau,u_\tau) &\approx (t^{n+1}_i - t^n_i, x^{n+1}_i - x^n_i,u^{n+1}_i - u^n_i),\\
(t_s,x_s,u_s) &\approx (t^n_{i+1} - t^n_i, x^n_{i+1} - x^n_i,u^n_{i+1} - u^n_i).
\end{aligned}
\end{equation}
Referring to \eqref{ux-ut} for the expressions of the $t$ and $x$ derivatives of $u$ in terms of the computational variable derivatives, and using \eqref{first order computation variable derivatives} we have that
\begin{equation}\label{first order derivative approximations}
\begin{aligned}
u_x &= \frac{t_\tau\, u_s - t_{s}\, u_\tau}{x_s t_\tau - x_\tau t_s} 
\approx \frac{(t^{n+1}_i-t^n_i)(u^n_{i+1}-u^n_i)-(t^n_{i+1}-t^n_i)(u^{n+1}_i-u^n_i)}{(x^n_{i+1}-x^n_i)(t^{n+1}_i-t^n_i)-(x^{n+1}_i-x^i_n)(t^n_{i+1}-t^n_i)},\\
u_t &= \frac{x_s\, u_\tau - x_\tau\, u_s}{x_s t_\tau - x_\tau t_s} 
\approx \frac{(x^n_{i+1}-x^n_i)(u^{n+1}_i-u^n_i)-(x^{n+1}_i-x^n_i)(u^n_{i+1}-u^n_i)}{(x^n_{i+1}-x^n_i)(t^{n+1}_i-t^n_i)-(x^{n+1}_i-x^i_n)(t^n_{i+1}-t^n_i)}.
\end{aligned}
\end{equation}
The latter expressions are first order forward approximations of the first order partial derivatives $u_x$ and $u_t$ on any mesh that satisfies
\[
\det \begin{bmatrix}
(t^{n+1}_i-t^n_i) & (t^n_{i+1} - t^n_i) \\
(x^{n+1}_i-x^n_i) & (x^n_{i+1} - x^n_i)
\end{bmatrix}= (x^n_{i+1}-x^n_i)(t^{n+1}_i-t^n_i)-(x^{n+1}_i-x_i^n)(t^n_{i+1}-t^n_i)\neq 0,
\]
the latter being a discrete version of the non-degeneracy condition \eqref{non-degeneracy condition}.  The procedure can be repeated to obtain approximations of higher order derivatives on arbitrary meshes.  For example, applying the implicit derivative operators \eqref{D_x,D_t} to the first order derivative expressions \eqref{ux-ut} one obtains formulas for the second order derivatives $u_{xx}$, $u_{xt}$, and $u_{tt}$ expressed in terms of the computational derivatives.  Substituting the approximations \eqref{first order computation variable derivatives} and the second order derivative approximations 
\[ 
\begin{aligned}
&t_{\tau \tau} \approx t^{n+2}_i - 2t^{n+1}_i + t^n_i,
& &x_{\tau \tau} \approx x^{n+2}_i - 2x^{n+1}_i + x^n_i,\\
&t_{\tau s} \approx t^{n+1}_{i+1} - t^{n+1}_i-t^n_{i+1} + x^n_i,
&\qquad &x_{\tau s} \approx x^{n+1}_{i+1} - x^{n+1}_i-x^n_{i+1} + x^n_i,\\
&t_{ss} \approx t^n_{i+2} - 2t^n_{i+1} + t^n_i,& & x_{ss} \approx x^n_{i+2} - 2x^n_{i+1} + x^n_i,
\end{aligned}
\]
into the formulas obtained yields discrete approximations for $u_{xx}$, $u_{xt}$, and $u_{tt}$ in the computational variables on an orthogonal grid with unit spacing.
\end{example}

\begin{definition}
A \emph{finite difference equation} is the zero locus of a discrete map $E\colon \mJ\sn \to \mathbb{R}$.  That is,
\[
E(N,x_N\sn,u_N\sn) = 0.
\]
\end{definition}

\begin{definition}\label{regular definition difference equations}
A finite difference equation $E\colon \mJ\sn \to \mathbb{R}$ is said to be \emph{regular} if the rank of the differential
\[
dE = \sum_K \bigg[\sum_{i=1}^p \pp{E}{x^i_{N+K}} dx^i_{N+K} + \sum_{\alpha=1}^q \pp{E}{u^\alpha_{N+K}}du^\alpha_{N+K}\bigg]
\]
is constant for all $N$ in the domain of definition of the equation.
\end{definition}

Finite difference equations can be studied as mathematical objects of interest in their own, \cite{E-2005,H-2014,LL-2011}.  In the following we are interested in finite difference equations that approximate differential equations. 

\begin{definition}
A finite difference equation $E(N,x\sn_N,u\sn_N)=0$ is said to be \emph{consistent} with the differential equation $\Delta(x,u\n)=0$ (or $\overline{\Delta}(s,x\n,u\n)=0)$ if in the continuous limit $(x_{N+K},u_{N+K}) \to (x_N,u_N)$,
\[
E(N,x\sn_N,u\sn_N) \to \Delta(x,u\n)\qquad (\text{or } E(N,x\sn_N,u\sn_N) \to \overline{\Delta}(s,x\n,u\n)).
\]
\end{definition}

\begin{remark}
The process of taking continuous limits is discussed in more details in Section \ref{symmetry-preserving schemes section}.
\end{remark}

\begin{definition}\label{numerical scheme definition}
Let $\Delta(x,u\n)=0$ be a differential equation with extended system $\{\overline{\Delta}(s,x\n,u\n)=0,\, \widetilde{\Delta}(s,x\n,u\n)=0\}$.  A \emph{numerical scheme} is a system of finite difference equations
\[
\overline{E}(N,x\sn_N,u\sn_N) = 0,\qquad \widetilde{E}(N,x\sn_N,u\sn_N)=0,
\]
where $\overline{E}(N,x\sn_N,u\sn_N) = 0$ approximates the differential equation 
\[
\Delta(x,u\n)=\overline{\Delta}(s,x\n,u\n)=0
\]
and the equations $\widetilde{E}(N,x\sn_N,u\sn_N)=0$ provide an approximation of the companion equations 
\[
\widetilde{\Delta}(s,x\n,u\n)=0.
\]
\end{definition}

Intuitively, the difference equations $\widetilde{E}(N,x\sn_N,u\sn_N)=0$ provide constraints on the mesh used to approximate the differential equation $\Delta=0$.  The latter should not yield any restrictions on the discrete dependent variables $u_N\sn$.

\begin{example}\label{KdV standard scheme example}
To illustrate Definition \ref{numerical scheme definition}, let us consider the KdV equation \eqref{KdV}.  Assume the equation is to be discretized on the orthogonal mesh
\begin{equation}\label{orthogonal mesh}
t^n = k \, n + t^0,\qquad x_i = h\, i + x_0,
\end{equation}
where $k,h>0$, $(n,i) \in \mathbb{Z}^2$, and $t^0$, $x_0$ are arbitrary constants.  The mesh \eqref{orthogonal mesh} can be encapsulated in a system of finite difference equations in different ways.  For example, it is not difficult to see that \eqref{orthogonal mesh} is the solution to the system of equations
\begin{equation}\label{orthogonal mesh equations 1}
\begin{aligned}
&t^{n+1}_i-t^n_i = k, &\qquad & x^{n+1}_i - x^n_i = 0,\\
&t^n_{i+1} - t^n_i = 0, & & x^n_{i+1} - x^n_i = h.
\end{aligned}
\end{equation}
The mesh \eqref{orthogonal mesh} is also a solution to
\begin{equation}\label{orthogonal mesh equations 2}
\begin{aligned}
&t^{n+1}_i-2t^n_i + t^{n-1}_i = 0, &\qquad & x^{n+1}_i - x^n_i = 0,\\
&t^n_{i+1} - t^n_i = 0, & & x^n_{i+1} - 2x^n_i + x^n_{i-1} = 0.
\end{aligned}
\end{equation}
The difference between the two systems of mesh equations is that in \eqref{orthogonal mesh equations 1} the time step $k$ and the spatial step $h$ are fixed by the system whereas in \eqref{orthogonal mesh equations 2} those steps corresponds to constants of integration.  In both cases, the KdV equation can be approximated by 
\begin{equation}\label{naive KdV discretization}
\frac{u^{n+1}_i - u^n_i}{k} + u^n_i \cdot \frac{u^n_{i+1} - u^n_{i-1}}{2h} + \frac{u^n_{i+2}-2u^n_{i+1} + 2u^n_{i-1}-u^n_{i-2}}{2h^3} = 0.
\end{equation}
The systems of equations \eqref{orthogonal mesh equations 1}--\eqref{naive KdV discretization} or \eqref{orthogonal mesh equations 2}--\eqref{naive KdV discretization} provide two examples of Definition \ref{numerical scheme definition}.  They also illustrate the fact that, in general, the equations $\widetilde{E}=0$ specifying the mesh are not unique.  
\end{example}

\section{Lie symmetries}\label{symmetry section}

Let $G$ be an $r$-dimensional Lie group, and let $\mM$ be a $d$-dimensional manifold with local coordinates $z=(z^1,\ldots,z^d)$.  In the following, the manifold $\mM$ can represent the submanifold jet spaces $\J\n$ or $\mJ\n$ or the discrete jet space $\mJ\sn$.  In the latter case, $\mM$ should in fact be called a \emph{lattifold} or \emph{lattice variety}, that is a manifold-like structure modeled on $\mathbb{Z}^p$, \cite{BM,HM-2008}.

\begin{definition}\label{transformation group definition}
A \emph{transformation group} acting on a manifold $\mM$ is given by a Lie group $G$ and a smooth map $\Phi\colon G\times \mM \to \mM$, such that $\Phi(g,z) = g\cdot z$,  which satisfies the following two properties
\begin{equation}\label{transformation group properties}
e \cdot z = z,\qquad g\cdot(h\cdot z) = (gh)\cdot z,\qquad \text{for all}\qquad z\in \mM,\; g,\,h\in G,
\end{equation}
and where $e\in G$ denotes the identity element.
\end{definition}

It follows from \eqref{transformation group properties} that the inverse of the transformation defined by the group element $g$ is given by the inverse group element $g^{-1}$.  Therefore $g$ induces a diffeomorphism from $\mM$ to itself.  

\begin{remark}
Definition \ref{transformation group definition} assumes that the group action is \emph{global}, meaning that $g\cdot z$ is defined for every $g \in G$ and every $z\in \mM$.  In practice, group actions may only be defined \emph{locally}, meaning that for a given $z\in \mM$, the transformation $g\cdot z$ is only defined for group elements $g$ sufficiently near the identity.  For a \emph{local transformation group}, the map $\Phi$ is defined on an open subset $\mathcal B$ with $\{e\} \times \mM \subset \mathcal{B} \subset G\times \mM$, and the conditions \eqref{transformation group properties} of Definition \ref{transformation group properties} are imposed wherever they are defined.
\end{remark}

In the following, we use capital letters to denote the image of a point under a group transformation.  For example,
\[
Z = g\cdot z\qquad\text{where}\qquad g\in G\qquad \text{and}\qquad z \in \mM.
\]
At the infinitesimal level, let $\g$ denote the Lie algebra of vector fields corresponding to the infinitesimal generators of the group action.  A vector field
\[
 \vv = \sum_{a=1}^d \> \zeta^a(z) \pp{}{z^a} 
\]
will be in $\g$ if it is tangent to the orbits of a one-parameter subgroup of transformations of $G$.   The \emph{flow} through the point $z\in \mM$ generated by a vector field $\vv \in \g$, is found by solving the initial value problem
\[
\frac{dZ^a}{d\epsilon} = \zeta^a(Z),\qquad Z^a(0) = z^a,\qquad a=1,\ldots,d.
\] 
The maximal integral curve is denoted $\exp[\epsilon \vv]\cdot z$, and is called the \emph{exponentiation} of the infinitesimal generator $\vv$.

\begin{definition}\label{symmetry group definition}
Let $G$ be a local Lie group of transformations acting on $\mM$.  The Lie group $G$ is a (local) \emph{symmetry group} of the (fully) regular equation\footnote{Depending whether $\mM$ represents $\J\n$, $\mJ\n$, or $\mJ\sn$, we refer to Definitions \ref{regular definition differential equations},  \ref{regular definition differential equations 2}, or \ref{regular definition difference equations} for the notion of (fully) regular equation.} $F(z)=0$ if and only if 
\[
F(g\cdot z)=0\qquad \text{whenever}\qquad F(z)=0,
\]
for all $g \in G$ such that the local action is defined.  Infinitesimally, a connected Lie group of transformations $G$ acting on $\mM$ is a local symmetry group of $F(z)=0$ if and only if
\begin{equation}\label{infinitesimal symmetry criterion}
\vv(F)\big|_{F=0}=0\qquad \text{for all}\qquad \vv \in \g.
\end{equation}
\end{definition}

\begin{remark}
Definition \ref{symmetry group definition} extends to systems of equations and more general local groups of transformations by including discrete transformations as well, \cite{BDP-2015,H-1998,H-2000,H-2000-2}. In the following we restrict all our considerations to Lie point symmetries and omit the interesting case of discrete symmetries.
\end{remark}

\subsection{Symmetries of differential equations}\label{symmetry of differential equations section}

Symmetries of differential equations are covered extensively in many excellent textbooks such as \cite{BA-2002,BCA-2010,H-2000,O-1993,O-2009,O-1982}.  We refer to these references for a more detailed exposition.


If $\mM = \J\n$, then the local group action is given by the \emph{prolonged action} $(X,U\n) = g\cdot (x,u\n)$ on the submanifold $\on$-jet $(x,u\n)$.  Let 
\begin{equation}\label{group action on M}
X^i=g\cdot x^i,\quad i=1,\ldots,p,\qquad U^\alpha = g\cdot u^\alpha,\quad \alpha =1,\ldots,q
\end{equation}
denote the local group action of $G$ on the manifold $M$ locally coordinatized by $(x,u)$.   To compute the prolonged action, we introduce  the \emph{implicit differentiation operators}, \cite{FO-1999},
\begin{equation}\label{implicit derivative operators}
D_{X^i} = \sum_{j=1}^p\> W^j_i D_{x^j},\qquad \text{where}\qquad (W^j_i) = \bigg(\pp{X^j}{x^i}\bigg)^{-1}
\end{equation}
denotes the entries of the inverse Jacobian matrix and 
\[
D_{x^j} = \pp{}{x^j} + \sum_{J} \sum_{\alpha=1}^q\> u^\alpha_{x^{J+e_j}} \pp{}{u^\alpha_{x^J}}
\] 
is the total differentiation operator with respect to the independent variable $x^j$.  In the above formula, $e_j = (0,\ldots,0,1,0,\ldots,0)\in \mathbb{R}^p$ denotes the unit vector with zeros everywhere except in the $j^\text{th}$ component. We note that the operators \eqref{implicit derivative operators} mutually commute
$$
[D_{X^i},D_{X^j}]=0,\qquad 1\leq i,\,j \leq p.
$$
Successively applying the implicit differentiation operators \eqref{implicit derivative operators} to $U^\alpha = g\cdot u^\alpha$ yields the expressions for the prolonged action
\[
U^\alpha_{X^J} = (D_{X^1})^{j_1}\cdots(D_{X^p})^{j_p} U^\alpha,\qquad \alpha=1,\ldots,q,\quad \#J \geq 0.
\]

At the infinitesimal level, let 
\begin{equation}\label{v}
\vv = \sum_{i=1}^p \> \xi^i(x,u)\pp{}{x^i} + \sum_{\alpha=1}^q\> \phi^\alpha(x,u)\pp{}{u^\alpha}
\end{equation}
denote an infinitesimal generator of the group action \eqref{group action on M}.  The prolongation of \eqref{v} to $\J\n$ is given by
\[
\vv\n = \sum_{i=1}^p\> \xi^i\pp{}{x^i} + \sum_{\alpha=1}^q\sum_{J} \> \phi^{\alpha;J} \pp{}{u^\alpha_{x^J}},
\]
where the prolonged vector field coefficients are defined recursively by the standard prolongation formula
\[
\phi^{\alpha;J+e_j} = D_{x^j}\phi^{\alpha;J} - \sum_{i=1}^p\> (D_{x^j}\xi^i)\, u^\alpha_{x^{J+e_i}}.
\]

Given a differential equation $\Delta(x,u\n)=0$, the Lie point symmetries of the equation are found from the \emph{infinitesimal invariance criterion} 
\begin{equation}\label{Delta infinitesimal symmetry criterion}
\vv\n(\Delta)\big|_{\Delta=0}=0\qquad \text{for all}\qquad \vv \in \g.
\end{equation}
The latter yields a differential equation in $x$, $u$ and the derivatives of $u$ with respect to $x$, as well as $\xi^i(x,u)$ and $\phi^\alpha(x,u)$ and their partial derivatives with respect to $x$ and $u$.  After eliminating any dependencies among the derivatives of the $u$'s due to the equation $\Delta(x,u\n)=0$, one can equate the coefficients of the remaining unconstrained partial derivatives of $u$ to zero.  This yields a system of linear partial differential equations for the coefficients $\xi^i$ and $\phi^\alpha$, called the \emph{determining equations} of the (maximal) Lie symmetry algebra.  The procedure for obtaining and solving the determining equations has been implemented in all major computer algebra systems such as {\sc Macsyma, Maple, Mathematica, MuMath} and {\sc Reduce}.  An extensive list of references on the subject can be found in \cite{CV-2000}.

\begin{example}\label{KdV symmetry computation example}
To illustrate the algorithm outlined above, we compute the infinitesimal generators of the KdV equation \eqref{KdV}.  Let $\vv = \xi(t,x,u)\partial_x + \eta(t,x,u) \partial_t + \phi(t,x,u)\partial_u$ denote a general vector field on $\mathbb{R}^3$.  The third order prolongation of $\vv$ is
\begin{multline*}
\vv^{(3)} = \xi \pp{}{x} + \eta \pp{}{t} + \phi \pp{}{u} + \phi^x \pp{}{u_x} + \phi^t \pp{}{u_t} + \phi^{xx}\pp{}{u_{xx}} + \phi^{xt}\pp{}{u_{xt}} + \phi^{tt}\pp{}{u_{tt}} \\
+ \phi^{xxx}\pp{}{u_{xxx}} + \phi^{xxt}\pp{}{u_{xxt}} + \phi^{xtt}\pp{}{u_{xtt}} +  \phi^{ttt}\pp{}{u_{ttt}},
\end{multline*}
where
\begin{align}
&\phi^t = D_t\phi - u_x D_t\xi - u_t D_t\eta,\nonumber \\
&\phi^x = D_x\phi - u_x D_x\xi - u_t D_x\eta, \nonumber \\
&\phi^{xx} = D_x(\phi^x) - u_{xx} D_x\xi - u_{xt} D_x \eta \label{explicit prolonged coefficients}\\
& \hskip 0.5cm= D_x^2\phi - u_x D_x^2\xi -u_t D_x^2 \eta - 2u_{xx} D_x\xi - 2u_{xt} D_x\eta, \nonumber \\
&\phi^{xxx} = D_x(\phi^{xx}) - u_{xxx} D_x\xi - u_{xxt} D_t \eta \nonumber \\
& \hskip 0.6cm= D_x^3\phi - u_x D_x^3 \xi - u_t D_x^3\eta - 3u_{xx} D_x^2\xi - 3u_{xt} D_x^2\eta - 3u_{xxx}D_x \xi - 3u_{xxt} D_x\eta. \nonumber
\end{align}
Applying the infinitesimal invariance criterion \eqref{infinitesimal symmetry criterion} to the KdV equation \eqref{KdV} we obtain
\begin{equation}\label{KdV infinitesimal criterion}
\phi^t + u \phi^x + u_x\phi + \phi^{xxx} = 0,
\end{equation}
where $u$ satisfies \eqref{KdV}.  Substituting the expressions \eqref{explicit prolonged coefficients} into \eqref{KdV infinitesimal criterion} and replacing $u_t$ by $-u u_x - u_{xxx}$, we obtain the determining equations of the Lie symmetry algebra, which we now solve.  Firstly, the coefficient of $u_{xxt}$ is $D_x \eta = \eta_x + u_x\eta_u$ which implies that $\eta_x = \eta_u = 0$.  In other words, $\eta = \eta(t)$ is a function of $t$ only\footnote{This is true for all evolution equations.}.  Secondly, the coefficient of $u_{xx}^2$ yields $\xi_u=0$ and thus $\xi=\xi(t,x)$, implying that the admitted Lie symmetries are projectable.  Next, the coefficient of $u_{xxx}$ gives $\eta_t - 3\xi_x=0$.  Integrating the latter with respect to $x$, we find that $\xi= \frac{1}{3}x\, \eta_t + \chi(t)$.  The coefficient of $u_{xx}$ implies that $\phi_{uu} = \phi_{xu} = 0$ so that $\phi= \sigma(t) u + \varphi(t,x)$.  Next the coefficient in $u_x$ yields the equation
\[
-\xi_t  + u (\eta_t - \xi_x) + \phi = 0.
\]
Substituting the expressions for $\xi$ and $\phi$, we find
\[
\sigma = -\frac{2}{3} \eta_t\qquad\text{and}\qquad \varphi = \frac{1}{3}x\,\eta_{tt}+\chi_t \qquad \text{so that}\qquad \phi = -\frac{2}{3} u\, \eta_t+ \frac{1}{3}x\,\eta_{tt} + \chi_t.
\]
Finally, the term with no derivatives of $u$ gives $\phi_t + \phi_{xxx} + u\, \phi_x = 0$, which after substitution yields
\[
-\frac{1}{3}u\,\eta_{tt} + \frac{1}{3}x\, \eta_{ttt} + \chi_{tt} = 0.
\]
Since $\eta=\eta(t)$ and $\chi=\chi(t)$ are functions of $t$, the latter equation holds for all $(t,x,u)$ provided that $\eta_{tt} = \chi_{tt} = 0$.  Therefore,
\[
\xi = c_1 + c_2\, t + c_3\, x,\qquad \eta=c_4 + 3 c_3\, t,\qquad \phi = c_2 - 2 c_3\, u,
\]
and the maximal Lie symmetry algebra is spanned by the four vector fields
\begin{equation}\label{KdV symmetry generators}
\begin{aligned}
\vv_1 &= \pp{}{x},&\longrightarrow&\hskip 0.5cm  \text{space translations}, \\
\vv_2 &= \pp{}{t},&\longrightarrow&\hskip 0.5cm  \text{time translations},\\
\vv_3 &= t\pp{}{x} + \pp{}{u},&\longrightarrow&\hskip 0.5cm \text{Galilean boosts}, \\ 
\vv_4 &= x\pp{}{x} + 3t \pp{}{t} - 2u \pp{}{u},&\longrightarrow&\hskip 0.5cm  \text{scalings}.
\end{aligned}
\end{equation}
\end{example}

\begin{exercise}\label{u_xx=0 exercise}
Show that the symmetry group of the ordinary differential equation $u_{xx} = 0$ is eight-dimensional, and generated by
\begin{equation}\label{u_xx=0 symmetry generators}
\begin{aligned}
&\pp{}{x},&\qquad & x\pp{}{x},&\qquad & u\pp{}{x},&\qquad & x^2\pp{}{x}+xu\pp{}{u},\\
&\pp{}{u},& & x\pp{}{u},& & u\pp{}{u},& & xu\pp{}{x}+u^2\pp{}{u}.
\end{aligned}
\end{equation}
Show that the corresponding group of local transformations is the projective group $SL(3,\mathbb{R})$ acting via fractional linear transformations
$$
X = \frac{\epsilon_1x+\epsilon_2u+\epsilon_3}{\epsilon_7x+\epsilon_8u+\epsilon_9},\qquad U = \frac{\epsilon_4x+\epsilon_5u+\epsilon_6}{\epsilon_7x+\epsilon_8u+\epsilon_9},\qquad 
\det\begin{bmatrix} \epsilon_1 & \epsilon_2 & \epsilon_3 \\ \epsilon_4 & \epsilon_5 & \epsilon_6 \\ \epsilon_7 & \epsilon_8 & \epsilon_9 \end{bmatrix} = 1,
$$
where $\epsilon_1,\ldots,\epsilon_9\in\mathbb{R}$ are group parameters.
\end{exercise}

\begin{exercise}
Consider the Schwarzian differential equation
\begin{equation}\label{Schwarzian equation}
\frac{u_x\, u_{xxx} - (3/2)u_{xx}^2}{u_x^2} = F(x),
\end{equation}
where $F(x)$ is an arbitrary function.
\begin{enumerate}
\item Find the determining equations for the vector fields spanning the maximal Lie symmetry algebra and show that a basis is given by
\begin{equation}\label{sl2 generators}
\vv_1 = \pp{}{u},\qquad \vv_2 = u\pp{}{u},\qquad \vv_3= u^2\pp{}{u}.
\end{equation}
\item Show that the corresponding local Lie group of transformations is
\begin{equation}\label{sl2 action}
X = x,\qquad U = \frac{au + b}{c u + d},\qquad \text{with}\qquad ad-bc=1.
\end{equation}
\item When $F(x)\equiv 0$ is identically zero, show that the maximal Lie symmetry algebra is four-dimensional and determine a basis. Also find the corresponding finite group transformations.
\end{enumerate}
\end{exercise}

\begin{exercise}\label{burgers symmetry exercise}
Show that the maximal Lie symmetry algebra of Burgers' equation
\begin{equation}\label{burgers equation}
u_t + uu_x = \nu u_{xx},
\end{equation}
where $\nu>0$ denotes the viscosity, is spanned by the vector fields
\begin{equation}\label{burgers symmetry generators}
\begin{aligned}
&\vv_1 = \pp{}{x},\quad\; \vv_2 = \pp{}{t},\quad\; \vv_3 = t \pp{}{x} + \pp{}{u},\\
&\vv_4 = x\pp{}{x} + 2t\pp{}{t} - u\pp{}{u},\quad\; \vv_5=tx \pp{}{x} + t^2 \pp{}{t} + (x - tu)\pp{}{u}.
\end{aligned}
\end{equation}
\end{exercise}

In the computational variable framework, the local transformation group $G$ acting on the manifold $M$ is trivially extended to the computational variables.  That is,
\[
g\cdot(s,x,u) = (s,g\cdot x, g\cdot u).
\]
The prolongation of an infinitesimal generator \eqref{v} to $\mathcal{J}\n$ is then simply given by
\[
\vv\n = \sum_{i=1}^p \sum_J\> D_s^J\xi^i\pp{}{x^i_{s^J}} + \sum_{\alpha=1}^q \sum_J\> D_s^J\phi^\alpha\pp{}{u^\alpha_{s^J}},
\]
where $D_s^J=(D_{s^1})^{j_1}\cdots (D_{s^p})^{j_p}$ denotes the total differentiation operator in the computational variables $s=(s^1,\ldots,s^p) $ with 
\[
D_{s^j} = \pp{}{s^j} + \sum_{i=1}^p \sum_J x^i_{s^{J+e_j}} \pp{}{x^i_{s^J}} + \sum_{\alpha=1}^q \sum_J u^\alpha_{s^{J+e_j}} \pp{}{u^\alpha_{s^J}},\qquad j=1,\ldots,p.
\]

\begin{definition}
Let $G$ be a symmetry group of the differential equation $\Delta(x,u\n)=0$.  The extended system of differential equations 
\begin{equation}\label{extended system}
\{\overline{\Delta}(s,x\n,u\n)=0,\, \widetilde{\Delta}(s,x\n,u\n)=0\}
\end{equation}
is said to be \emph{$G$-compatible} if $G$ is a symmetry group of \eqref{extended system}.  That is,
\[
\begin{cases}\overline{\Delta}(s,g\cdot x\n,g\cdot u\n)=0,\\ \widetilde{\Delta}(s,g\cdot x\n,g\cdot u\n)=0,\end{cases} \qquad \text{whenever}\qquad 
\begin{cases}\overline{\Delta}(s,x\n,u\n)=0,\\ \widetilde{\Delta}(s,x\n,u\n)=0,\end{cases}
\]
and where the prolonged action is defined.  At the infinitesimal level,
\[
\vv\n(\overline{\Delta})\big|_{\{\overline{\Delta}=0,\widetilde{\Delta}=0\}}=0\qquad \text{and}\qquad
\vv\n(\widetilde{\Delta})\big|_{\{\overline{\Delta}=0,\widetilde{\Delta}=0\}}=0
\]
for all infinitesimal generators $\vv \in \g$.
\end{definition}

\begin{example}\label{KdV symmetry example}
Recall from Example \ref{KdV symmetry computation example} that the KdV equation \eqref{KdV} is invariant under a four-dimensional maximal Lie symmetry group whose associated algebra of infinitesimal generators is spanned by the vector fields \eqref{KdV symmetry generators}.  In the computational variables $(\tau,s)$ introduced in Example \ref{KdV in computational variables}, the first prolongation of the infinitesimal generators \eqref{KdV symmetry generators} is given by
\begin{align*}
\vv_1^{(1)} &= \pp{}{x},\qquad \vv_2^{(1)} = \pp{}{t},\qquad
\vv_3^{(1)} = t \pp{}{x} + \pp{}{u} + t_s \pp{}{x_s} + \pp{}{u_s} + t_\tau \pp{}{x_\tau} + \pp{}{u_\tau},\\
\vv_4^{(1)} &= x \pp{}{x} + 3t \pp{}{t} - 2u \pp{}{u} + x_s \pp{}{x_s} + 3t_s \pp{}{t_s} - 2u_s \pp{}{u_s} \\
& + x_\tau \pp{}{x_\tau} + 3t_\tau \pp{}{t_\tau} - 2u_\tau \pp{}{u_\tau}.
\end{align*}
By direct computation, it is not hard to verify that for the differential equation \eqref{KdV Lagrangian coordinates}
\[
\vv_1^{(2)}[\overline{\Delta}]\big|_{\overline{\Delta}=0} = \vv_2^{(2)}[\overline{\Delta}]\big|_{\overline{\Delta}=0} = \vv_3^{(2)}[\overline{\Delta}]\big|_{\overline{\Delta}=0} = \vv_4^{(2)}[\overline{\Delta}]\big|_{\overline{\Delta}=0} = 0.
\]
Therefore, equation \eqref{KdV Lagrangian coordinates} is invariant under the symmetry group of the KdV equation.  Also,
\[
\vv^{(2)}_\kappa(t_s) =0,\qquad  \vv^{(2)}_\kappa(t_{\tau\tau})= 0,\qquad \vv^{(2)}_\kappa\bigg(\frac{x_\tau}{t_\tau}-u\bigg) = 0,\qquad \kappa=1,2,3,4,
\]
whenever 
\[
t_s = 0,\qquad t_{ss} = 0,\qquad \frac{x_\tau}{t_\tau} - u=0.
\] 
Therefore, the companion equations \eqref{KdV t companion equations}, \eqref{KdV x companion equation} are invariant under the symmetry group of the KdV equation.  The extended system of differential equations \eqref{KdV t companion equations}, \eqref{KdV x companion equation}, \eqref{KdV Lagrangian coordinates} is therefore $G$-compatible with the symmetry group of the KdV equation.
\end{example}

\subsection{Symmetries of finite difference equations}

As in Section \ref{symmetry of differential equations section}, let $G$ be a local Lie group of transformations acting smoothly on the manifold $M$.  The induced action on the discrete $n$-jet $(N,x\sn_N,u\sn_N)$ is given by the product action
\begin{equation}\label{discrete group action}
g\cdot (N,x\sn_N,u\sn_N) = (N, \ldots, g\cdot x_{N+K},\ldots,g\cdot u_{N+K},\ldots),
\end{equation}
where each point $(x_{N+K},u_{N+K})$ is transformed by the same group transformation $g$.  At the infinitesimal level, given the vector field
\begin{equation}\label{discrete vector field}
\vv = \sum_{i=1}^p\> \xi^i_N\pp{}{x^i_N} + \sum_{\alpha=1}^q\> \phi^\alpha_K\pp{}{u^\alpha_N},
\end{equation}
where $\xi^i_N = \xi^i(x_N,u_N)$ and $\phi^\alpha_N = \phi^\alpha(x_N,u_N)$, the prolonged vector field is given by
\[
\vv\sn = \sum_{K}\bigg[\sum_{i=1}^p\> \xi^i_{N+K}\pp{}{x^i_{N+K}} + \sum_{\alpha=1}^q\> \phi^\alpha_{N+K}\pp{}{u^\alpha_{N+K}}\bigg],
\]
which is obtained by adding copies of $\vv$ evaluated at the different points in the discrete jet $(N,x_N\sn,u_N\sn)$.

\begin{remark}
The above considerations can be generalized by allowing the group action \eqref{discrete group action} or the infinitesimal generator \eqref{discrete vector field} to depend on the multi-index $N$.  For example, in \eqref{discrete vector field}, the vector field coefficients could be functions of $N$ so that $\xi_N^i = \xi^i(N,x_N,u_N)$ and $\phi^\alpha_N = \phi^\alpha(N,x_N,u_N)$.  When constructing symmetry-preserving schemes, this more general case does not occur as the transformation group that one considers is the group of point symmetries of the differential equation $\Delta(x,u\n)=0$, which only contains point transformations in the $x$, $u$ variables.
\end{remark}

Using the infinitesimal invariance criterion
\begin{equation}\label{finite difference infinitesimal symmetry criterion}
\vv\sn(E)\big|_{E=0}=0\qquad \text{for all}\qquad \vv \in \g,
\end{equation}
the symmetries of finite difference equations can be computed in a manner similar to the differential case.  Equation \eqref{finite difference infinitesimal symmetry criterion} yields a finite difference equation for the vector field coefficients $\xi^i_N$ and $\phi^\alpha_N$.  Since the invariance condition \eqref{finite difference infinitesimal symmetry criterion} only has to hold on the solution of the difference equation, one must eliminate any dependencies among $(x_N,u_N)$ and their shifts due to the equation $E(N,x_N\sn,u_N\sn)=0$.  Differentiating the resulting equation with respect to the remaining variables sufficiently many times, one obtains a system of differential equations for the vector field coefficients.  Once the differential equations are solved, one will, in general, have to substitute the solution into the original difference equation for the vector field coefficients, or an intermediate equation obtained along the way, and solve the resulting equation to obtain the symmetry generators.   

\begin{example}
As a first example, let us compute the admitted infinitesimal generators of the ordinary difference equation
\begin{equation}\label{second order linear homogeneous ordinary difference equation}
u_{i+2} = a(i) u_{i+1} + b(i) u_i\qquad \text{where}\qquad a(i)b(i)\neq 0\quad\forall\; i \in \mathbb{Z}.
\end{equation}
Let 
\[
\vv = \phi_i \pp{}{u_i}
\]
be a vector field, where we allow $\phi_i=\phi(i,u_i)$ to depend on the discrete index $i \in \mathbb{Z}$.  Applying the infinitesimal invariance criterion \eqref{finite difference infinitesimal symmetry criterion} we obtain the equation
\begin{equation}\label{infinitesimal discrete invariance criterion example}
\phi(i+2,a(i) u_{i+1} + b(i) u_i) = a(i) \phi(i+1,u_{i+1}) + b(i) \phi(i,u_i),
\end{equation}
where we replaced $u_{i+2}$ by the right-hand side of \eqref{second order linear homogeneous ordinary difference equation}.  Applying the differential operator $\frac{1}{b(i)}\partial_{u_i} - \frac{1}{a(i)}\partial_{u_{i+1}}$ to \eqref{infinitesimal discrete invariance criterion example} we obtain the differential--difference equation
\begin{equation}\label{new constraint 1}
-\phi^\prime(i+1,u_{i+1}) + \phi^\prime(i,u_i) = 0,
\end{equation}
where the prime notation means differentiation with respect to the second entry of the function.  Differentiating \eqref{new constraint 1} with respect to $u_i$ we obtain
\[
\phi^{\prime\prime}(i,u_i)=0.
\]
Integrating this equation once, we find that 
\begin{equation}\label{new constraint 2}
\phi^\prime(i,u_i) = \alpha(i),
\end{equation}
for some arbitrary function $\alpha(i)$.  Substituting \eqref{new constraint 2} in \eqref{new constraint 1} yields $\alpha(i+1) = \alpha(i)$.  Thus, $\alpha(i) = c$ is constant.  Integrating \eqref{new constraint 2}, we obtain that  
\begin{equation}\label{phi solution}
\phi(i,u_i) = c\, u_i + \beta(i).
\end{equation}
Substituting \eqref{phi solution} in \eqref{infinitesimal discrete invariance criterion example} we conclude that $\beta(i)$ must be a solution of the equation $\beta_{i+2} = a(i) \beta_{i+1} + b(i) \beta_i$.  Thus, the Lie algebra of infinitesimal symmetry generators is spanned by
\begin{equation*}
\begin{aligned}
\vv_1 &= u_i \pp{}{u_i},&\longrightarrow&\hskip 0.5cm  \text{dilations}, \\
\vv_\beta &= \beta(i) \pp{}{u_i},&\longrightarrow&\hskip 0.5cm  \text{linear superposition of solutions}.
\end{aligned}
\end{equation*}
\end{example}

\begin{example}
As a second example, we consider the autonomous discrete potential Korteweg-de Vries equation (dpKdV)
\begin{equation}
u^{n+1}_{i+1} = u^n_i + \frac{1}{u^n_{i+1} - u^{n+1}_i},
\end{equation}
which can be found in the work of Hirota, \cite{H-1977}.  Let 
\[
\vv = \phi^n_i \pp{}{u^n_i},\qquad \phi^n_i = \phi(i,n,u^n_i)
\]
be a vector field.  Implementing the infinitesimal invariance criterion \eqref{finite difference infinitesimal symmetry criterion}, we obtain the equation
\begin{equation}\label{LSC}
\phi^{n+1}_{i+1} = \phi^n_i + \frac{\phi^{n+1}_i - \phi^n_{i+1}}{(u^n_{i+1}-u^{n+1}_i)^2},
\end{equation}
where
\[
\phi^{n+1}_{i+1} = \phi\bigg(i+1,n+1,u^n_i+\frac{1}{u^n_{i+1} - u^{n+1}_i}\bigg).
\]
Applying the operator $\partial_{u^n_{i+1}}+\partial_{u^{n+1}_i}$ yields
\begin{equation}\label{phi prime difference equation}
\phi^\prime(i,n+1,u^{n+1}_i) - \phi^\prime(i+1,n,u^n_{i+1}) = 0.  
\end{equation}
Differentiating with respect to $u^n_{i+1}$ gives
\[
\phi^{\prime\prime}(i+1,n,u^n_{i+1})=0\qquad \text{so that}\qquad \phi^n_i = \alpha(i,n)u^n_i + \beta(i,n).
\]
Substituting $\phi^n_i$ in \eqref{phi prime difference equation} we obtain the difference equation
\[
\alpha(i,n+1) - \alpha(i+1,n)=0,
\]
which implies that $\alpha(i,n) = \gamma (i+n)$.  Substituting $\phi^n_i$ in \eqref{LSC} yields the constraints
\begin{gather*}
\gamma(i+n+2) = -\gamma(i+n+1) = \gamma(i+n),\\\beta(i+1,n+1)=\beta(i,n),\qquad \beta(i+1,n) = \beta(i,n+1),
\end{gather*}
which imply that
\[
\phi^n_i = c_1 (-1)^{i+n}u^n_i + c_2(-1)^{i+n} + c_3,
\]
where $c_1$, $c_2$, $c_3$ are arbitrary constants.  We conclude that the Lie algebra of infinitesimal symmetry generators is spanned by 
\[
\vv_1 = (-1)^{i+n} u^n_i \pp{}{u^n_i},\qquad
\vv_2 =  (-1)^{i+n}\pp{}{u^n_i},\qquad
\vv_3 = \pp{}{u^n_i}.
\]
These vector fields satisfy the commutation relations
\[
[\vv_1,\vv_2]=-\vv_3,\qquad [\vv_1,\vv_3]=-\vv_2,\qquad [\vv_2,\vv_3]=0,
\]
which are isomorphic to the commutation relations of the pseudo-euclidean Lie algebra $\mathfrak{e}(1,1)$, \cite{TW-1998}.
\end{example}

For more examples, we refer the reader to \cite{H-2000-2,H-2014,LTW-2000,LTW-2001,LW-2006}.


\section{Invariants}\label{invariants section}

Intuitively, an invariant is a quantity that remains unchanged under the action of a group of local transformations.  In this section, we review two methods for constructing  invariants.  The first approach is based on Lie's infinitesimal invariance criterion which leads to systems of first order partial differential equations that can be solved using the method of characteristics.  The second approach uses the novel theory of equivariant moving frames.  In this framework, invariants are obtained by solving a system of nonlinear algebraic equations.  Remarkably, the latter can be solved for a wide variety of group actions.  

\subsection{Lie's infinitesimal approach} \label{invariants - symmetry generators section}

As in Section \ref{symmetry section}, we consider the differential and finite difference cases simultaneously by considering an $r$-parameter local Lie group $G$ acting on $\mM$, which can represent either $\J\n$, $\mJ\n$ or $\mJ\sn$.

\begin{definition}
A function $I\colon \mM \to \mathbb{R}$ is said to be a \emph{$G$-invariant} if
\begin{equation}\label{group invariance}
I(g\cdot z) = I(z) \qquad \text{for all}\qquad g\in G
\end{equation}
where the action is defined.  At the infinitesimal level, $I\colon \mM \to \mathbb{R}$ is an invariant if
\begin{equation}\label{infinitesimal invariance}
\vv(I) = 0\qquad \text{for all}\qquad \vv \in \g.
\end{equation}
\end{definition}

\begin{remark}\label{invariant remark}
The notion of an invariant is more restrictive than that of an invariant equation.  The invariance of an equation only has to hold on its solution space  whereas the invariance of a function must hold on its domain of definition.  
\end{remark}

Finding invariants from the group invariance condition \eqref{group invariance} can be difficult as the group action is generally nonlinear.  One of the key insights of Sophus Lie was to work with the infinitesimal invariance condition \eqref{infinitesimal invariance} as the latter is a linearized version of the nonlinear problem.  Let 
\begin{equation}\label{g basis}
\vv_\kappa = \sum_{a=1}^d\> \zeta^a_\kappa(z) \pp{}{z^a},\qquad \kappa = 1,\ldots,r=\text{dim}\,\g,
\end{equation}
be a basis of the Lie algebra $\mathfrak g$ of infinitesimal generators of the Lie group action $G$.  To find the functions $I\colon \mM \to \mathbb{R}$ invariant under the group action $G$, we require that the infinitesimal invariance criterion \eqref{infinitesimal invariance} holds for each basis element \eqref{g basis}.  This yields the system of first order linear partial differential equations
\begin{equation}\label{basis symmetry criterion}
\sum_{a=1}^d\> \zeta^a_\kappa(z) \pp{I}{z^a} = 0,\qquad \kappa = 1,\ldots,r.
\end{equation}
The latter is solved using the method of characteristics.  The corresponding characteristic system of ordinary differential equations is
\begin{equation}\label{characteristic system}
\frac{dz^1}{\zeta^1_\kappa(z)} = \frac{dz^2}{\zeta^2_\kappa(z)} = \cdots = \frac{dz^d}{\zeta^d_\kappa(z)},\qquad \kappa=1,\ldots,r,
\end{equation}
and, in the generic case,  the system of equations \eqref{characteristic system} yields a complete set of  $\text{dim}\,\mM -  r$ functionally independent invariants
\[
I^\nu(z),\qquad \nu=1,\ldots,\text{dim}\,\mM - r.
\]

\begin{definition}
A set of invariants $\mathbf{I}_c = \{\ldots,I^\nu(z),\ldots\}$ is said to be \emph{complete} if any  invariant function $I\colon \mM \to \mathbb{R}$ can be expressed in terms of those invariants.  That is,
\[
I(z) = F(\ldots,I^\nu(z),\ldots).
\]
\end{definition}

Most textbooks on symmetries and differential equations cover Lie's infinitesimal method of computing differential invariants,  \cite{BA-2002,BCA-2010,H-2000,O-1993,O-2009,O-1982}. Differential invariants are fundamental objects in mathematics and have many applications.  They occur in geometry as the curvature of curves, surfaces, and submanifolds, \cite{G-1977}, they are used in differential equations to reduce the order of ordinary differential equations and find invariant solution of partial differential equations, \cite{BA-2002,BCA-2010,H-2000,O-1993,O-1982}, their signature manifold is used to solve local equivalence problems, \cite{G-1989,K-1989,O-2009,V-2013}, geometric flows of differential invariants are closely related to completely integrable equations, \cite{BO-2010}, and have applications in computer vision, \cite{K-2009}, climate and turbulence modeling, \cite{BDP-2014}, and much more.

In the finite difference situation, since the Lie group $G$ acts trivially on the multi-index $N$, the components of $N$ will always provide $p$ invariants.  Solving the infinitesimal invariance criterion 
\begin{equation}\label{basis finite difference symmetry criterion}
\vv_\kappa \sn(I) = \sum_{K} \bigg[\sum_{i=1}^p\> \xi^i_{\kappa;N+K} \pp{I}{x_{N+K}} + \sum_{\alpha=1}^q\> \phi^\alpha_{\kappa;N+K} \pp{I}{u^\alpha_{N+K}}\bigg] = 0,
\end{equation}
where $\kappa = 1,\ldots,r,$ will, in the generic case, produce $\dim \mJ\sn - p - r$ difference invariants $I^\nu(x_N\sn,u_N\sn)$ independent of the multi-index $N$.

\begin{example}\label{ex:sl2 invariants}

To illustrate the application of the infinitesimal invariance criterion \eqref{basis finite difference symmetry criterion}, we consider the special linear group $SL(2,\mathbb{R})$ acting on $M=\mathbb{R}^2 = \{(x,u)\}$ by the fractional linear action \eqref{sl2 action} with infinitesimal generators \eqref{sl2 generators}.  For future reference, we consider the order three discrete jet space $\mJ^{[3]}$ with coordinates
\begin{equation}\label{sl2 lattice}
(i,x_{i-1},u_{i-1},x_i,u_i,x_{i+1},u_{i+1},x_{i+2},u_{i+2}).
\end{equation}
To compute a complete set of finite difference invariants on $\mJ^{[3]}$,  we prolong the infinitesimal generators \eqref{sl2 generators} to $\mJ^{[3]}$:
\begin{equation}\label{sl2 discrete prolongation}
\begin{aligned}
\vv_1^{[3]} &= \pp{}{u_{i-1}} + \pp{}{u_i} + \pp{}{u_{i+1}} + \pp{}{u_{i+2}}, \\
\vv_2^{[3]} &= u_{i-1} \pp{}{u_{i-1}} + u_i \pp{}{u_i} + u_{i+1} \pp{}{u_{i+1}} + u_{i+2} \pp{}{u_{i+2}}, \\
\vv_3^{[3]} &= u_{i-1}^2 \pp{}{u_{i-1}} + u_i^2 \pp{}{u_i} + u_{i+1}^2 \pp{}{u_{i+1}} + u_{i+2}^2 \pp{}{u_{i+2}}.
\end{aligned}
\end{equation}
Omitting the trivial invariant given by the index $i$, we expect $(\text{dim}\,\mJ^{[3]}-1)-\text{dim}\,\g = 8-3=5$ functionally independent invariants.  Clearly, four of them are given by
\begin{equation}\label{sl2 x invariants}
x_{i-1},\qquad x_i,\qquad x_{i+1},\qquad x_{i+2}.
\end{equation}
To find the remaining functionally independent invariant $I=I(u_{i-1},u_i,u_{i+1},u_{i+2})$, we first solve the differential equation
\[
\vv_1^{[3]}(I) = \pp{I}{u_{i-1}} + \pp{I}{u_i} + \pp{I}{u_{i+1}} + \pp{I}{u_{i+2}}=0
\]
using the method of characteristics.  The corresponding characteristic system of ordinary differential equations is given by
\[
du_{i-1} = du_i = du_{i+1}=du_{i+2}.
\]
The three functionally independent solutions are
\begin{equation}\label{I variables}
I_{i-1} = u_{i} - u_{i-1},\qquad I_i = u_{i+1} - u_i,\qquad I_{i+1} = u_{i+2} - u_{i+1}.
\end{equation}
The functions \eqref{I variables} form a complete set of difference invariants for the infinitesimal generator $\vv_1$. To proceed further, we notice that any function $I\colon \mJ^{[3]}\to \mathbb{R}$ invariant under all three infinitesimal generators \eqref{sl2 discrete prolongation} must necessarily be a function of the invariants \eqref{I variables},  that is $I=I(I_{i-1},I_i,I_{i+1})$.   Thus to find the functions that are simultaneously invariant under $\vv_1^{[3]}$ and $\vv_2^{[3]}$, we must now restrict the vector field $\vv_2^{[3]}$ to the variables \eqref{I variables}.  The result is
\[
\vv_2^{[3]} = I_{i-1} \pp{}{I_{i-1}} + I_i \pp{}{I_i} + I_{i+1} \pp{}{I_{i+1}}.
\]
Thus, the characteristic system associated with the differential equation
\[
\vv_2^{[3]}(I) = I_{i-1} \pp{I}{I_{i-1}} + I_i \pp{I}{I_i} + I_{i+1} \pp{I}{I_{i+1}}=0
\]
is
\[
\frac{dI_{i-1}}{I_{i-1}} = \frac{dI_i}{I_i} = \frac{dI_{i+1}}{I_{i+1}}.
\]
The two functionally independent solutions are
\begin{equation}\label{J variables}
J_i = \frac{I_{i-1}}{I_i},\qquad J_{i+1} = \frac{I_i}{I_{i+1}}.
\end{equation}
Therefore, any invariant function $I$ must be expressible in terms of \eqref{J variables}.  That it, $I = I(J_i,J_{i+1})$. The restriction of the vector field $\vv_3^{[3]}$ to the variables \eqref{J variables} yields
\[
\vv_3^{[3]} = -I_i\bigg( J_i[J_i+1]\pp{}{J_i} + [1+J_{i+1}]\pp{}{J_{i+1}}\bigg).
\]
Thus, the equation $\vv_3^{[3]}(I)=0$ becomes
\[
J_i[J_i+1]\pp{I}{J_i} + [1+J_{i+1}]\pp{I}{J_{i+1}}=0.
\]
Solving the characteristic system
\[
\frac{dJ_i}{J_i[J_i+1]} = \frac{dJ_{i+1}}{1+J_{i+1}}
\]
we find that the cross-ratio
\begin{equation}\label{R}
R_i = \frac{J_i}{(1+J_i)(1+J_{i+1})} = \frac{(u_i-u_{i-1})(u_{i+2}-u_{i+1})}{(u_{i+1}-u_{i-1})(u_{i+2}-u_i)}
\end{equation}
is an invariant of the $SL(2,\mathbb{R})$ product action on $\mJ^{[3]}$.
\end{example}

\begin{exercise}\label{invariants Lie's approach exercise}
Continuing Exercise \ref{burgers symmetry exercise}, introduce the discrete points $(t^n_i,x^n_i,u^n_i)$, where $(n,i)\in\mathbb{Z}^2$.  
\begin{enumerate}
\item
Verify that the equation
\begin{equation}\label{t constraints exercise}
t^n_{i+1}-t^n_i=0,
\end{equation}
is invariant.  Therefore, Burgers' equation can be invariantly discretized on a mesh with horizontal time layers, the discrete time $t^n$ being only a function of $n\in \mathbb{Z}$. 

\item \label{invariants Lie's approach exercise part b}
Compute a complete set of difference invariants on the lattice
\begin{center}
\tikzstyle{block} = [draw,  
rectangle, minimum height=10pt, minimum width=10pt]
\begin{tikzpicture}
\clip (0,0) rectangle (11, 3.25);
\draw (0,1) -- (13,1);
\draw(2.5,1) node[] {$\bullet$};
\draw(2.5,0.5) node {$(t^n,x^n_{i-1},u^n_{i-1})$};
\draw(6,1) node[] {$\bullet$};
\draw(6,0.5) node {$(t^n,x^n_i,u^n_i)$};
\draw(9,1) node[] {$\bullet$};
\draw(9,0.5) node {$(t^n,x^n_{i+1},u^n_{i+1})$};
\draw (0,3) -- (13,3);
\draw(2.75,3) node[] {$\bullet$};
\draw(2.75,2.5) node {$(t^{n+1},x^{n+1}_{i-1},u^{n+1}_{i-1})$};
\draw(6.5,3) node[] {$\bullet$};
\draw(6.5,2.5) node {$(t^{n+1},x^{n+1}_i,u^{n+1}_i)$};
\draw(9.5,3) node[] {$\bullet$};
\draw(9.5,2.5) node {$(t^{n+1},x^{n+1}_{i+1},u^{n+1}_{i+1})$};
\end{tikzpicture}
\end{center}
using Lie's infinitesimal method.
\end{enumerate}
\end{exercise}

\subsection{Moving frame approach}\label{invariants - moving frame section}

The method of equivariant moving frames is a new theoretical formulation of Cartan's method of moving frames, \cite{BCGGG-1991,G-1989,K-1989,O-2009}.  In this novel framework, moving frames are no longer constrained by frame bundles or connections and can thereby be extended to discrete geometry.  The theory of equivariant moving frames for local Lie group actions was first presented in \cite{FO-1999} and then extended to infinite-dimensional Lie pseudo-group actions in \cite{OP-2008,OP-2009}.  For a comprehensive introduction we refer the reader to the textbook \cite{M-2010}.  In the discrete setting, the theoretical foundations have been expounded in \cite{BM,O-2001}.  

As in the previous sections, our starting point is an $r$-dimensional local Lie group of transformations $G$ acting on the $d$-dimensional manifold $\mM$.

\begin{definition}
A \emph{right moving frame} is a $G$-equivariant map $\rho\colon \mM \to G$.  The $G$-equivariance means that
\[
\rho(g\cdot z) = \rho(z)g^{-1}.  
\]
\end{definition}

\begin{remark}
It is also possible to consider left moving frames.  Given a right moving frame $\rho\colon \mM \to G$, a left moving frame $\overline{\rho}\colon\mM \to G$ is simply given by group inversion, $\overline{\rho} = \rho^{-1}$.  Thus, a left moving frame $\overline{\rho}$ is a $G$-equivariant map satisfying $\overline{\rho}(g\cdot z) = g  \overline{\rho}(z)$.
\end{remark}

To guarantee the existence of a moving frame, the group action must satisfy certain regularity assumptions.

\begin{definition}
A Lie group $G$ is said to act \emph{freely} at $z$ if the isotropy group
\[
G_{z} = \{g \in G\,|\, g\cdot z = z\}
\]
is trivial, i.e.\ $G_{z}=\{ e \}$.  The group action is \emph{locally free} at $z$ if the isotropy group is discrete. The action is (locally) free on $\mM$ if it is (locally) free at all $z \in \mM$.
\end{definition}

When the action is (locally) free, the dimension of the group orbits is constant and equal to $r=\text{dim } G$.

\begin{definition}
A Lie group action is said to be \emph{regular} if the orbits form a regular foliation.
\end{definition}

The main existence theorem for moving frames is given by the following proposition.

\begin{proposition}
If the action of $G$ on $\mM$ is locally free and regular, then a moving frame locally exists on $\mM$.
\end{proposition}

\begin{remark}
Let $\mathcal{V}$ be a connected open submanifold of $\mM$ where a moving frame exists.  By restricting $\mM$ to $\mathcal{V}$, we can always assume that a moving frame is globally defined on $\mM$.  
\end{remark}

In practice, the construction of a moving frame is based on the choice of a cross-section $\mathcal{K}$ to the group orbits.  For simplicity, we assume that $\mathcal{K}$ is a coordinate cross-section, which means that it is specified by fixing some of the coordinates of $z \in \mM$ to constant values:
\begin{equation}\label{cross-section}
\mathcal{K} = \{ z^{a_\kappa} = c^\kappa\;|\; \kappa=1,\ldots,r=\text{dim}\,G \}.
\end{equation}
When the action is free and regular, the right moving frame at $z$ is the unique group element $\rho(z)$ sending $z$ onto the cross-section $\mathcal{K}$, that is $\rho(z)\cdot z \in \mathcal{K}$.  The expressions for the right moving frame are obtained by solving the \emph{normalization equations}
\begin{equation}\label{normalization equations}
g\cdot z^{a_\kappa} = c^\kappa,\qquad \kappa=1,\ldots,r,
\end{equation}
for the group parameters $g=\rho(z)$.   Given a right moving frame, there is a systematic mechanism for constructing invariants known as the \emph{invariantization} procedure.

\begin{definition}
The \emph{invariantization} of a function $F(z)$ is the invariant
\begin{equation}\label{invariantization}
\iota(F)(z) = F(\rho(z) \cdot z).
\end{equation}
\end{definition}

The fact that \eqref{invariantization} is an invariant follows from the $G$-equivariance of the right moving frame:
\[
\iota(F)(g\cdot z) = F(\rho(g\cdot z) \cdot g \cdot z) = F(\rho(z)\cdot g^{-1} \cdot g \cdot z) = F(\rho(z) \cdot z) = \iota(F)(z).
\]

Geometrically, $\iota(F)$ is the unique invariant that agrees with $F$ on the cross-section $\mathcal{K}$.  In particular, the invariantization of an invariant $I$ is the invariant itself, $\iota(I) = I$.  Therefore, the invariantization map $\iota$ defines a canonical projection (depending upon the moving frame) from the space of  functions to the space of invariants.

The invariantization of the components of $z$ is of particular interest.  The invariants $\iota(z^a) = \rho(z)\cdot z^a$, $a=1,\ldots,d,$ are called \emph{normalized invariants}. By the moving frame construction, the invariantization of the component functions defining the cross-section \eqref{cross-section} yields constant invariants, $\iota(z^{a_\kappa})=c^\kappa$. These are called \emph{phantom invariants}.  The following proposition explains why the normalized invariants are important.

\begin{proposition}\label{complete set proposition}
The normalized invariants $\iota(z^a)$, $a=1,\ldots,d$, form a complete set of invariants on $\mM$.
\end{proposition}

Proposition \ref{complete set proposition} follows from the \emph{replacement principle}.  If $I = I(z)$ is an invariant, since $\iota(I)=I$, it follows that 
\[
I(z) = \iota(I)(z) = I(\iota(z)).
\]
In other words, the invariant $I(z)$ can be expressed as a function of the normalized invariants by replacing $z$ with the invariants $\iota(z)$.

For an arbitrary manifold $\mM$, the group action of $G$ on $\mM$, does not have to be free.  On the other hand, when $\mM$ is either $\J\sn$, $\mJ\n$ or $\mJ\sn$, it is always possible, under some mild assumptions, to choose $\ell$ large enough so that the prolonged action becomes (locally) free.  To state the result precisely, we need the following technical definitions.

\begin{definition}
Let $G$ be a local Lie group of transformations acting on the manifold $M$.  The \emph{isotropy subgroup} of a subset $S$ of $M$ is the subgroup
\[
G_S = \{g \in G\;|\; g\cdot S = S\}.
\]
The \emph{global isotropy subgroup} of a subset $S$ of $M$ is the subgroup
\[
G_S^\star = \{g\in G\;|\; g\cdot z = z \text{ for all } z \in S\}.
\]
\end{definition}

\begin{definition}
A local Lie group of transformations $G$ is said to act \emph{effectively on subsets} if, for any open subset $U \subset M$, $G_U^\star =\{e\}$.  The local Lie group acts \emph{locally effectively on subsets} if, for any open subset $U \subset M$, $G_U^\star$ is a discrete subgroup of $G$.
\end{definition}

In the differential case, the following theorem due to Ovsiannikov, \cite{O-1982}, and corrected by Olver, \cite{O-2000}, states that if a group acts (locally) effectively on subsets, then its prolonged action will eventually become free.

\begin{proposition}\label{eventually free action proposition}
If a local Lie group of transformations $G$ acts (locally) effectively on subsets of $M$, then there exists $\ell_0$ such that for all $\ell \geq \ell_0$, the prolonged action of $G$ acts locally freely on an open dense subset $\rm{V}\n \subset \J\n$ (or $\mathcal{V}\n \subset \mJ\n$).
\end{proposition}

The discrete version of Proposition \ref{eventually free action proposition} was proved in \cite{B-2002}. 

\begin{example}
We now implement the moving frame construction for the projective action \eqref{sl2 action} on the order three submanifold jet space 
\[
\J^{(3)}=\{(x,u,u_x,u_{xx},u_{xxx})\}.
\]  
We must therefore compute the prolonged action up to the third derivative.  Since the independent variable $x$ is an invariant of the action \eqref{sl2 action}, the implicit derivative operator \eqref{implicit derivative operators} is
\[
D_X = D_x = \pp{}{x} + u_x\pp{}{u}+ u_{xx} \pp{}{u_x}+u_{xxx}\pp{}{u_{xx}} + \cdots,
\]
and the prolonged action, up to order 3, is
\begin{align*}
& U_X = D_X(U) = \frac{u_x}{(cu+d)^2},\\
& U_{XX} = D_X^2(U) = \frac{u_{xx}}{(cu+d)^2}-\frac{2cu_x^2}{(cu+d)^3},\\
& U_{XXX} = D_X^3(U) = \frac{u_{xxx}}{(cu+d)^2} - \frac{6cu_xu_{xx}}{(cu+d)^3}+\frac{6c^2u_x^3}{(cu+d)^4}.
\end{align*}
Assuming $u_x\neq 0$, we construct a moving frame by choosing the cross-section
\begin{equation}\label{sl2 differential cross-section}
\mathcal{K} = \{u=0,\; u_x=\epsilon = \text{sign}(u_x),\; u_{xx}=0\}.
\end{equation}
Solving the normalization equations
\[
U = \frac{au+b}{cu+d} = 0,\qquad U_X = \frac{u_x}{(cu+d)^2} = \epsilon,\qquad U_{XX} = \frac{u_{xx}}{(cu+d)^2}-\frac{2cu_x^2}{(cu+d)^3} = 0,
\]
for the group parameters and using the unitary constraint $ad-bc=1$, we obtain the right moving frame
\begin{equation}\label{sl2 differential moving frame}
a=\frac{1}{|u_x|^{1/2}},\qquad b=-\frac{u}{|u_x|^{1/2}},\qquad c=\frac{u_{xx}}{2|u_x|^{3/2}},\qquad d=\frac{2u_x^2-u u_{xx}}{2|u_x|^{3/2}}.
\end{equation}
Invariantizing $\epsilon u_{xxx}$, we obtain the Schwarzian derivative
\[
\epsilon\, \iota(u_{xxx}) = \frac{u_x u_{xxx} - (3/2)u_{xx}^2}{u_x^2}.
\]
\end{example}

\begin{exercise}\label{differential moving frame exercise}
Referring to Exercise \ref{invariants Lie's approach exercise}:
\begin{enumerate}
\item \label{differential moving frame exercise part a} Find the one-parameter group action induced by each of the infinitesimal generators~\eqref{burgers symmetry generators}.
\item \label{differential moving frame exercise part b} Construct a moving frame on $\J^{(1)} = \{(t,x,u,u_t,u_x)\}$.
\item \label{differential moving frame exercise part c} Compute the normalized invariant $\iota(u_{xx})$.
\end{enumerate}
\end{exercise}

\begin{example}\label{ex:sl2 invariants moving frames}
We now reconsider Example \ref{ex:sl2 invariants} using the method of moving frames.  The product action on $\mJ^{[3]}$ is
\begin{gather*}
X_{i-1}=x_{i-1},\qquad X_i = x_i,\qquad X_{i+1}=x_{i+1},\qquad X_{i+2} = x_{i+2},\\
U_{i-1} = \frac{au_{i-1}+b}{cu_{i-1}+d},\quad\; U_i = \frac{au_i+b}{cu_i+d},\quad\; 
U_{i+1} = \frac{au_{i+1}+b}{cu_{i+1}+d},\quad\; U_{i+2} = \frac{au_{i+2}+b}{cu_{i+2}+d}.
\end{gather*}
In the following, we let
\[
\epsilon_i = \text{sign}\bigg(\frac{u_{i+1}-u_{i-1}}{(u_i-u_{i-1})(u_{i+1}-u_i)} \bigg).
\]
Then, a cross-section to the group orbits is given by
\begin{equation}\label{sl2 discrete cross-section}
\mathcal{K} = \{ u_{i-1} = \epsilon_i,\; u_i \to \infty,\, u_{i+1} = 0\},
\end{equation}
where we let $u_i$ tend to infinity.  Solving the normalization equations
\[
U_{i-1} = \epsilon_i,\qquad U_i \to \infty,\qquad U_{i+1} = 0,
\]
we obtain the right moving frame
\[
a = - \frac{1}{c(u_{i+1}-u_i)},\qquad b = -\frac{u_{i+1}}{c(u_{i+1}-u_i)},\qquad d = -c u_i,
\]
where
\[
c = \pm \sqrt{\bigg|\frac{u_{i+1}-u_{i-1}}{(u_{i+1}-u_i)(u_i-u_{i-1})}\bigg|}.
\]
Invariantizing $\epsilon_i u_{i+2}$ we obtain the same difference invariant as in \eqref{R}:
\[
\epsilon_i\, \iota(u_{i+2}) = R_i = \frac{(u_{i+2}-u_{i+1})(u_i - u_{i-1})}{(u_{i+2}-u_i)(u_{i+1}-u_{i-1})}.
\]

The latter could also be derived from the replacement principle.  Invariantizing \eqref{R} we find that 
\begin{multline*}
R_i = \iota(R_i) = \frac{(\iota(u_i) - \iota(u_{i-1}))(\iota(u_{i+2})-\iota(u_{i+1}))}{(\iota(u_{i+1}-\iota(u_{i-1}))(\iota(u_{i+2})-\iota(u_i))} \\
=\frac{(\iota(u_i) - \epsilon_i)\iota(u_{i+2})}{-\epsilon_i(\iota(u_{i+2})-\iota(u_i))}\underset{\iota(u_i)\to \infty}{\longrightarrow} \epsilon_i\, \iota(u_{i+2}). 
\end{multline*}
\end{example}

\section{Weakly invariant equations}\label{weakly invariant equation section}

As observed in Remark \ref{invariant remark}, the notion of an invariant function is more restrictive than that of an invariant equation.  This brings us to distinguish two types of invariant equations.

\begin{definition}
An equation $F(z) = 0$ is said to be \emph{weakly invariant} if it is invariant only on its solution space.  That is
\[
F(g\cdot z) = 0\qquad \text{provided}\qquad F(z)=0
\]
and the action is defined.  An equation $F(z) = 0$ is said to be \emph{strongly invariant} if the function $F\colon \mM \to \mathbb{R}$ is $G$-invariant.  That is,
\[
F(g\cdot z) = F(z)\qquad \text{for all}\qquad g \in G
\]
where the action is defined.
\end{definition}

\begin{remark}
We note that a weakly invariant equation can, sometimes, be made strongly invariant by appropriately multiplying the equation by a certain \emph{relative invariant}.  We recall that a relative invariant of weight $\mu$ is a function $R(z)$ which satisfies $R(g\cdot z) = \mu(g,z)R(z)$. Indeed, if a weakly invariant equation $F(z)=0$ is such that $F(g\cdot z) = \mu(g,z) F(z)$, with $\mu(g,z) \neq 0$, then multiplying the equation by a relative invariant $R(z)\neq 0$ of weight $1/\mu$ yields the strongly invariant equation $R(z)F(z) =0$.

As a simple example, let $\mM = \mJ^{[1]} = \{(i,x_i,x_{i+1},u_i,u_{i+1})\}$, and consider the product action
\[
X_i = \lambda x_i+a,\qquad X_{i+1} = \lambda x_{i+1}+a,\qquad U_i = \lambda u_i + b,\qquad U_{i+1} = \lambda u_{i+1}+b,
\]
where $\lambda > 0$ and $a,b \in \mathbb{R}$.  Then the equation
\begin{equation}\label{weakly invariant equation exampleT}
u_{i+1} - u_i = 0
\end{equation}
is weakly invariant as $g\cdot u_{i+1} - g\cdot u_i = \lambda(u_{i+1}-u_i)$.  Dividing equation \eqref{weakly invariant equation exampleT} by the relative invariant $h_i = x_{i+1}-x_i$, one obtains the equivalent strongly invariant equation
\[
\frac{u_{i+1} - u_i}{x_{i+1}-x_i} = 0.
\]
\end{remark}

We now explain how to systematically search for weakly invariant equations.  As always, we assume that $G$ is an $r$-dimensional Lie group acting locally on a $d$-dimensional manifold $\mM$.

\subsection{Lie's infinitesimal approach}

A weakly invariant equation is found by searching for a submanifold $\mS\subset \mM$, defined as the zero locus of an equation $W(z)=0$, where the isotropy group is non-trivial.  To find such a submanifold we consider a basis of infinitesimal generators \eqref{g basis} and introduce the corresponding Lie matrix.

\begin{definition}
The \emph{Lie matrix} is the  $r\times d$ matrix whose components are given by the coefficients of the infinitesimal generators \eqref{g basis}:
\begin{equation}\label{Lie matrix}
\mathbf{L}(z) = 
\begin{bmatrix}
\zeta^1_1(z) & \cdots & \zeta^d_1(z) \\
\vdots & & \vdots \\
\zeta^1_r(z) & \cdots & \zeta^d_r(z)
\end{bmatrix}.
\end{equation}
\end{definition}

\begin{proposition}
The dimension of the group orbit through $z \in \mM$ is equal to the rank of the Lie matrix $\mathbf{L}(z)$.
\end{proposition}

\begin{proposition}\label{weakly invariant equation proposition}
Let $0\leq k \leq r$.  The set of points 
\[
\mS_k = \{ z \in \mM\,|\, \text{rank}\,\mathbf{L}(z) = k\}
\] 
is invariant under the action of $G$.  The number of functionally independent invariants on $\mS_k$ is given by the formula
\[
\text{dim}\,\mM - \text{rank}\, \mathbf{L}|_{\mS_k} = d - k.
\]
\end{proposition}

The sets of points $\mS_k$ where the rank of the Lie matrix $\mathbf{L}$ is not maximal, i.e. $\text{rank}\, \mathbf{L} < r$, are described by equations of the form $W(z)=0$.  By Proposition \ref{weakly invariant equation proposition}, these equations are weakly invariant.  Therefore, weakly invariant equations are found by searching for submanifolds where the rank of the Lie matrix is not maximal.

\begin{example}\label{weakly invariant equation example}
To illustrate the above considerations, we consider the Lie algebra of vector fields
\begin{equation}\label{5d infinitesimal generators}
\vv_1 = \pp{}{x},\qquad \vv_2 = \pp{}{u},\qquad \vv_3=x\pp{}{x},\qquad \vv_4=x\pp{}{u},\qquad \vv_5=u\pp{}{u},
\end{equation}
and search for weakly invariant equations on the discrete jet space
\begin{equation}\label{u_xx discrete jet space}
\mJ^{[2]} = \{(i,x_{i-1},x_i,x_{i+1},u_{i-1},u_i,u_{i+1})\}.
\end{equation}
The prolongation of the infinitesimal generators \eqref{5d infinitesimal generators} to $\mJ^{[2]}$ is given by
\begin{align*}
&\vv_1^{[2]} = \pp{}{x_{i-1}} + \pp{}{x_i} + \pp{}{x_{i+1}}, \\
&\vv_2^{[2]} = \pp{}{u_{i-1}} + \pp{}{u_i} + \pp{}{u_{i+1}}, \\
&\vv_3^{[2]} = x_{i-1}\pp{}{x_{i-1}} + x_i \pp{}{x_i} + x_{i+1}\pp{}{x_{i+1}}, \\
&\vv_4^{[2]} = x_{i-1}\pp{}{u_{i-1}} + x_i \pp{}{u_i} + x_{i+1} \pp{}{u_{i+1}},\\
&\vv_5^{[2]} = u_{i-1}\pp{}{u_{i-1}} + u_i \pp{}{u_i} + u_{i+1} \pp{}{u_{i+1}},
\end{align*}
and the corresponding Lie matrix is
\[
\mathbf{L} = 
\begin{bmatrix}
1 & 1 & 1 & 0 & 0 & 0 \\ 
0 & 0 & 0 & 1 & 1 & 1 \\
x_{i-1} & x_i & x_{i+1} & 0 & 0 & 0 \\
0 & 0 & 0 & x_{i-1} & x_i & x_{i+1} \\
0 & 0 & 0 & u_{i-1} & u_i & u_{i+1}
\end{bmatrix}.
\]
Assuming that $h_i = x_{i+1} - x_i \neq 0$ and $h_{i-1} = x_i - x_{i-1} \neq 0$, the Lie matrix can be row reduced to
\[
\mathbf{L} \sim \begin{bmatrix}
1 & 1 & 1 & 0 & 0 & 0 \\ 
0 & 0 & 0 & 1 & 1 & 1 \\
-h_{i-1} & 0 & h_i & 0 & 0 & 0 \\
0 & 0 & 0 & -h_{i-1} & 0 & h_i \\
0 & 0 & 0 & W & 0 & 0
\end{bmatrix},
\]
where $W = h_{i-1}(u_{i+1}-u_i)-h_i(u_i-u_{i-1})$.  Therefore, when
\begin{equation}\label{W=0}
h_{i-1}(u_{i+1}-u_i)-h_i(u_i-u_{i-1}) = 0
\end{equation}
the rank of the Lie matrix is not maximal and \eqref{W=0} yields a weakly invariant equation.
\end{example}

\subsection{Moving frame approach}

As explained in Section \ref{invariants - moving frame section}, a moving frame $\rho\colon \mM \to G$ exists provided the group action is free.  In terms of the Lie matrix \eqref{Lie matrix}, this occurs where the rank of $\mathbf{L}(z) = r = \text{dim}\,G$ is maximal.  Therefore, submanifolds where the rank of the Lie matrix is not maximal occur where a moving frame does not exist.  In those situations it is still possible to construct \emph{partial moving frames}, \cite{O-2011,OV-2015,V-2013}.  Intuitively, a partial moving frame is the $G$-equivariant map that one obtains when some of the group parameters cannot be normalized during the normalization procedure.  Given a partial moving frame, the invariantization map is still defined as in \eqref{invariantization}, and a complete set of normalized difference invariants can still be constructed.  

In applications, partial moving frames naturally occur as one attempts to solve the normalizing equations \eqref{normalization equations}.   The solution to the normalization equations will, in general, require some non-degeneracy conditions to hold and submanifolds where those constraints are not satisfied will determine weakly invariant equations.


\begin{example}\label{partial moving frame example}
We now reconsider Example \ref{weakly invariant equation example} using the equivariant moving frame method.  The group of transformations induced by the infinitesimal generators is given by
\[
X = \lambda x + a,\qquad U = \alpha u + \beta x + b,
\]
where $\lambda > 0$, $\alpha > 0$, and $a, b, \beta \in \mathbb{R}$.   The product action on the discrete jet space \eqref{u_xx discrete jet space} is
\begin{equation}\label{u_xx product action}
\begin{aligned}
&X_{i-1} = \lambda x_{i-1} + a,& & U_{i-1} = \alpha u_{i-1} + \beta x_{i-1} + b,\\
&X_i = \lambda x_i + a,& & U_i = \alpha u_i + \beta x_i + b,\\
&X_{i+1} = \lambda x_{i+1} + a, &\qquad & U_{i+1} = \alpha u_{i+1} + \beta x_{i+1} + b.
\end{aligned}
\end{equation}
Starting the normalization process, we first set $X_i=0$ and $U_i=0$.  Solving the normalization equations
\[
0 = X_i = \lambda x_i + a,\qquad 0 = U_i = \alpha u_i + \beta x_i + b,
\]
we obtain
\begin{equation}\label{ab normalization}
a = -\lambda x_i,\qquad b = -\alpha u_i - \beta x_i.
\end{equation}
Introducing the notation
\[
h_i = x_{i+1}-x_i,\qquad h_{i-1} = x_i-x_{i-1},\qquad \Delta u_i = u_{i+1}-u_i,\qquad \Delta u_{i-1} = u_i - u_{i-1},
\]
the substitution of the group normalizations \eqref{ab normalization} into the product action \eqref{u_xx product action} yields 
\begin{equation}\label{u_xx partially normalized prolonged action}
\begin{aligned}
&X_{i-1} = -\lambda\, h_{i-1}, & & U_{i-1} = -\alpha\, \Delta u_{i-1} - \beta\, h_{i-1},\\
&X_{i+1} = \lambda\, h_i,&\qquad & U_{i+1} = \alpha\, \Delta u_i + \beta\, h_i.
\end{aligned}
\end{equation}

At this stage, assuming that $h_{i-1} > 0$ (and similarly $h_i > 0$) we can set $X_{i-1} = -1$, which leads to the group normalization
\begin{equation}\label{lambda normalization}
\lambda = \frac{1}{h_{i-1}}.
\end{equation}
Substituting \eqref{lambda normalization} into \eqref{u_xx partially normalized prolonged action} yields the difference invariant
\[
H_i = \iota(x_{i+1}) = \frac{h_i}{h_{i-1}}.
\]
To normalize the remaining group parameters $\alpha$ and $\beta$ in
\[
\begin{bmatrix} U_{i-1} \\ U_{i+1} \end{bmatrix} = \begin{bmatrix} -\Delta u_{i-1} & -h_{i-1} \\ \Delta u_i & h_i \end{bmatrix}
\begin{bmatrix} \alpha \\ \beta \end{bmatrix},
\]
it is necessary for the coefficient matrix to be invertible.  On the other hand, if the matrix is not invertible, that is, if
\begin{equation}\label{zero determinant}
0 = \det \begin{bmatrix} -\Delta u_{i-1} & -h_{i-1} \\ \Delta u_i & h_i \end{bmatrix} = -h_i\, \Delta u_{i-1} + h_{i-1}\, \Delta u_i,
\end{equation}
one can only construct a partial moving frame and \eqref{zero determinant} is a weakly invariant equation, which is identical to \eqref{W=0}.  When \eqref{zero determinant} holds, we can normalize either $U_{i+1}$ or $U_{i-1}$.  Let $U_{i-1}=0$, then
\[
\beta = - \alpha \frac{\Delta u_{i-1}}{h_{i-1}},
\]
and one obtains the partial moving frame
\[
a = -\frac{x_i}{h_{i-1}},\qquad b = \frac{\alpha}{h_{i-1}}\det\begin{bmatrix} u_i & x_i \\ u_{i-1} & x_{i-1} \end{bmatrix},\qquad 
\lambda = \frac{1}{h_{i-1}},\qquad \beta = - \alpha \frac{\Delta u_{i-1}}{h_{i-1}}.
\]
Finally, we note that
\[
\iota(u_{i+1}) = \alpha\bigg[\Delta u_i - \frac{h_i}{h_{i-1}} \Delta u_{i-1}\bigg] = 0,
\]
by virtue of \eqref{zero determinant}.
\end{example}

\section{Symmetry-preserving numerical schemes}\label{symmetry-preserving schemes section}

At this point, given a Lie group of local transformations $G$ acting on $\mM$, we have everything needed to construct $G$-invariant equations.  As introduced in Section \ref{weakly invariant equation section}, a $G$-invariant equation $F(z)=0$ will either be weakly invariant or strongly invariant.  To obtain strongly invariant equations, the first step consists of computing a complete set of invariants $\mathbf{I}_c$ using either Lie's infinitesimal approach or the moving frame method. Once a complete set of invariants $\mathbf{I}_c$ has been computed, a strongly invariant equation $0=F(z)=\widetilde{F}(\mathbf{I}_c)$ is simply obtained by combining invariants from $\mathbf{I}_c$.  To obtain weakly invariant equations, one simply has to use one of the two procedures outlined in Section~\ref{weakly invariant equation section}.  

When $\mM = \J\n$, the above procedure will produce all the differential equations $\Delta(x,u\n)=0$ admitting $G$ as a symmetry group.  Similarly, when $\mM = \mJ\n$, one obtains all the differential equations $\overline{\Delta}(s,x\n,u\n)=0$ invariant under the prolonged action of $G$. Obtaining these differential equations is referred to as the \emph{inverse problem of group classification}. Given a $G$-invariant differential equation $\Delta(x,u\n)=0$, the procedure can also be used to construct an extended system of equations $\{\overline{\Delta}(s,x\n,u\n)=0,\widetilde{\Delta}(s,x\n,u\n)=0\}$ that is $G$-compatible with $\Delta(x,u\n)=0$.  

In the following, we are mainly interested in the case when $\mM = \mJ\sn$.  Given a differential equation $\Delta(x,u\n)=0$ with symmetry group $G$, we want to construct a system of finite difference equations that approximates the differential equation, specifies constraints on the mesh, and preserves the symmetry group $G$.  This is now obviously done by finding an appropriate collection of strongly invariant and weakly invariant difference equations, which, in the continuous limit, converge to the differential equation.  To find an approximation of the differential equation $\Delta(x,u\n)=0$, the first step consists of computing a complete set of difference invariants using either Lie's infinitesimal approach or the moving frame method and to consider their Taylor expansion.   Then one searches for a combination of these Taylor expansions that will, in the continuous limit, converge to the differential equation, and thereby provide a finite difference approximation of the equation.  This step will not always work.  It is possible that the difference invariants cannot be combined in such a way to converge to the differential equation in the continuous limit.  When this is the case, one should search for a weakly invariant equation $W(x_N\sn,u_N\sn)=0$ that converges to $\Delta=0$.  If the later fails, one can try to add more points in the lattice.  It is not clear yet if all invariant differential equations admit at least one symmetry-preserving scheme.  Differential equations with infinite-dimensional symmetry groups are particularly challenging.  To this day, no one has been able to systematically construct symmetry-preserving schemes for such equations.

In parallel, one also searches for a set of strongly and/or weakly invariant difference equations that will constrain the mesh on which the differential equation is approximated.  Mesh equations do not always have to be included. If they are avoided, this leads to \emph{invariant meshless discretization schemes}, \cite{B-2013}, which tend to be more complicated numerical schemes as they can operate on an arbitrary collection of nodes where the solution is sought numerically. When, included, the mesh equations will influence how the continuous limit should be taken in the above paragraph.  The number of equations specifying the mesh is not a priori fixed.  The only requirements are that those equations should not impose any constraints on the dependent variables $u_N\sn$, that they should be compatible, and have the appropriate continuous limit.  The latter can mean different things depending on the point of view used.  For example, in \cite{BCW-2006,BRW-2008,CRW-2016,DKW-2000,DKW-2003,DKW-2004,LMW-2015,LW-2006}, the discrete indices $N+K$ are assumed to be fixed and in the continuous limit, the points $(x_{N+K},u_{N+K})$ converge to $(x_N,u_N)$ for all $K \in \mathbb{Z}^p$.  With this perspective, the mesh equations will converge, in the continuous limit, to identities such as $0=0$. Alternatively, if one regards the discrete index $N=(n^1,\ldots,n^p)$ as sampling the computational variables $s=(s^1,\ldots,s^p)$, one can take the limit in the index variables, \cite{BM,RV-2013,RV-2015}. To this end, we introduce the variation parameters $\epsilon = (\epsilon^1,\ldots,\epsilon^p)\in [0,1]^p$.  One can then write the multi-index $N+K$ as
\[
N+K = N+ \epsilon\cdot K|_{\epsilon = (1,\ldots,1)}=(n^1+\epsilon^1k^1,\ldots,n^p+\epsilon^p k^p)|_{\epsilon = (1,\ldots,1)}.
\]
Letting $\epsilon \to 0^+$, one has that
\[
\lim_{\epsilon \to 0^+} N+ \epsilon\cdot K = N.
\]
By introducing the variation parameters $\epsilon=(\epsilon^1,\ldots,\epsilon^p)$  in the mesh equations and letting $\epsilon \to 0^+$, the latter will now converge to the companion equations \eqref{companion equations}.  

We now illustrate the above procedure for constructing symmetry-preserving numerical schemes by considering three examples.

\begin{example}
As our first example, we construct a symmetry-preserving scheme for the Schwarzian differential equation \eqref{Schwarzian invariant scheme}, whose symmetry group is given by the fractional linear action \eqref{sl2 action}.  The difference invariants on the lattice \eqref{sl2 lattice} are given by the index $i$, the discrete $x$-variables \eqref{sl2 x invariants}, and the cross-ratio \eqref{R}.  These invariants are sufficient to construct a symmetry-preserving scheme of the Schwarzian equation.  We begin by specifying the mesh equation, as this step is the easiest.  Clearly, we can set
\[
x_{i+1}-x_i = h,
\]
where $h>0$ is a positive constant.  From the mesh equation, it follows that 
\[
x_{i-1} = x_i - h,\qquad x_{i+1} = x_i + h,\qquad x_{i+2} = x_i + 2h.
\]
Therefore, the Taylor expansions of $u_{i-1}$, $u_{i+1}$, and $u_{i+2}$ centered at $x_i$ are
\begin{equation}\label{Taylor expansions}
\begin{aligned}
u_{i-1} &= u(x_{i-1}) = u - h u_x + \frac{h^2}{2} u_{xx} - \frac{h^3}{6} u_{xxx} + \mO(h^4),\\
u_{i+1} &= u(x_{i+1}) = u + h u_x + \frac{h^2}{2} u_{xx} + \frac{h^3}{6} u_{xxx} + \mO(h^4),\\
u_{i+2} &= u(x_{i+2}) = u + 2h u_x + 2h^2 u_{xx} + \frac{4h^3}{3} u_{xxx} + \mO(h^4),
\end{aligned}
\end{equation}
where the function $u$ and its derivatives are evaluated at $x_i$.  Substituting the Taylor expansions \eqref{Taylor expansions} in the difference invariant \eqref{R}, we obtain
\[
R_i = \frac{h^2 u_x^2 + h^3 u_x u_{xx} + h^4(\frac{4}{3} u_x u_{xxx} - \frac{3}{4} u_{xx}^2) + \mO(h^5)}{4 h^2 u_x^2 + 4h^3 u_x u_{xx} + \frac{10}{3} h^4 u_x u_{xxx} + \mO(h^5)}.
\]
Therefore,
\[
 \frac{4-1/R_i}{h^2} = \frac{2 u_x u_{xxx} - 3 u_{xx}^2}{u_x^2} + \mO(h)
\]
and an invariant approximation of the Schwarzian equation \eqref{Schwarzian equation} is given by
\begin{equation}\label{Schwarzian invariant scheme}
\frac{1}{h^2}\bigg[2 - \frac{1}{2R_i}\bigg] = F(x_i).
\end{equation}
\end{example}

\begin{example}
As a second example, we consider the second order ordinary differential equation
\begin{equation}\label{u_xx=0}
u_{xx} = 0.
\end{equation}
As seen in Exercise \ref{u_xx=0 exercise}, the infinitesimal symmetry algebra is spanned by the vector fields \eqref{u_xx=0 symmetry generators}. In the following, we construct a symmetry-preserving scheme of the differential equation \eqref{u_xx=0} invariant under the five-dimensional symmetry subgroup generated by \eqref{5d infinitesimal generators} on the discrete jet space $\mJ^{[2]} = \{(i,x_{i-1},x_i,x_{i+1},u_{i-1},u_i,u_{i+1})\}.$  Since $\dim \mJ^{[2]} - \dim \g = 7 - 5 = 2$, other than the index $i$, we expect one more difference invariant.  Solving the system of first order partial differential equations 
\begin{align*}
\vv_1^{[2]}(I) &= \pp{I}{x_{i-1}} + \pp{I}{x_i} + \pp{I}{x_{i+1}} = 0, \\
\vv_2^{[2]}(I) &= \pp{I}{u_{i-1}} + \pp{I}{u_i} + \pp{I}{u_{i+1}} = 0, \\
\vv_3^{[2]}(I) &= x_{i-1}\pp{I}{x_{i-1}} + x_i \pp{I}{x_i} + x_{i+1}\pp{I}{x_{i+1}} = 0, \\
\vv_4^{[2]}(I) &= x_{i-1}\pp{I}{u_{i-1}} + x_i \pp{I}{u_i} + x_{i+1} \pp{I}{u_{i+1}} = 0, \\
\vv_5^{[2]}(I) &= u_{i-1}\pp{I}{u_{i-1}} + u_i \pp{I}{u_i} + u_{i+1} \pp{I}{u_{i+1}} = 0,
\end{align*}
we obtain the invariant 
\begin{equation}\label{H_n}
H_i = \frac{x_{i+1}-x_i}{x_i-x_{i-1}} = \frac{h_i}{h_{i-1}}.
\end{equation}
We note that this invariant was also found in our construction of a partial moving frame in Example \ref{partial moving frame example}.  Clearly, it is not possible to approximate \eqref{u_xx=0} using only the invariant \eqref{H_n}. We therefore search for weakly invariant equations.  In Example \ref{weakly invariant equation example} (and Example \ref{partial moving frame example}) we found that $W = h_{i-1}(u_{i+1}-u_i)-h_i(u_i-u_{i-1})=0$ is a weakly invariant equation.  Since the product of a weakly invariant equation $W(x_i^{[2]},u_i^{[2]})=0$ by a nonzero difference function $F(i,x_i^{[2]},u_i^{[2]}) \neq 0$ remains weakly invariant,
\[
\frac{2W}{h_i h_{i-1}(h_i + h_{i-1})} = \frac{2}{x_{i+1}-x_{i-1}}\bigg(\frac{u_{i+1}-u_i}{x_{i+1}-x_i} - \frac{u_i - u_{i-1}}{x_i-x_{i-1}} \bigg)=0
\]
is weakly invariant equation, and happens to approximate the differential equation $u_{xx}=0$.  As for the mesh equation, we set $H_i = f(i)$, with $f(i) > 0$ for all $i$, and obtain
\[
x_{i+1}-(1+f(i))x_i + f(i) x_{i-1} = 0.
\]
\end{example}

\begin{exercise}\label{order 1 ode exercise}
(This exercise was taken from \cite{RW-2004}.) The first order ordinary differential equation
\begin{equation}\label{order 1 ode}
u^\prime = A^\prime(x) u + B^\prime(x) e^{A(x)}
\end{equation}
is invariant under the infinitesimal symmetry generators
\[
\vv_1 = e^{A(x)}\pp{}{u},\qquad \vv_2 = [u-B(x) e^{A(x)}]\pp{}{u}.
\]
Working on the discrete jet space
\[
\mJ^{[1]} = \{(i,x_i,x_{i+1},u_i,u_{i+1})\}:
\]
\begin{enumerate}
\item Show that, other than the index $i$, the only two invariants are $x_i$ and $x_{i+1}$.
\item \label{order 1 ode exercise part b} Find a weakly invariant difference equation.
\item Write down a symmetry-preserving scheme for \eqref{order 1 ode}.
\end{enumerate}
\end{exercise}

\begin{example}\label{KdV invariant scheme}
The standard discretization of the KdV equation on an orthogonal mesh given in \eqref{naive KdV discretization} is not invariant under the Galilean boosts
\[
X=x+vt,\qquad T = t,\qquad U = u + v,\qquad v\in\mathbb{R}.
\]
Indeed, under this transformation, the second term in \eqref{naive KdV discretization} is transformed to 
\[
(u_i^n + v)\cdot \frac{u^n_{i+1} - u^n_{i-1}}{2h},
\]
while the other two terms remain unchanged.  Thus, the discretization \eqref{naive KdV discretization} does not preserve all the symmetries of the equation.  It is not difficult to see that the discretization \eqref{naive KdV discretization} is only invariant under shifts and dilations.  

We now proceed to the construction of a symmetry-preserving scheme for the KdV equation.  Introducing the multi-index $N=(n,i)$, let
\[
(t_N,x_N,u_N) = (t^n_i,x^n_i,u^n_i)
\]
as in Example \ref{KdV standard scheme example}.  Recall that the symmetry generators of the KdV equation were found in \eqref{KdV symmetry generators}.  Clearly, the ratios
\[
\frac{t^n_{i+1} - t^n_i}{t^{n+1}_i-t^n_i}\qquad\text{and}\qquad \frac{t^{n+1}_i - t^n_i}{t^n_i-t^{n-1}_i}
\]
are invariant under space and time translations, Galilean boosts, and scalings.  Therefore, we can use these two invariants to fix invariant constraints on the discretization of the time variable $t$ by setting $\frac{t^n_{i+1} - t^n_i}{t^{n+1}_i-t^n_i}=0$ and $\frac{t^{n+1}_i - t^n_i}{t^n_i-t^{n-1}_i}=1$. These two equations are equivalent to 
\begin{equation}\label{KdV t mesh equations}
t^n_{i+1} - t^n_i = 0,\qquad t^{n+1}_i - 2t^n_i + t^{n-1}_i=0.
\end{equation}
The latter imply that a symmetry-preserving scheme for the KdV equation can be formulated on a mesh with flat, equally spaced time layers:
\[
t^n = k\, n + t^0.
\]
From here on, we assume that \eqref{KdV t mesh equations} hold.  The prolongation of the vector fields \eqref{KdV symmetry generators} to the points in the stencil depicted in Figure \ref{KdV stencil} is given by

\begin{figure}[!ht]
\begin{center}
\tikzstyle{block} = [draw, 
rectangle, minimum height=10pt, minimum width=10pt]
\begin{tikzpicture}
\clip (-1,0) rectangle (14, 3.25);
\draw (0,1) -- (13,1);
\draw(1,1) node[] {$\bullet$};
\draw(1,0.5) node {$(t^n,x^n_{i-2},u^n_{i-2})$};
\draw(3.5,1) node[] {$\bullet$};
\draw(3.5,0.5) node {$(t^n,x^n_{i-1},u^n_{i-1})$};
\draw(6,1) node[] {$\bullet$};
\draw(6,0.5) node {$(t^n,x^n_i,u^n_i)$};
\draw(9,1) node[] {$\bullet$};
\draw(9,0.5) node {$(t^n,x^n_{i+1},u^n_{i+1})$};
\draw(12,1) node[] {$\bullet$};
\draw(12,0.5) node {$(t^n,x^n_{i+2},u^n_{i+2})$};
\draw (0,3) -- (13,3);
\draw(0.5,3) node[] {$\bullet$};
\draw(0.5,2.5) node {$(t^{n+1},x^{n+1}_{i-2},u^{n+1}_{i-2})$};
\draw(3.75,3) node[] {$\bullet$};
\draw(3.5,2.5) node {$(t^{n+1},x^{n+1}_{i-1},u^{n+1}_{i-1})$};
\draw(6.5,3) node[] {$\bullet$};
\draw(6.5,2.5) node {$(t^{n+1},x^{n+1}_i,u^{n+1}_i)$};
\draw(9.5,3) node[] {$\bullet$};
\draw(9.5,2.5) node {$(t^{n+1},x^{n+1}_{i+1},u^{n+1}_{i+1})$};
\draw(12.5,3) node[] {$\bullet$};
\draw(12.5,2.5) node {$(t^{n+1},x^{n+1}_{i+2},u^{n+1}_{i+2})$};

\end{tikzpicture}
\caption{\footnotesize{Stencil for the KdV equation.}}\label{KdV stencil}
\end{center}
\end{figure}
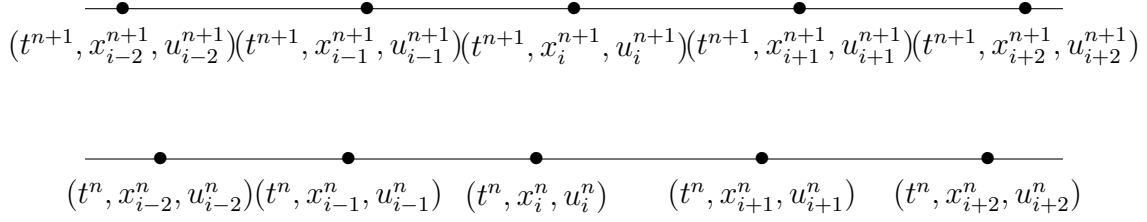

\begin{align}
&\vv_1 = \sum_{l=0}^1 \pp{}{t^{n+l}},\qquad
\vv_2 = \sum_{l=0}^1 \sum_{j=-2}^2  \pp{}{x^{n+l}_{i+j}},\qquad
\vv_3 = \sum_{l=0}^1 \sum_{j=-2}^2  t^{n+l}\pp{}{x^{n+l}_{i+j}} + \pp{}{u^{n+l}_{i+j}},\nonumber\\
&\vv_4 = \sum_{l=0}^1 \bigg[3 t^{n+l} \pp{}{t^{n+l}} + \sum_{j=-2}^2  x^{n+l}_{i+j} \pp{}{x^{n+l}_{i+j}} - 2u^{n+l}_{i+j}\pp{}{u^{n+l}_{i+j}}\bigg].\label{KdV notation}
\end{align}
To simplify the notation, we introduce
\begin{equation}\label{notation}
k = t^{n+1}-t^n,\qquad h^n_i = x^n_{i+1} - x^n_i,\qquad \sigma^n_i = x^{n+1}_i - x^n_i,\qquad Du^n_i = \frac{u^n_{i+1}-u^n_i}{h^n_i},
\end{equation}
for the spacings and elementary first order discrete $x$-derivatives.  Applying the infinitesimal invariance criterion \eqref{basis symmetry criterion} and solving the corresponding system of first order partial differential equations, we obtain the following 18 functionally independent invariants
\begin{align}\label{eq:DifferenceInvariantsKdV}
\begin{split}
&H^{n+l}_{i+j} = \frac{h^{n+l}_{i+j-1}}{h^{n+l}_{i+j}},\qquad l=0,1,\quad j=-1,0,1,\\
&I^n_i = \frac{h^{n+1}_i}{h^n_i},\qquad J^n_i = \frac{(h^n_i)^3}{k},\qquad L^n_i = \frac{\sigma^n_i - k\cdot u^n_i}{h^n_i},
\qquad T^n_i = (u^{n+1}_i - u^n_i)(h^n_i)^2,\\
&K^{n+l}_{i+j} = k\cdot Du^{n+l}_{i+j},\qquad l=0,1,\quad j=-2,-1,0,1.
\end{split}
\end{align}
Introducing the invariant quantity
\[
Q^n_i = H^n_{i+1}\bigg(\frac{K^n_{i+1}-K^n_i}{1+H^n_{i+1}}\bigg) - \bigg(\frac{K^n_i - K^n_{i-1}}{1+H^n_i}\bigg),
\]
an invariant numerical scheme for the KdV equation (together with the mesh equations \eqref{KdV t mesh equations}) is given by
\begin{equation}\label{6 point invariant scheme}
T^n_i - J^n_i\cdot L^n_i \bigg(\frac{K^n_i + K^n_{i-1}}{2}\bigg) + Q^n_i + \frac{Q^n_{i-1}}{(H^n_i)^2}=0.
\end{equation}
Introducing the third order discrete $x$-derivative
\begin{equation}\label{D^3u}
D^3 u^n_i = \frac{2}{h^n_i}\bigg[\bigg(\frac{Du^n_{i+1}-Du^n_i}{h^n_{i+1}+h^n_i}\bigg) - \bigg(\frac{Du^n_i - Du^n_{i-1}}{h^n_i + h^n_{i-1}}\bigg) \bigg],
\end{equation}
the explicit expression of the invariant scheme \eqref{6 point invariant scheme} is
\begin{equation}\label{coordinate 6 point invariant scheme}
\frac{u^{n+1}_i - u^n_i}{k} + \bigg(u^n_i - \frac{\sigma^n_i}{k}\bigg)\frac{Du^n_i + Du^n_{i-1}}{2} + \frac{1}{2}\big[D^3u^n_i + D^3u^n_{i-1}\big]=0.
\end{equation}

A more appropriate invariant numerical scheme can be realized on the entire ten point lattice.  The latter is given by
\begin{multline*}
T^n_i - J^n_i\cdot L^n_i \bigg(\frac{K^n_i+K^n_{i-1}+K^{n+1}_i+K^{n+1}_{i-1}}{4}\bigg) \\
+ \frac{1}{2}\bigg[ \frac{1}{I^n_i}\bigg (Q^{n+1}_i + \frac{Q^n_{i-1}}{I^n_i (H^{n+1}_i)^2}\bigg) + Q^n_i + \frac{Q^n_{i-1}}{(H^n_i)^2} \bigg] = 0.
\end{multline*}
Explicitly,
\begin{multline}\label{coordinate 10 point invariant scheme}
\frac{u^{n+1}_i - u^n_i}{k} + \bigg(u^n_i - \frac{\sigma^n_i}{k}\bigg)\frac{Du^n_i + Du^n_{i-1} + Du^{n+1}_i+Du^{n+1}_{i-1}}{4} \\
+ \frac{1}{4}\big[D^3u^{n+1}_i+D^3u^{n+1}_{i-1} + D^3u^n_i + D^3 u^n_{i-1} \big]=0.
\end{multline}
To use the scheme \eqref{coordinate 6 point invariant scheme} or \eqref{coordinate 10 point invariant scheme}, the \emph{grid velocity} 
\[
\frac{\sigma^n_i}{k} = \frac{x^{n+1}_i - x^n_i}{k}
\]
must be specified in an invariant manner to preserve the symmetries of the KdV equation.  One possibility is to set 
\begin{equation}\label{KdV x mesh equation}
L^n_i = 0\qquad \text{so that}\qquad \frac{\sigma^n_i}{k} = u^n_i.
\end{equation}
Together, the equations \eqref{KdV t mesh equations}, \eqref{coordinate 6 point invariant scheme} (or \eqref{coordinate 10 point invariant scheme}), and \eqref{KdV x mesh equation} provide a numerical approximation of the extended system of differential equations \eqref{KdV t companion equations}, \eqref{KdV x companion equation}, \eqref{KdV Lagrangian coordinates} for the KdV equation.  The latter scheme can perform poorly as there is no built-in mechanism preventing the clustering of grid points as the numerical integration proceeds. Alternatives to using \eqref{KdV x mesh equation} to obtain the position of the grid points at the next time level will be presented in Section~\ref{numerical simulations section}.  Using \emph{adaptive moving mesh methods}, we will construct more reliable invariant numerical schemes.


We conclude this example by discussing the continuous limit of the numerical scheme \eqref{coordinate 6 point invariant scheme} with mesh equations \eqref{KdV t mesh equations}, \eqref{KdV x mesh equation}.   Let us introduce the variation parameters $(\epsilon,\delta)$ so that
\[
n+l = n+l \epsilon\big|_{\epsilon=1}\qquad \text{and}\qquad i + j = i+j \delta\big|_{\delta = 1}.
\]
In the numerical scheme \eqref{coordinate 6 point invariant scheme}, \eqref{KdV t mesh equations}, \eqref{KdV x mesh equation}, let
\begin{gather*}
u^{n+1}_i - u^n_i=\frac{u^{n+\epsilon}_i - u^n_i}{\epsilon},\qquad u^n_{i+1} - u^n_i = \frac{u^n_{i+\delta} - u^n_i}{\delta},\\
u^n_{i+1} - 2u^n_i + u^n_{i-1}= \frac{u^n_{i+\delta} - 2 u^n_i + u^n_{i-\delta}}{\delta^2},\\
u^n_{i-2} - 3u^n_{i+1} + 3u^n_i - u^n_{i-1} = \frac{u^n_{i+2\delta} - 3u^n_{i+\delta} + 3u^n_i - u^n_{i-\delta}}{\delta^3}
\end{gather*}
with $\epsilon = 1$ and $\delta = 1$, and similarly for the differences in $t$ and $x$.  Then, as $\epsilon \to 0$ and $\delta \to 0$,
\begin{gather*}
\frac{u^{n+\epsilon}_i - u^n_i}{\epsilon} \to u_\tau,\qquad \frac{u^n_{i+\delta} - u^n_i}{\delta} \to u_s,\\
\frac{u^n_{i+\delta} - 2 u^n_i + u^n_{i-\delta}}{\delta^2} \to u_{ss} ,\qquad
\frac{u^n_{i+2\delta} - 3u^n_{i+\delta} + 3u^n_i - u^n_{i-\delta}}{\delta^3} \to u_{sss},
\end{gather*}
and the numerical scheme \eqref{coordinate 6 point invariant scheme}, \eqref{KdV t mesh equations}, \eqref{KdV x mesh equation}, converges to \eqref{KdV Lagrangian coordinates}, \eqref{KdV t companion equations}, and \eqref{KdV x companion equation}, respectively.

Alternatively, if one lets $t^{n+l}_{i+j}\to t^n_i$, $x^{n+l}_{i+j}\to x^n_i$, $u^{n+l}_{i+j}\to u^n_i$, without introducing the variation parameters $(\epsilon,\delta)$, then after multiplying equation \eqref{KdV x mesh equation} by $k$, the mesh equations \eqref{KdV t mesh equations}, \eqref{KdV x mesh equation} converge to the identity $0=0$, while \eqref{coordinate 6 point invariant scheme} converges to the original KdV equation \eqref{KdV}.
\end{example}

\begin{exercise}\label{Lie symmetry-preserving scheme exercise}
Using the difference invariants computed in Exercise \ref{invariants Lie's approach exercise} part (\ref{invariants Lie's approach exercise part b}), construct a symmetry-preserving scheme for Burgers' equation \eqref{burgers equation}.
\end{exercise}

When using the method of equivariant moving frames to construct symmetry-preserving numerical schemes, it is possible to avoid the step where one has to search for a combination of the difference invariants that will approximate the differential equation $\Delta(x,u\n)=0$. In general, for this to be the case, some care needs to be taken when constructing the discrete moving frame $\rho\sn\colon \mJ\sn \to G$.  The latter has to be \emph{compatible} with a continuous moving frame $\rho\n\colon \J\n \to G$.  By this, we mean that the discrete moving frame $\rho\sn$ must converge, in the continuous limit, to the moving frame $\rho\n$.   If $\mK\sn$ is the cross-section used to define $\rho\sn$ and $\mK\n$ is the cross-section defining $\rho\n$, then the moving frame $\rho\sn$ will be compatible with $\rho\n$ if, in the coalescing limit, $\mK\sn$ converges to $\mK\n$.  As shown in \cite{BM}, discrete compatible moving frames can be constructed by using the Lagrange interpolation coordinates on $\mJ\sn$, although in applications these can lead to complicated expressions that may limit the scope of the method.  It is frequently preferable to fix invariant constraints on the mesh, and then consider finite difference approximations of the derivatives compatible with the mesh.  On a nonuniform mesh, these expressions can be obtained by following the procedure of Example \ref{derivative discretization example}.

Given a compatible moving frame $\rho\sn$ an invariant approximation of the differential equation $\Delta(x,u\n)=0$ is simply obtained by invariantizing any finite difference approximation $E(N,x_N\sn,u_N\sn)=0$, compatible with the mesh.  In particular, the equation $E(N,x_N\sn,u_N\sn)=0$ does not have to be invariant.  

We now illustrate the construction of symmetry-preserving numerical schemes using the method of moving frames. 

\begin{example}
As our first example, we revisit Example \ref{Schwarzian invariant scheme} using the moving frame machinery.  In Example \ref{ex:sl2 invariants moving frames}, we constructed a discrete moving frame for the symmetry group of the Schwarzian differential equation.  This moving frame was constructed using the cross-section \eqref{sl2 discrete cross-section}, which is not compatible with the cross-section \eqref{sl2 differential cross-section} used to define a moving frame in the differential case.  Indeed, in \eqref{sl2 differential cross-section} we have $u=0$, while in the discrete case we let $u_i\to \infty$.  Therefore,  one should not expect that the invariantization of the standard scheme
\begin{equation}\label{Schwarzian standard scheme}
\frac{u_{i+2}-3u_{i+1}+3u_i - u_{i-1}}{(u_{i+1}-u_i)h^2} - \frac{3}{2} \bigg(\frac{u_{i+1}-2u_i + u_{i-1}}{(u_{i+1}-u_i)h}\bigg)^2 = F(x_i)
\end{equation}
will provide an invariant approximation of the Schwarzian equation \eqref{Schwarzian equation}.  Indeed, the invariantization of \eqref{Schwarzian standard scheme}, yields the inconsistent equation
\[
-\frac{9}{h^2} = F(x_i).
\]
In this case one has to combine the normalized invariants $x_{i-1}$, $x_i$, $x_{i+1}$, $x_{i+2}$, and the cross-ratio $\epsilon_i \iota(u_{i+2})=R_i$ as in Example \ref{Schwarzian invariant scheme} to obtain the invariant numerical scheme \eqref{Schwarzian invariant scheme}.  

For the invariantization of \eqref{Schwarzian standard scheme} to give a meaningful invariant discretization, we need to construct a discrete moving frame compatible with \eqref{sl2 differential moving frame}.  Working on the uniform mesh
\[
x_{i+1} - x_i = h,
\]
we introduce the finite difference derivatives
\begin{gather*}
Du_i = \frac{u_{i+1}-u_i}{h},\qquad D^2u_i = \frac{u_{i+1}-2u_i + u_{i-1}}{h^2},\\
D^3u_i = \frac{u_{i+2}-3u_{i+1}+3u_i - u_{i-1}}{h^3}.
\end{gather*}
To obtain a compatible discrete moving frame, we consider a finite difference approximation of the cross-section \eqref{sl2 differential cross-section} used in the differential case.  Namely, let
\[
\mK = \{u_i=0,\; Du_i = \epsilon_i,\; D^2u_i = 0\}
\]
where
\[
\epsilon_i = \text{sign}(u_{i+1}-u_i)(u_i-u_{i-1})(u_{i+2}-u_{i-1}). 
\]
The latter is equivalent to
\[
\mK = \{u_{i-1} = -h\epsilon_i,\; u_i = 0,\; u_{i+1} = h \epsilon_i\}.
\]
Solving the normalization equations
\[
-h\epsilon_i = U_{i-1} = \frac{au_{i-1}+b}{cu_{i-1}+b},\qquad
0=U_i = \frac{au_i+b}{cu_i+d},\qquad 
h \epsilon_i=U_{i+1} = \frac{au_{i+1}+b}{cu_{i+1}+b},
\]
and using the unitary constraint $ad-bc=1$, we obtain the moving frame
\[
a = \frac{D^2u_i}{2c\,Du_i\, Du_{i-1}},\qquad 
b=-\frac{u_i D^2u_i}{2c\,Du_i\,Du_{i-1}},\qquad
d= c\cdot \frac{u_{i+1}\,Du_{i-1}-u_{i-1}\,Du_i}{Du_i-Du_{i-1}},
\]
where
\[
c^2= \frac{\epsilon_i (D^2u_i)^2}{2\,Du_i\, Du_{i-1}(Du_i+Du_{i-1})}.
\]
We note that the right-hand side of the last equality is nonnegative by definition of $\epsilon_i$. Invariantizing the non-invariant scheme \eqref{Schwarzian standard scheme}, we obtain the invariant discretization
\begin{equation}\label{Schwarzian invariantized scheme}
\frac{(u_{i+2}-u_i)(u_{i+1}-u_{i-1})}{h^3[(u_{i+2}-u_{i-1})Du_i-(u_{i+2}-u_{i+1})Du_{i-1}]}-\frac{2}{h^2} = F(x_i).
\end{equation}
The latter can be written using cross-ratios in a form similar to the invariant scheme \eqref{Schwarzian invariant scheme}.  After some simplifications, the invariant scheme \eqref{Schwarzian invariantized scheme} is equivalent to
\[
\frac{1}{h^2}\bigg[\frac{1}{\overline{R}_i-R_i} - 2\bigg] = F(x_i),
\]
where
\[
\overline{R}_i = \frac{(u_{i+2}-u_{i-1})(u_{i+1}-u_i)}{(u_{i+2}-u_i)(u_{i+1}-u_{i-1})}\qquad \text{and}\qquad
R_i = \frac{(u_{i+2}-u_{i+1})(u_i-u_{i-1})}{(u_{i+2}-u_i)(u_{i+1}-u_{i-1})}.
\]
\end{example}

\begin{example}
As our final example, we consider the invariant discretization of the KdV equation.  The group action induced by the infinitesimal generators \eqref{KdV symmetry generators} is given by
\begin{equation}\label{eq:SymmetryGroupKdVEquation}
X = \lambda x + vt + a,\qquad T = \lambda^3 t + b,\qquad U = \frac{u}{\lambda^2} + v,
\end{equation}
where $a, b, v \in \mathbb{R}$ and $\lambda \in \mathbb{R}^+$.  As in Example \ref{KdV invariant scheme}, to simplify the computations, we assume that \eqref{KdV t mesh equations} holds.   When this is the case, it follows from \eqref{first order derivative approximations} that 
\[
u_t \approx \frac{u^{n+1}_i-u^n_i}{k} - \frac{\sigma^n_i}{k}Du^n_i,\qquad u_x \approx Du^n_i,
\]
where we use the notation that was introduced in \eqref{notation}.  For better numerical accuracy, we let
\begin{equation}\label{u_t and u_x approximations}
u_t \approx \Delta_tu^n_i =\frac{u^{n+1}_i-u^n_i}{k} - \frac{\sigma^n_i}{k}\cdot \frac{Du^n_i+Du^n_{i-1}}{2},\qquad u_x \approx \Delta_xu^n_i = \frac{Du^n_i+Du_{i-1}^n}{2}.
\end{equation}
Also, recalling formula \eqref{D^3u}, we let 
\begin{equation}\label{u_xxx approximation}
u_{xxx} \approx \Delta_x^3 u^n_i = \frac{1}{2}[D^3 u^n_i + D^3 u^n_{i-1}].
\end{equation}

Implementing the discrete moving frame construction, we choose the cross-section
\[
\mathcal{K} = \{x^n_i=0,\; t^n=0,\; u^n_i = 0,\; \Delta_x u^n_i=1\},
\]
which is compatible with the cross-section $\{x=0,\, t=0,\, u=0,\, u_x=1\}$ that one could use to construct a moving frame in the continuous case.  Solving the normalization equations
\begin{gather*}
\lambda x^n_i + vt^n + a = 0,\qquad \lambda^3 t^n + b=0,\qquad \lambda^{-2}u^n_i + v = 0, \\
\frac{\lambda^{-2}u^n_{i+1}+v}{\lambda x^n_{i+1}+vt^n+a}+\frac{\lambda^{-2}u^n_{i-1}+v}{\lambda x^n_{i-1}+vt^n+a}=2,
\end{gather*}
for the group parameters $a$, $b$, $v$, $\lambda$, we obtain the right moving frame
\begin{equation}\label{KdV moving frame}
\begin{aligned}
&a = -x^n_i (\Delta_x u^n_i)^{1/3} + \frac{t^n u^n_i}{(\Delta_x u^n_i)^{2/3}},\qquad b=-t^n\, \Delta_x u^n_i,\\
&v = -\frac{u^n_i}{(\Delta_x u^n_i)^{2/3}},\qquad \lambda=(\Delta_x u^n_i)^{1/3}.
\end{aligned}
\end{equation}
To obtain an invariant scheme, we approximate the KdV equation using \eqref{u_t and u_x approximations} and \eqref{u_xxx approximation},
\begin{equation}\label{coordinate 6 point invariant scheme 2}
\Delta_t u^n_i + u^n_i \cdot \Delta_xu^n_i + \Delta_x^3 u^n_i = 0.
\end{equation}
and invariantize the resulting scheme.  Since the latter is already invariant, the scheme remains the same.   We note that the scheme \eqref{coordinate 6 point invariant scheme 2} is the same as \eqref{coordinate 6 point invariant scheme}.
\end{example}

\begin{exercise}\label{difference moving frame exercise}
Referring to Exercise \ref{invariants Lie's approach exercise}:
\begin{enumerate}
\item \label{difference moving frame exercise part a}  Construct a discrete moving frame on the stencil
\[
\{(n,i,t^n,t^{n+1},x^n_{i-1},x^n_i,x^n_{i+1},x^{n+1}_{i-1},x^{n+1}_i,x^{n+1}_{i+1},u^n_{i-1},u^n_i,u^n_{i+1},u^{n+1}_{i-1},u^{n+1}_i,u^{n+1}_{i+1})\}
\]
compatible with the differential moving  frame found in Exercise \ref{differential moving frame exercise} part (\ref{differential moving frame exercise part b}).
\item \label{difference moving frame exercise part b} Invariantize the discrete approximation
\[
u_{xx} \approx D^2u^n_i = \frac{2}{h^n_i + h^n_{i-1}}[Du^n_i - Du^n_{i-1}].
\]
\item \label{difference moving frame exercise part c} Write a symmetry-preserving scheme for Burgers' equation \eqref{burgers equation}.
\end{enumerate}
\end{exercise}

\section{Numerical simulations}\label{numerical simulations section}

In this section we present some numerical simulations using the invariant numerical schemes derived in Section \ref{symmetry-preserving schemes section}.

\subsection{Schwarzian ODE}

\begin{figure}[!ht]
\centering
\captionsetup{width=0.85\textwidth}
\includegraphics[scale=0.6]{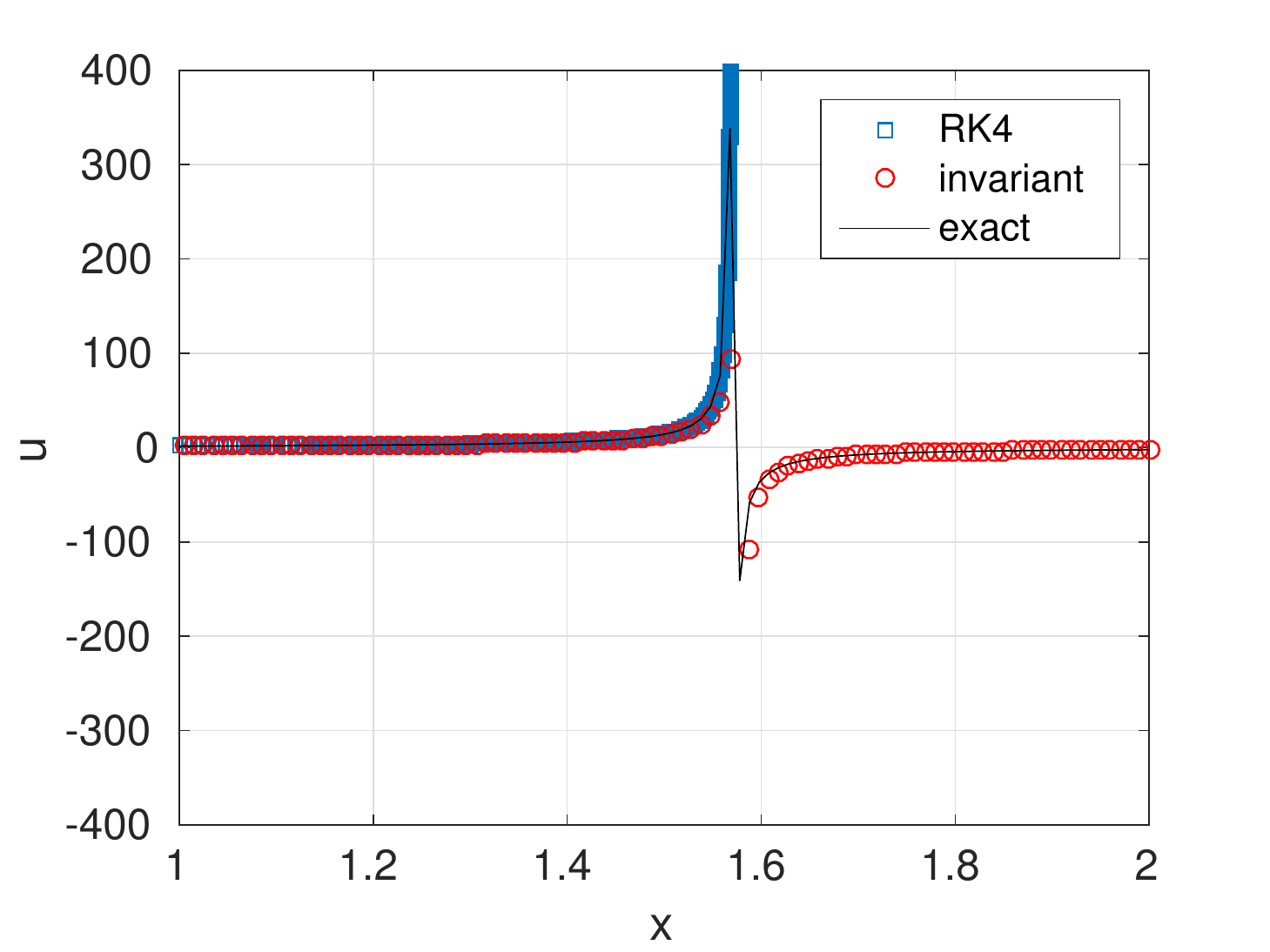}
\caption{\footnotesize{Numerical integration of the Schwarzian ODE~\eqref{Schwarzian equation} with $F(x)=2$. \textbf{Blue:} Non-invariant fourth order adaptive RK method. \textbf{Red:} Invariant first order method. \textbf{Black:} Exact solution.}}
\label{fig:NumericalResultsSchwarzianODE}
\end{figure}

We begin with the Schwarzian ODE~\eqref{Schwarzian equation} with $F(x)=2$. In other words, we consider the differential equation
\begin{equation}\label{particular schwarzian equation}
\frac{u_x\, u_{xxx} - (3/2)u_{xx}^2}{u_x^2} = 2.
\end{equation}
By the Schwarz' Theorem, \cite{M-2010}, the general solution of~\eqref{particular schwarzian equation} is
\[
u(x) = \frac{a \sin x + b \cos x}{c\sin x + d \cos x}\qquad \text{with}\qquad ad-bc \neq 0.
\]
Choosing $a=d=1$ and $b=c=0$, we obtain the particular solution $u(x) = \tan x$.  We now aim to obtain this particular solution numerically using the invariant scheme \eqref{Schwarzian invariant scheme} and a standard non-invariant scheme, and compare the results.  For the standard method, we choose the explicit fourth order adaptive Runge--Kutta solver \texttt{ode45} as provided by {\sc Matlab}. On the surface, this appears to be an unfair comparison since the invariant scheme~\eqref{Schwarzian invariant scheme} is only first order accurate. However, preserving geometric properties can give a numerical scheme a distinct advantage, even if it is only of relatively low order. This is verified in Figure~\ref{fig:NumericalResultsSchwarzianODE}. 

The relative error tolerance controlling the step size in the (non-invariant) adaptive Runge--Kutta method was set to $10^{-12}$. Despite this extremely small tolerance, the numerical solution diverges at the point $x=\pi/2$ where the solution has a vertical asymptote. On the other hand, the invariant method, with a step size of $h = 0.01$, is able to integrate beyond this singularity and follows the exact solution $u(x)=\tan x$ very closely.  For the conceptually related case where $F(x)=\sin x$ in \eqref{Schwarzian equation}, see~\cite{BCW-2006}.

\subsection{Korteweg--de Vries equation}\label{KdV section}

As reviewed in the previous sections, the earliest examples of invariant numerical schemes for evolution equations almost exclusively rested on the discretization of their associated Lagrangian form. However, the use of fully Lagrangian techniques for discretizing differential equations is not common due to their tendency to cluster grid points in certain areas of the computational domain and to poorly resolve the remaining parts of the domain.  Even more problematic, Lagrangian numerical methods regularly lead to mesh tangling, especially in the case of several space dimensions. 

\begin{figure}[!ht]
\centering
\captionsetup{width=0.85\textwidth}
\begin{subfigure}[b]{0.40\textwidth}
  \centering
  \includegraphics[width=\linewidth]{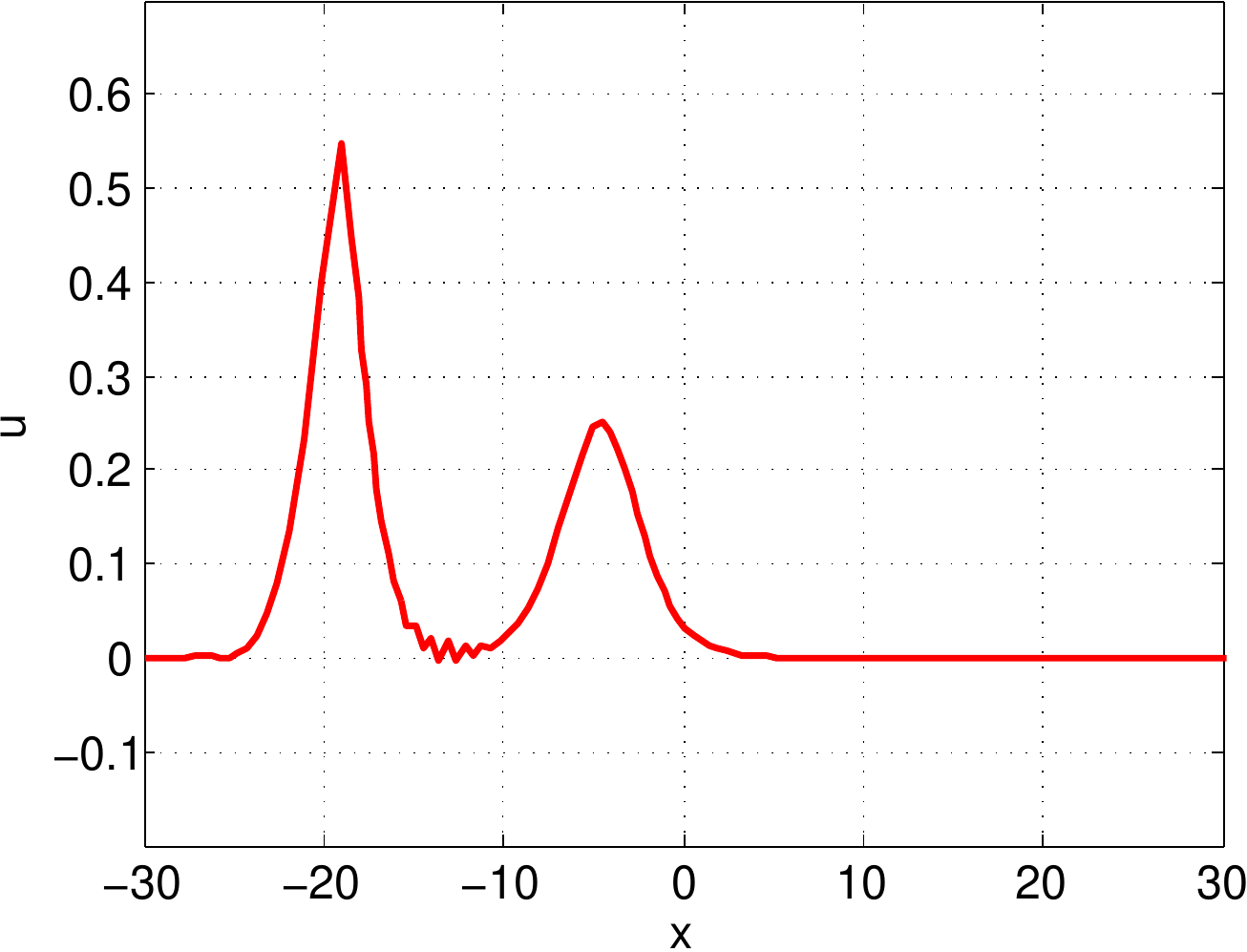}
\end{subfigure}\qquad
\begin{subfigure}[b]{0.40\textwidth}
  \centering
  \includegraphics[width=\linewidth]{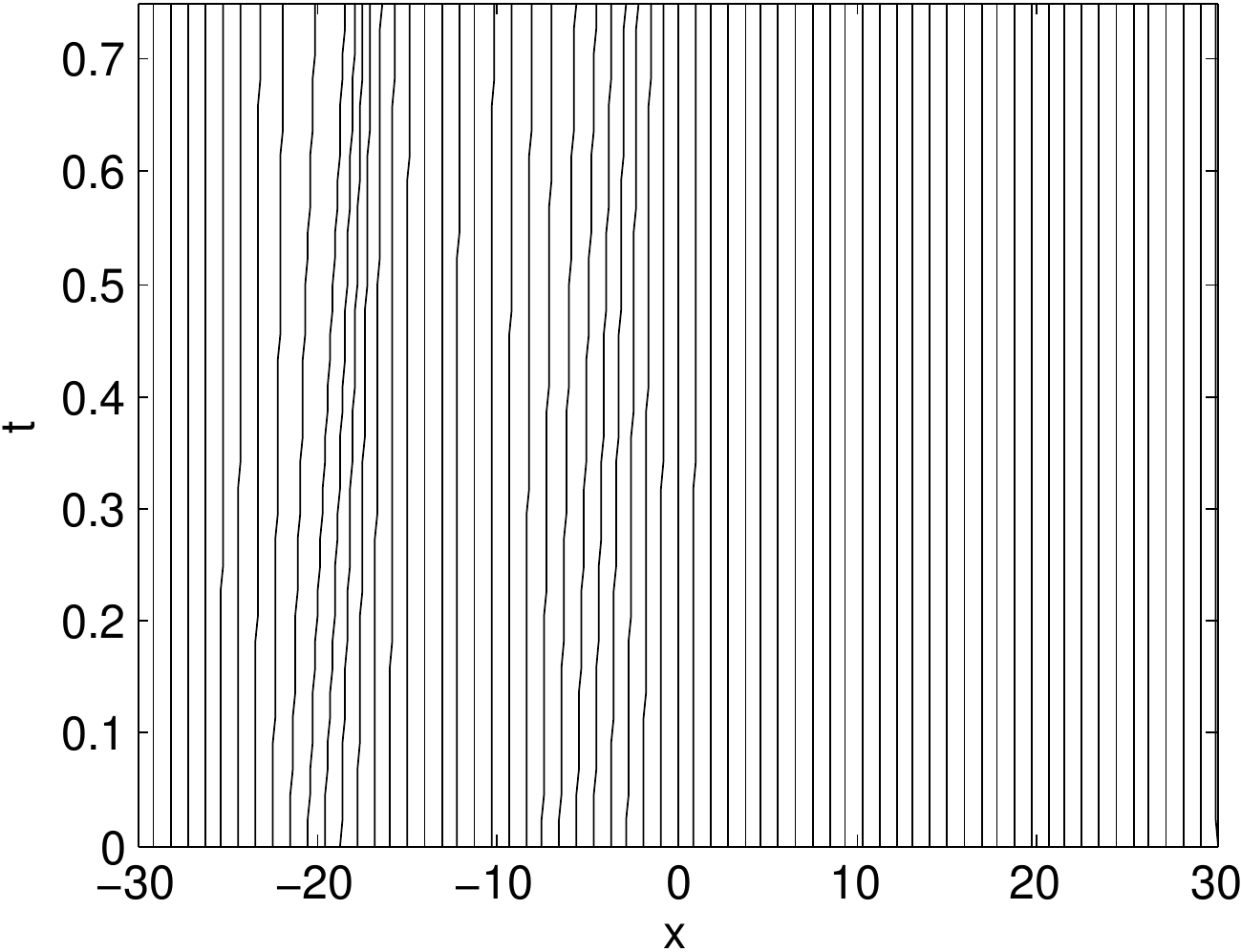}
\end{subfigure}
\caption{\footnotesize{\textbf{Left:} Numerical solution for the KdV equation using the invariant Lagrangian scheme~\eqref{coordinate 10 point invariant scheme} and~\eqref{KdV x mesh equation}. \textbf{Right:} Associated evolution of the mesh points.}}
\label{fig:NumericalResultsKdVLagrangian}
\end{figure}

For the KdV equation this basic problem is readily demonstrated using the invariant Lagrangian scheme given by \eqref{coordinate 10 point invariant scheme} and~\eqref{KdV x mesh equation}.  To do so, we numerically implement this scheme using as initial condition a double soliton solution of the form
\begin{equation}\label{double soliton}
 u(t,x)=\frac12c_1\, \textup{sech}^2\left(\frac{\sqrt{c_1}}{2}(x+a_1-c_2t)\right)+\frac12c_2\, \textup{sech}^2\left(\frac{\sqrt{c_2}}{2}(x+a_2-c_2t)\right),
\end{equation}
where $c_1, c_2\in\mathbb{R}$ are the phase velocities of the individual solitons and $a_1,a_2\in\mathbb{R}$ are the initial displacements. In our numerical simulations, we set $a_1=20$, $a_2=5$, $c_1=1$, $c_2=0.5$, and restricted the computational domain to the interval $[-30,30]$ discretized with a total of $N=128$ spatial grid points. The time step~$k$ was chosen to be proportional to $h^3$, $k\propto h^3$, and the final integration time for the Lagrangian experiment was $t=0.75$. The result of the computation is presented in Figure~\ref{fig:NumericalResultsKdVLagrangian}. 

It is visible from the evolution of the mesh points in Figure~\ref{fig:NumericalResultsKdVLagrangian} (right) that mesh tangling (here: crossing of mesh lines) is about to occur. This leads to numerical instability that is visible as high wavenumber oscillations between the two solitons, which sets in almost immediately after the start of the numerical integration. In other words, the invariant Lagrangian scheme is unsuitable for practical applications in virtually all relevant numerical situations, since no inherent control over the evolution of the mesh points is build into the scheme.

It is striking to observe that, despite the well-known shortcomings of Lagrangian numerical schemes, the latter have played a prominent role in the field of symmetry-preserving discretization.  This can be explained by the fact that for most papers on the subject, invariant numerical schemes were constructed more as an example highlighting the possibility of deriving symmetry-preserving schemes of differential equations rather than as a tool for practical numerical experiments. Indeed, it is fair to say that even now, the numerical analysis of invariant discretization schemes is still lacking rigor.

In order to make the invariant numerical scheme~\eqref{coordinate 10 point invariant scheme} 
practical, a different invariant grid equation has to be derived. Possible strategies include the use of \textit{invariant evolution--projection schemes} and \textit{invariant adaptive numerical schemes}. 

The invariant evolution--projection scheme conceptually builds upon the invariant Lagrangian scheme. The main idea of this approach is to use the invariant numerical scheme and the invariant mesh equations only over a single time step, and use an interpolation scheme to project the solution of the differential equation defined at the new spatial grid points back to the initial (typically uniform) spatial grid. The entire procedure is invariant if the interpolation scheme used is invariant~\cite{BN-2013}. The main appeal of this method is that it enables the use of invariant numerical schemes on rectangular meshes.

It is readily verified that classical interpolation methods such as linear, quadratic, cubic or spline interpolations are all invariant under the maximal symmetry group of the KdV equation. The main reason is that all these schemes are polynomials in terms like $(x_{i+1}^{n+1}-x_i^{n+1})$ or $(\hat x^{n+1}-x_i^{n+1})/(x_{i+1}^{n+1}-x_i^{n+1})$, where $\hat x^{n+1}$ is the interpolation point, which are invariant under spatial and temporal shifts and Galilean boosts.  The invariance under scaling transformations follows from the consistency of the interpolation scheme.  For example, consider the linear interpolation given by
\[
 u(\hat x^{n+1}_i)=u_i^{n+1}+\frac{u^{n+1}_{i+1}-u^{n+1}_i}{x^{n+1}_{i+1}-x^{n+1}_i}(\hat x^{n+1}_i-x_i^{n+1}).
\]
Then, under the action of the KdV symmetry group given by~\eqref{eq:SymmetryGroupKdVEquation}, the linear interpolation formula remains invariant. In other words, it follows that
\[
 U(\widehat X^{n+1}_i)=U_i^{n+1}+\frac{U^{n+1}_{i+1}-U^{n+1}_i}{X^{n+1}_{i+1}-X^{n+1}_i}(\widehat X^{n+1}_i-X_i^{n+1}).
\]
For more details and examples, see~\cite{BN-2013,BN-2014}.

\begin{figure}[!ht]
\centering
\captionsetup{width=0.85\textwidth}
  \includegraphics[width=0.45\linewidth]{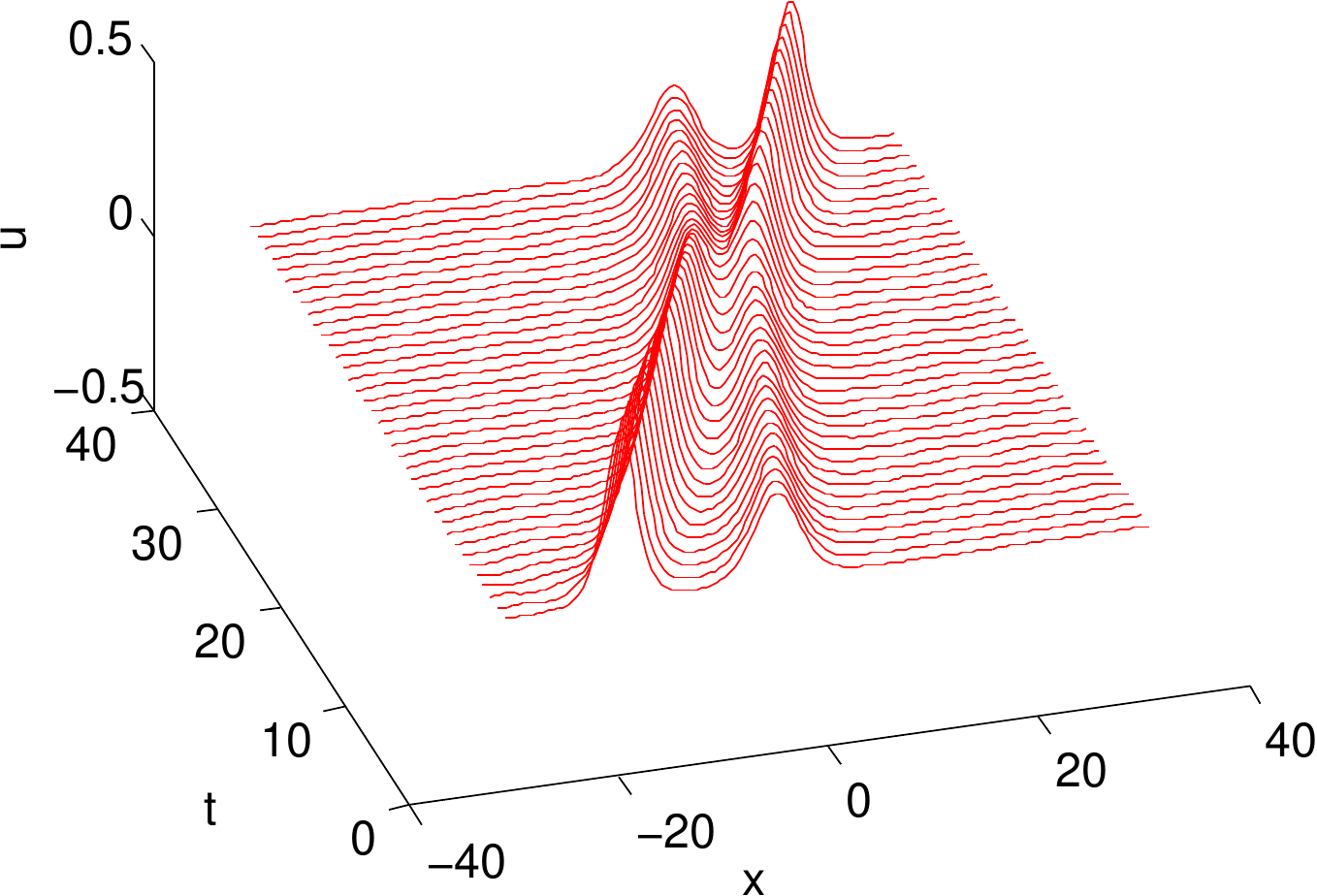}
\caption{\footnotesize{Numerical solution for the KdV equation using the invariant scheme~\eqref{coordinate 10 point invariant scheme}, \eqref{KdV x mesh equation}  augmented with a cubic spline interpolation after every step to project the solution back to the original uniform mesh.}}
\label{fig:NumericalResultsKdVEvolveProject}
\end{figure}

In Figure~\ref{fig:NumericalResultsKdVEvolveProject} we present the numerical results for the scheme~\eqref{coordinate 10 point invariant scheme}, \eqref{KdV x mesh equation} using the double soliton \eqref{double soliton} as initial condition on the interval $[-30, 30]$. As opposed to the previous simulation, we now introduce a cubic spline interpolation at each time integration to project the solution back to the original space grid.  The number of discrete spatial points is as before, that is $N=128$, and the final integration time is $t=40$. As it can be seen, the two solitons interact with each other and remain unchanged after their collision, which is properly captured by the invariant evolution--projection scheme. We point out though that the scheme is rather dissipative, with the amplitudes of the solitons slowly decreasing over time. While in the present example dissipation can be seen as a disadvantage, this dissipation can be essential in hyperbolic problems that involve shock solutions.  For these shock solutions, numerical simulations usually require schemes, such as upwind or Lax--Friedrich and Lax--Wendroff methods, which exhibit artificial dissipation.

A second possibility for completing the invariant numerical scheme~\eqref{coordinate 10 point invariant scheme} without using Lagrangian methods rests on moving mesh methods. Without going into great details, we present here an \emph{invariant $r$-adaptive scheme} for the KdV equation (for more information, see~\cite{BCW-2015}). In $r$-adaptive numerical schemes a fixed number of grid points is redistributed so that points automatically move to regions where higher resolution is required, for example near shocks.  Therefore, $r$-adaptive numerical methods are particularly important for hyperbolic problems.  We refer to \cite{HR-2011} for a comprehensive review of such methods. 

For one-dimensional problems, $r$-adaptive moving meshes on the interval $[a,b]$ are uniquely determined through the \textit{equidistribution principle}, which in differential form reads
\begin{equation}\label{eq:equidistributionPrinciple}
 (\delta(t,x)x_s)_s=0
\end{equation}
with boundary conditions $x(t,0)=a$ and $x(t,1)=b$. In \eqref{eq:equidistributionPrinciple}, the function $\delta$ is called the \textit{mesh density function} or \emph{monitor function}. Its role is to control the areas where grid points should concentrate or de-concentrate. It is typically linked to the solution of the physical differential equation. For example, the arc-length type mesh density function is
\begin{equation}\label{eq:arcLengthMonitorFunction}
 \delta=\sqrt{1+\alpha u_x^2},
\end{equation}
where $\alpha\in\mathbb{R}$ is a constant adaptation parameter. 

To complete the invariant scheme for the KdV equation, we discretize~\eqref{eq:equidistributionPrinciple} and~\eqref{eq:arcLengthMonitorFunction} using the difference invariants given in~\eqref{eq:DifferenceInvariantsKdV} or using the invariantization map induced by the discrete moving frame \eqref{KdV moving frame}. In particular, it turns out that the straightforward discretization
\begin{align}\label{eq:invariantAdaptiveSchemeKdV}
\begin{split}
& \frac{\delta^n_{i+1}+\delta^n_i}{2}(x_{i+1}^{n+1}-x_i^{n+1})- \frac{\delta^n_{i}+\delta^n_{i-1}}{2}(x_{i}^{n+1}-x_{i-1}^{n+1})=0,\\
& \delta^n_i=\sqrt{1+\alpha \left(k\frac{u^n_{i+1}-u^n_i}{x^n_{i+1}-x^n_i}\right)^2},
\end{split}
\end{align}
is invariant under the maximal Lie symmetry group of the KdV equation. 

\begin{figure}[!ht]
\centering
\captionsetup{width=0.85\textwidth}
\begin{subfigure}[b]{0.40\textwidth}
  \centering
  \includegraphics[width=\linewidth]{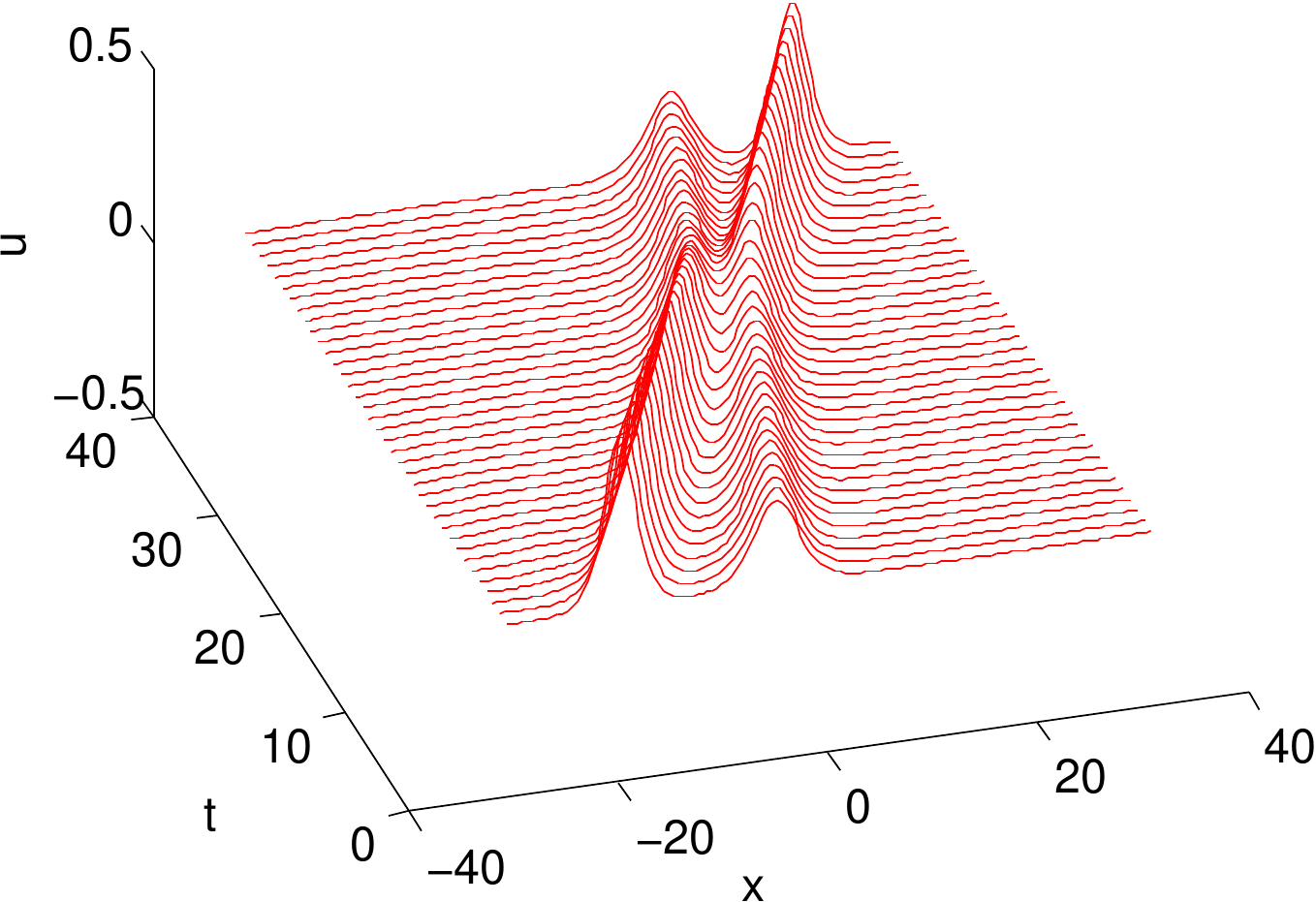}
\end{subfigure}\qquad
\begin{subfigure}[b]{0.40\textwidth}
  \centering
  \includegraphics[width=\linewidth]{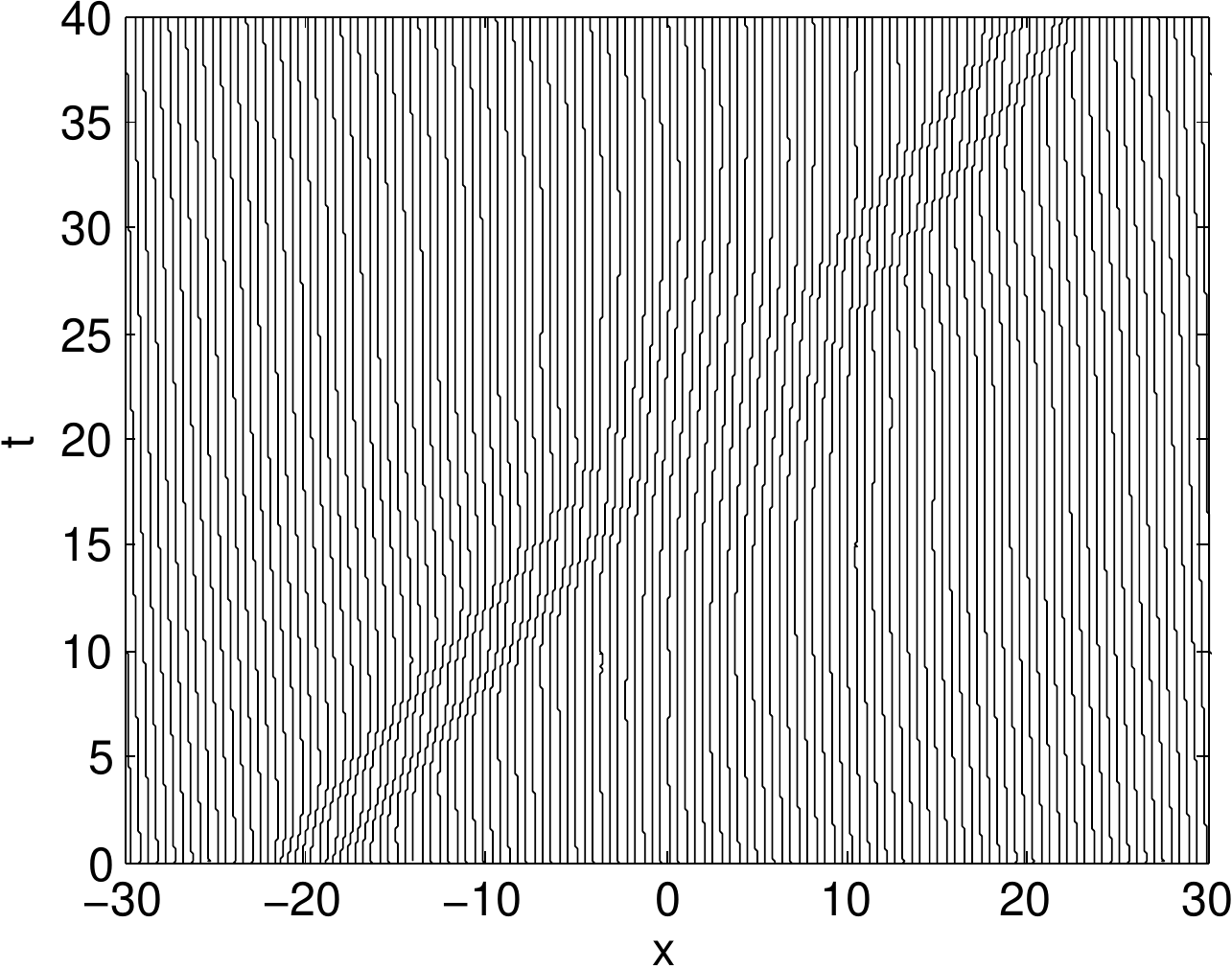}
\end{subfigure}
\caption{\footnotesize{\textbf{Left:} Numerical solution for the KdV equation using the invariant adaptive scheme~\eqref{coordinate 10 point invariant scheme}, \eqref{eq:invariantAdaptiveSchemeKdV}. \textbf{Right:} Associated evolution of the mesh points.}}
\label{fig:NumericalResultsKdVAdaptive}
\end{figure}

In Figure~\ref{fig:NumericalResultsKdVAdaptive} we present the numerical solution for the KdV equation using the invariant adaptive scheme~\eqref{coordinate 10 point invariant scheme} with~\eqref{eq:invariantAdaptiveSchemeKdV} and the same double soliton initial condition \eqref{double soliton} as in the previous simulation.  The final integration time was  again chosen to be $t=40$ and the adaptation parameter was set to $\alpha=10$.

It is readily seen from Figure~\ref{fig:NumericalResultsKdVAdaptive} that the invariant adaptive scheme \eqref{coordinate 10 point invariant scheme}, \eqref{eq:invariantAdaptiveSchemeKdV} again does not suffer from the shortcomings observed for the Lagrangian scheme \eqref{coordinate 10 point invariant scheme}, \eqref{KdV x mesh equation}. In particular, no mesh tangling occurs. The associated adaptive mesh suitably tracks the position of the solitons and remains almost uniform away from the two waves, although the adaptation is relatively weak since the solution does not exhibiting overly steep gradients. An advantage of the invariant $r$-adaptive scheme over the invariant evolution--projection scheme is that the amplitudes of the solitons are not damped during the adaptation strategy.

In general, using invariant adaptive schemes has the merit of combining a geometric numerical method with a well-proven numerical strategy for dynamically redistributing the points in a mesh. In particular, this technique works for all evolution equations that are invariant under the Galilean group, which are virtually all equations of classical hydrodynamics, including the shallow-water equations, the Euler equations and the Navier--Stokes equations. Since shock waves are physically important solutions in these models, invariant adaptive schemes are of high practical relevance in this field.

\subsection{Burgers' equation}

In this section we construct a new numerical scheme for Burgers' equation \eqref{burgers equation} invariant under the four-parameter symmetry group 
\begin{equation}\label{burgers symmetry subgroup}
X= e^{\epsilon_4} x + \epsilon_3 t + \epsilon_1,\qquad T = e^{2\epsilon_4}t+\epsilon_2,\qquad U=e^{-\epsilon_4}u+\epsilon_3, \qquad 
\epsilon_1, \epsilon_2, \epsilon_3, \epsilon_4 \in \mathbb{R}
\end{equation}
generated by the vector fields~$\vv_1$, $\vv_2$, $\vv_3$, $\vv_4$ given in \eqref{burgers symmetry generators}. We exclude the one-parameter group of transformations generated by~$\vv_5$ since in numerical simulations the evolution of the time variable $t$ should always be strictly increasing, and allowing the inversion transformations generated by $\vv_5$ would enable one to reverse the time direction, which is not desirable from a numerical standpoint.

\begin{remark}
As an exercise, the reader is invited to adapt the constructions below by including the inversion transformations generated by $\vv_5$.  This has never been attempted and could potentially lead to interesting new results!
\end{remark}

Due to the similarities between the symmetry subgroup action \eqref{burgers symmetry subgroup} and the KdV symmetry group \eqref{KdV symmetry generators}, the underlying symmetry-preserving schemes for Burgers' equation are conceptually similar to the invariant schemes constructed before for the KdV equation.   Though, one important differences between the two equations is that solutions to Burgers' equation can develop very steep gradients (although remaining smooth for all times provided that $\nu\ne0$). This is particularly the case if $\nu$ approaches zero. Hence, grid adaptation is of practical relevance for this equation.

In~\cite{K-2008}, an invariantization for the Crank--Nicolson scheme for Burgers' equation was proposed. However, we note that the Crank--Nicolson scheme is implicit and thus in the case where~$\nu$ is small it might not be the most efficient way of solving Burgers' equation since an explicit scheme should then suffice.  In the following, we propose a new scheme which draws some ideas from \textit{high-resolution finite volume methods}, \cite{L-2002}.   It is well-known that high order schemes, such as the Lax--Wendroff method, lead to oscillations in the numerical solution near shocks, whereas low order schemes, such as the upwind method, develop no such oscillations but exhibit an excessive amount of numerical viscosity. The idea in the high-resolution method is thus to use a high order method away from the shock and a low resolution method near the shock. The transition between the two regions is accomplished through the use of \textit{flux/slope limiters}.

To formulate an invariant finite volume type method for Burgers' equation, we rewrite~\eqref{burgers equation} in the form  
\begin{equation}\label{eq:BurgersFluxForm}
  u_t+\widetilde{f}_x=0,\qquad  \widetilde{f}=\frac12u^2-\nu u_x.
\end{equation}
We now discretize~\eqref{eq:BurgersFluxForm} on a moving mesh, which, as we have previously seen, is enough to guarantee invariance under Galilean transformations. As in Examples \ref{change of variable example} and \ref{KdV in computational variables}, we introduce the computational variables $(\tau,s)$ and let $t=t(\tau)=k\tau+t^0$ and $x=x(\tau,s)$.  Then, a suitable conservative form of Burgers' equation in the computational variables $(\tau, s)$ is given by
\begin{equation}\label{eq:BurgersFluxFormComputationalVariables}
(x_su)_\tau+k\left(\frac12u^2-\nu\frac{u_s}{x_s}-\frac{ux_\tau}{k}\right)_s=(x_su)_\tau+k f_s=0,
\end{equation}
see also~\cite{HR-2011}. We then discretize the flux~$f$ in two different ways, once using the second order centered difference method (high resolution) and once using the first order upwind method (low resolution). In doing so, we observe, as in~\cite{BCW-2015}, that the invariance under Galilean transformations requires us to discretize~\eqref{eq:BurgersFluxFormComputationalVariables} in such a way that all spatial derivatives are evaluated using the same finite difference discretizations. With that said, the high order discretization of~\eqref{eq:BurgersFluxFormComputationalVariables} is
$
 \Delta_\tau(x_su)^{\rm high}+k\, \Delta_sf^{\rm high}=0,
$
with
$$
\begin{aligned}
 &\Delta_\tau (x_su)^{\rm high}=(h_i^{n+1}+h_{i-1}^{n+1})u_i^{n+1}-(h_i^n+h_{i-1}^n)u_i^n,\\
 &\Delta_sf^{\rm high}=\frac12\left[(u_{i+1}^n)^2-(u_{i-1}^n)^2\right]-\nu(Du_i^n-Du_{i-1}^n)-\left(\frac{\sigma_{i+1}^n}{k}u_{i+1}^n-\frac{\sigma_{i-1}^n}{k}u_{i-1}^n\right).
\end{aligned}
$$
On the other hand, the low order discretization of~\eqref{eq:BurgersFluxFormComputationalVariables} is
$
 \Delta_\tau (x_su)^{\rm low}+k\,\Delta_sf^{\rm low}=0,
$
where
\[
\Delta_\tau(x_su)^{\rm low}=\begin{cases}
h_{i-1}^{n+1}u_i^{n+1}-h_{i-1}^nu_i^n, & u_i^n\ge0,\\ h_i^{n+1}u_i^{n+1}-h_i^nu_i^n, & u_i^n<0,
\end{cases}
\]
and
\[  
\Delta_sf^{\rm low}=\begin{cases}
\Delta_sf^{\rm low}_\ge, & u_i^n\ge0,\\ \Delta_sf^{\rm low}_< ,& u_i^n<0,
\end{cases}
\]
with
\begin{align*}
 & \Delta_sf^{\rm low}_\ge=\cfrac12\left[(u_{i}^n)^2-(u_{i-1}^n)^2\right]-\nu(Du_{i-1}^n-Du_{i-2}^n)-\bigg(\cfrac{\sigma_i^n}{k}u_i^n-\cfrac{\sigma_{i-1}^n}{k}u_{i-1}^n\bigg),\\
 & \Delta_sf^{\rm low}_<=\cfrac12\left[(u_{i+1}^n)^2-(u_i^n)^2\right]-\nu(Du_{i+1}^n-Du_i^n)-\bigg(\cfrac{\sigma_{i+1}^n}{k}u_{i+1}^n-\cfrac{\sigma_i^n}{k}u_i^n\bigg).
\end{align*}
The invariant high-resolution method is obtained by dynamically selecting the regions of the domain where the high order and low order methods are used. For this purpose, we introduce the ratio
\[
 \theta^n_i=\frac{\Delta u^n_{I-1}}{\Delta u^n_{i-1}},\qquad \text{where}\qquad \Delta u^n_{i-1} = u^n_i - u^n_{i-1},
\]
and $I=i-1$ if $u_i^n\ge0$ and $I=i+1$ if $u_i^n<0$. Geometrically, the quantity $\theta^n_i$ measures the smoothness of the solution over the interval $[x_{i-1},x_i]$. This ratio is, by its definition, invariant under the symmetry subgroup \eqref{burgers symmetry subgroup}, and therefore so is any function of~$\theta^n_i$.

We proceed to discretize~\eqref{eq:BurgersFluxFormComputationalVariables} by considering
%
\begin{equation}\label{eq:InvariantFiniteVolumeBurgersEquation}
 \Delta_\tau (x_su)+k \Delta_sf=0,
\end{equation} 
with
\begin{align*}
 &\Delta_\tau(x_su) =\Delta_\tau(x_su)^{\rm low}-\Phi(\theta^n_{i-1})\left[\Delta_\tau(x_su)^{\rm low}-\Delta_\tau(x_su)^{\rm high}\right],\\
 &\Delta_sf =\Delta_s f^{\rm low}-\Phi(\theta^n_{i-1})\left[\Delta_sf^{\rm low}-\Delta_sf^{\rm high}\right],
\end{align*}
and where, for the \textit{flux limiter function} $\Phi(\theta^n_i)$, we choose the so-called minmod-limiter, $\Phi(\theta^n_i)=\max\{0,\min(1,\theta^n_i)\}$. For further discussions on flux limiters, see~\cite{L-2002}. To complete the invariant finite volume type scheme for Burgers' equation, we use the same grid adaptation strategy as for the KdV equation to obtain the spatial step size $\sigma^n_i=x^{n+1}_i - x^n_i$ as time evolves.

\begin{figure}[!ht]
\centering
\captionsetup{width=0.85\textwidth}
\begin{subfigure}[b]{0.40\textwidth}
  \centering
  \includegraphics[width=\linewidth]{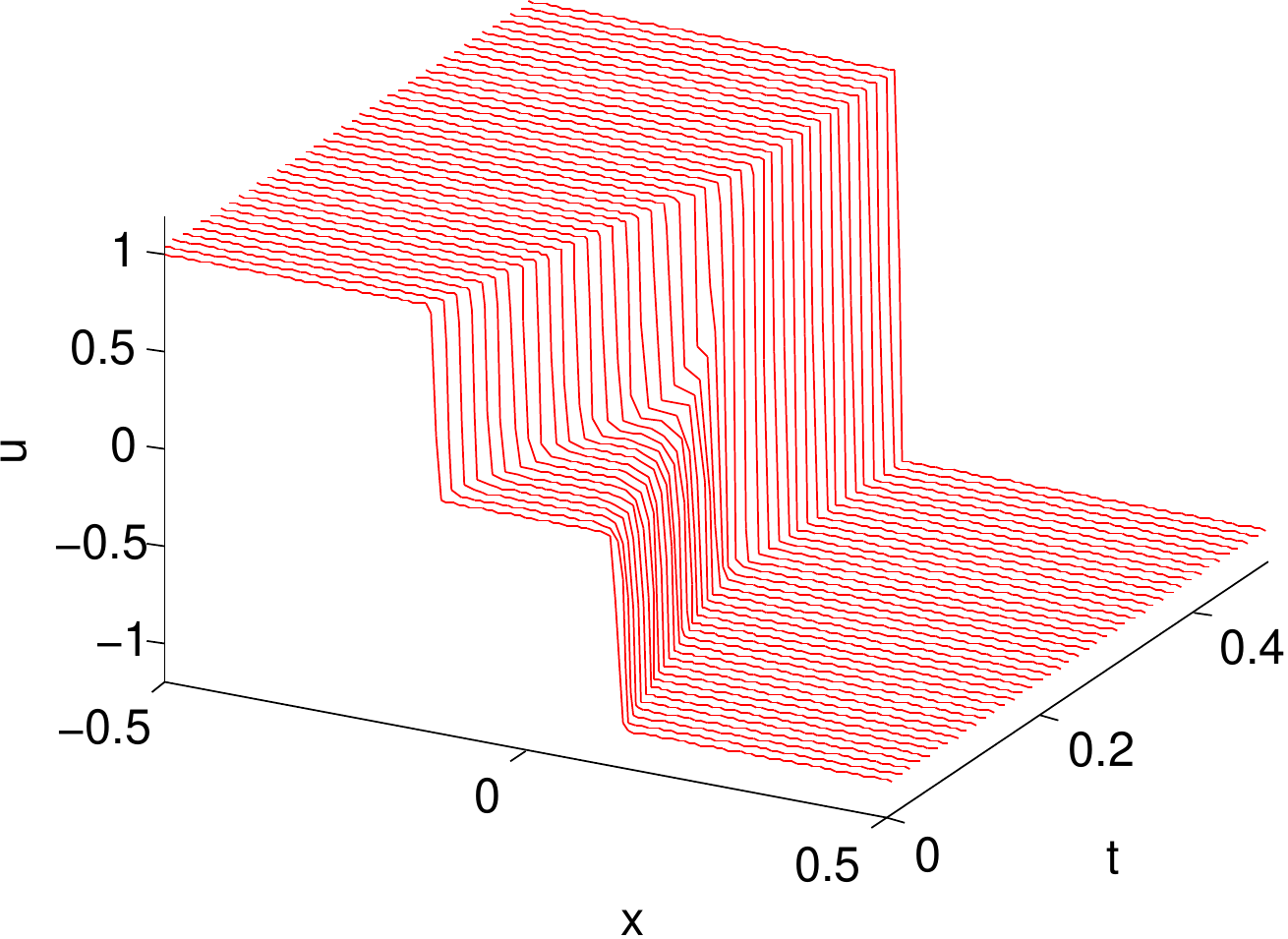}
\end{subfigure}\qquad
\begin{subfigure}[b]{0.40\textwidth}
  \centering
  \includegraphics[width=\linewidth]{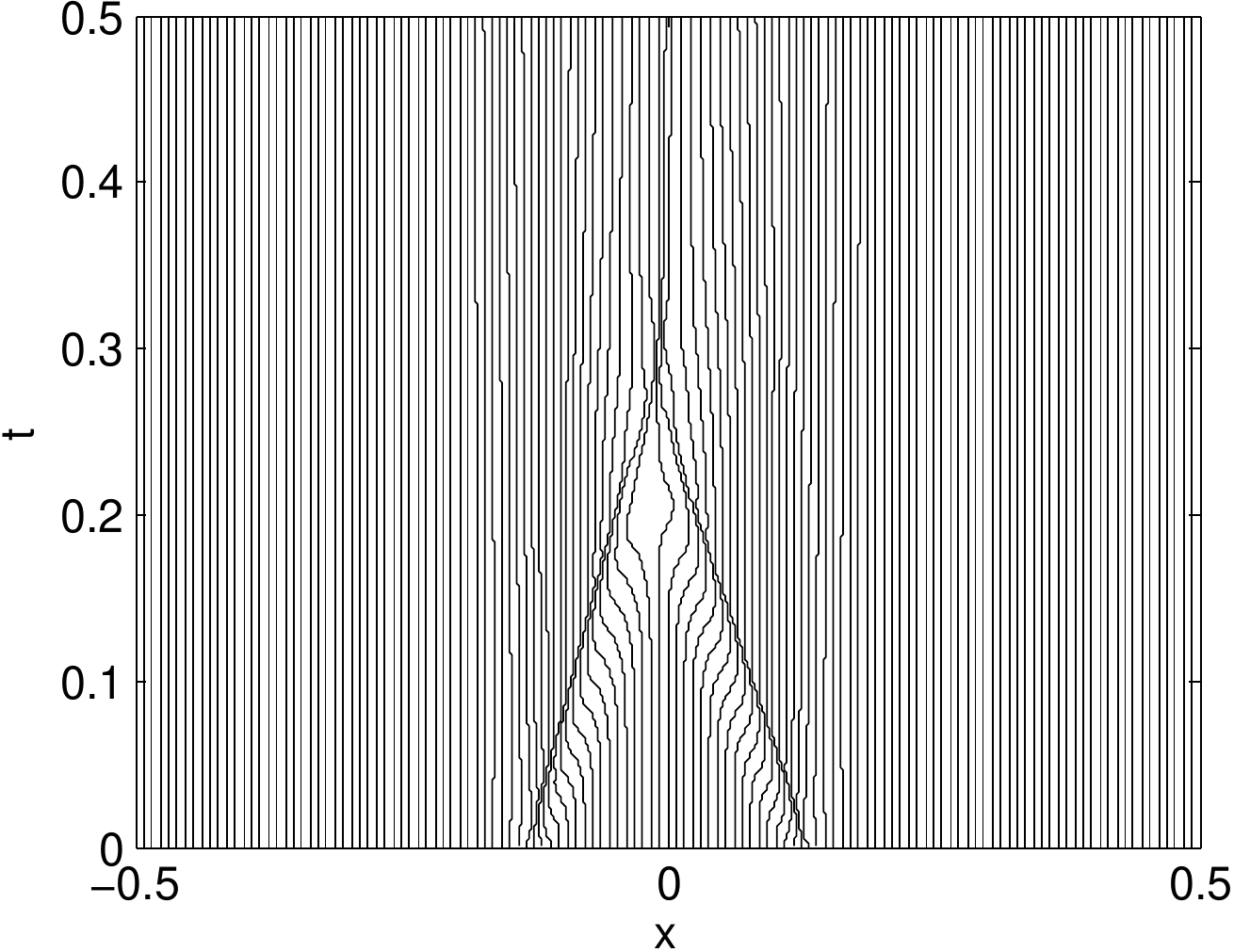}
\end{subfigure}
\caption{\footnotesize{\textbf{Left:} Numerical solution of Burgers' equation using the invariant adaptive scheme~~\eqref{eq:invariantAdaptiveSchemeKdV}, \eqref{eq:InvariantFiniteVolumeBurgersEquation}. \textbf{Right:} Corresponding evolution of the mesh points.}}
\label{fig:NumericalResultsBurgersAdaptive}
\end{figure}

As a numerical example, we carry out an experiment similar to the one given in~\cite{K-2008} for the exact solution
\[
 u(x)=-\frac{\sinh\left(\frac{x}{2\nu}\right)}{\cosh\left(\frac{x}{2\nu}\right)+\exp\left(-\frac{(c+t)}{4\nu}\right)},
\]
where $c\in\mathbb{R}$. We discretize the spatial domain $[-0.5,0.5]$ with $N=128$ grid points using Dirichlet boundary conditions, and choose the time step~$k$ to be proportional to $h^2$, $k\propto h^2$. The final integration time is $t=0.5$ and for numerical purposes $c=0.25$ and the viscosity was set to $\nu=0.001$. The adaptation parameter $\alpha$ in the arc-length type mesh density function in~\eqref{eq:invariantAdaptiveSchemeKdV} was set to $\alpha=0.5$. The respective numerical results are depicted in Figure~\ref{fig:NumericalResultsBurgersAdaptive}.

Unlike in the numerical simulations for the KdV equation presented in Section \ref{KdV section}, Figure~\ref{fig:NumericalResultsBurgersAdaptive} clearly demonstrates the need for an adaptive moving mesh. While this is implied from the structure of the numerical solution, it is remarkable that the requirement for a moving mesh is already encoded in the structure of the symmetry group of Burgers' equation. Hence, numerically preserving symmetries can be seen as a geometrical justification for using $r$-adaptive numerical methods. Moreover, due to the use of a high-resolution finite volume type scheme, no unphysical oscillations around the shock is observed.

\section{Conclusion}\label{conclusion section}

To recapitulate, let us summarize the algorithm for constructing symmetry-preserving finite difference schemes.  Given a differential equation $\Delta(x,u\n)=0$:
\begin{enumerate}
\item Use the infinitesimal invariance criterion \eqref{Delta infinitesimal symmetry criterion} to determine a basis of infinitesimal symmetry generators. \label{symmetry step}
\item Choose a lattice on which the differential equation is to be discretized.
\item When possible, in particular when discretizing a partial differential equation, impose obvious invariant constraints on the mesh.  This step is not necessary but if implemented it can, in general, simplify the implementation of the remaining steps.
\item Use either Lie's infinitesimal approach or the moving frame method to compute a complete set of difference invariants, and, if necessary, to find weakly invariant difference equations.  When using the moving frame method, one has to exponentiate the infinitesimal generators found in Step \ref{symmetry step} to obtain the connected component of the group of local symmetry transformations.
\item Combine the difference invariants and the weakly invariant equations in such a way to obtain an approximation of the differential equation $\Delta(x,u\n)=0$ and (possibly) constraints on the mesh.  If using the moving frame method, invariantize a finite difference approximation of the differential equation compatible with the mesh to obtain an invariant approximation of the differential equation.
\end{enumerate}

The basic algorithm for constructing symmetry-preserving numerical schemes is now fairly well-understood.  Below are some open problems and comments for the interested reader.

\begin{itemize}
\item Many important differential equations in mathematical physics admit an infinite-dimensional symmetry group.  Such equations include the Davey--Stewartson equations, \cite{CW-1988}, Liouville's equation, the Kadomtsev--Petviashvili equation, \cite{DKLW-1986}, the Infeld--Rowland equation, \cite{FW-1993}, the Euler equations, \cite{O-1993}, and many other equations from fluid dynamics, \cite{II-2011}. Implementing the above algorithm for infinite-dimensional symmetry groups remains a challenge. One particularity of these groups is that as new points are added to the stencil, new group parameters appear, which does not occur in the finite-dimensional case.  To avoid this difficulty, one possibility is to consider finite-dimensional subgroups of the infinite-dimensional symmetry group and implement the algorithm above, \cite{LMW-2015,LMW-2015-2}.  Another possibility, which preserves the infinite-dimensional nature of the group action, is to discretize the Lie pseudo-group action, \cite{RV-2015}. 

\item In the last 25 years, a great deal of efforts has been devoted to constructing symmetry-preserving finite difference numerical schemes.  With the emergence of finite element methods, \cite{R-2006}, and discrete exterior calculus, \cite{AFW-2010}, it would be interesting to extend the above symmetry-preserving algorithm to these settings as well.  Further extensions to finite volume and spectral methods should also be considered.

\item As with any geometric integrator, one of the motivations for developing symmetry-preserving schemes is to obtain better long term numerical results.   As we saw in Section \ref{numerical simulations section}, and as observed in the literature, \cite{BCW-2006,BRW-2008,CRW-2016,KO-2004}, symmetry-preserving schemes for ordinary differential equations perform extremely well, particularly near singularities.  For first order ordinary differential equation, it is even possible to construct symmetry-preserving schemes that will approximate exactly the solution of the original equation, \cite{RW-2004}.  On the other hand, the numerical improvements for partial differential equations are not as clear, \cite{B-2013,BCW-2015,CH-2010,K-2008,LMW-2015,RV-2013,RV-2015}.  In many cases, they tend to be comparable to standard schemes. Now that the theoretical foundations are on firm grounds, one of the main challenges in the field of symmetry-preserving schemes is to investigate the numerical properties of invariant schemes and understand why and when these schemes give better numerical results.

\item Most partial differential equations invariantly discretized to date have been evolutionary equations (such as the KdV and Burgers' equations).  Much more work, especially from the numerical side, has to be devoted to the invariant discretization of other types of partial differential equations, such as the wave equation, Laplace's equation, and the Sine--Gordon equation.  In particular, constructing symmetry-preserving schemes compatible with given boundary conditions is an important avenue of research.

\item Symmetries are usually not the only geometric properties that a differential equation admits. Other, equally important properties such as a Hamiltonian structure or conservation laws might be present as well. Developing geometric integrators that will preserve more than just one geometric property at the time is an important research direction to pursue.


\end{itemize}

\section*{Acknowledgement}
The research of the first author is supported in part by a Tier 2 NSERC Canada Research Chair grant.  The authors would like to thank the organizers of the ASIDE summer school for inviting them to give a series of lectures on continuous symmetries of discrete equations. We also thank Peter J.\ Olver for his comments on our lecture notes.


\begin{appendices}
\section{Answers to selected exercises}
\begin{description}
\item[Exercise \ref{invariants Lie's approach exercise}, part (\ref{invariants Lie's approach exercise part b}):]  A complete set of difference invariants is given by the 9 invariants
\begin{align*}
&I_1 = \frac{h^n_i}{h^n_{i-1}},\qquad I_2 = \frac{h^{n+1}_i}{h^{n+1}_{i-1}}, &\qquad &I_3=\frac{h^n_i h^{n+1}_i}{k^n},\\
&I_4 = h^n_ih^n_{i-1}(Du^n_i - Du^n_{i-1}),
& &I_5 = h^{n+1}_ih^{n+1}_{i-1}(Du^{n+1}_i - Du^{n+1}_{i-1}),\\
&I_6 = h^n_i \bigg(\frac{\sigma^n_i}{k^n} - u^n_i\bigg), 
& &I_7 = h^{n+1}_i \bigg(\frac{\sigma^n_i}{k^n} - u^{n+1}_i\bigg),\\
&I_8 = (h^n_i)^2\bigg(Du^n_i + \frac{1}{k^n}\bigg),
& &I_9 = (h^{n+1}_i)^2 \bigg(Du^{n+1}_i - \frac{1}{k^n}\bigg),
\end{align*}
where
\[
h^n_i = x^n_{i+1}-x^n_i,\qquad k^n = t^{n+1}-t^n, \qquad
\sigma^n_i = x^{n+1}_i - x^n_i,\qquad Du^n_i = \frac{u^n_{i+1}-u^n_i}{h^n_i}.
\]

\item[Exercise \ref{differential moving frame exercise} part (\ref{differential moving frame exercise part a}):]  The one-parameter group actions are
\begin{align*}
\exp[\epsilon_1\vv_1]\cdot(t,x,u) &= (x+\epsilon_1,t,u),\\
\exp[\epsilon_2 \vv_2]\cdot(t,x,u) &= (x,t+\epsilon_2,u),\\
\exp[\epsilon_3 \vv_3]\cdot(t,x,u) &= (x+\epsilon_3 t, t, u+\epsilon_3),\\
\exp[\epsilon_4 \vv_4]\cdot(t,x,u) &= (e^{\epsilon_4}x,e^{2\epsilon_4}t,e^{-\epsilon_4}u),\\
\exp[\epsilon_5 \vv_5]\cdot(t,x,u) &= \bigg(\frac{x}{1-\epsilon_5 t},\frac{t}{1-\epsilon_5 t}, (1-\epsilon_5 t) u +\epsilon_5 x\bigg).
\end{align*}
\item[Exercise \ref{differential moving frame exercise} part (\ref{differential moving frame exercise part b}):]  Working on the open dense set $\mathcal{V}^{(1)} = \{(t,x,u,u_t,u_x) \in \J^{(1)}\,|\, uu_x + u_t \neq 0\}$, the right moving frame corresponding to the cross-section 
\[
\mK = \{ t=0,\; x=0,\; u=0,\; u_x = 0,\; u_t = 1\}
\]
is
\[
\epsilon_1 = -x,\qquad \epsilon_2 = -t,\qquad \epsilon_3 = -u,\qquad e^{\epsilon_4} = (uu_x+u_t)^{1/3},\qquad \epsilon_5 = -u_x.
\]

\item[Exercise \ref{differential moving frame exercise} part (\ref{differential moving frame exercise part c}):]  The invariantization of $u_{xx}$ yields the differential invariant
\[
\iota(u_{xx}) = \frac{u_{xx}}{uu_x+u_t}.
\]

\item[Exercise \ref{order 1 ode exercise} part (\ref{order 1 ode exercise part b}):]  A weakly invariant equation is given by 
\[
u_{i+1} e^{-A(x_{i+1})} - u_i e^{-A(x_i)} - B(x_{i+1}) + B(x_i) = 0.
\]

\item[Exercise \ref{Lie symmetry-preserving scheme exercise}:]  Along with the equations \eqref{t constraints exercise}, we can add the mesh equation $I_6 = 0$ (refer to the solution of Exercise \ref{invariants Lie's approach exercise}, part (\ref{invariants Lie's approach exercise part b})).  On this mesh, the differential equation can be approximated by 
\[
-I_7 I_3 = \frac{2\nu I_4 I_1I_2}{1+I_1}.
\]
Explicitly,
\[
\frac{h^{n+1}_i h^{n+1}_{i-1}}{h^n_i h^n_{i-1}}\cdot \frac{u^{n+1}_i-u^n_i}{k^n} = \nu D^2u^n_i.
\]
Using the mesh equation $\sigma^n_i = k^n u^n_i$, we obtain the explicit scheme
\[
(1+k^n Du^n_i)(1+k^n Du^n_{i-1})\bigg(\frac{u^{n+1}_i-u^n_i}{k^n}\bigg) = \nu D^2u^n_i
\]
with $t^n_{i+1}=t^n_i$ and $\sigma^n_i = k^n u^n_i$.

\item[Exercise \ref{difference moving frame exercise} part (\ref{difference moving frame exercise part a}):]  A compatible discrete cross-section is given by
\[
\mK = \bigg\{t^n = 0,\; x^n_i=0,\; u^n_i = 0,\; \frac{u^n_{i+1}}{x^n_{i+1}} + \frac{u^n_{i-1}}{x^n_{i-1}}=0,\; u^{n+1}_i = t^{n+1}\bigg\}.
\]
The corresponding discrete moving frame is
\begin{gather*}
\epsilon_1 = -x^n_i,\qquad \epsilon_2 = -t^n,\qquad \epsilon_3 = -u^n_i,\\
e^{\epsilon_4} = [(1+k^n\Delta_x u^n_i)(\Delta_t u^n_i + u^{n+1}_i \Delta_x u^n_i)]^{1/3},\qquad
\epsilon_5=-\Delta_x u^n_i,
\end{gather*}
where
$$
\Delta_x u^n_i =\frac{Du^n_i+Du^n_{i-1}}{2}\qquad\text{and}\qquad \Delta_tu^n_i = \frac{u^{n+1}_i-u^n_i}{k^n} - \frac{\sigma^n_i}{k^n}\cdot\Delta_xu^n_i.
$$

\item[Exercise \ref{difference moving frame exercise} part (\ref{difference moving frame exercise part b}):] 
The invariantization yields the finite difference invariant
\[
\iota(D^2u^n_i) = \cfrac{D^2u^n_i}{(1+k^n\Delta_x u^n_i)(\Delta_t u^n_i + u^{n+1}_i\Delta_x u^n_i)}.
\]

\item[Exercise \ref{difference moving frame exercise} part (\ref{difference moving frame exercise part c}):] 
Invariantizing $\Delta_t u^n_i + u^n_i \Delta_x u^n_i = \nu D^2u^n_i$ we obtain the invariant scheme
\[
(1+k^n\Delta_x u^n_i)(\Delta_t u^n_i + u^{n+1}_i\Delta_x u^n_i) = \nu D^2u^n_i.
\]
\end{description}
\end{appendices}

\end{document}